\documentclass[a4paper]{article}
\usepackage[utf8]{inputenc}
\usepackage[english]{babel} 
\usepackage{amsmath,amsfonts,amsthm,amssymb,xfrac} 
\usepackage{calrsfs}
\usepackage{arydshln}
\usepackage{graphicx}
\usepackage{esint}
\usepackage{booktabs}
\usepackage{chngpage}
\usepackage{todonotes}
\usepackage{graphicx}
\usepackage{enumerate}
\usepackage[margin=2cm]{geometry}
\usepackage{epstopdf}
\usepackage{listings}
\usepackage{lscape}
\usepackage{float, tikz, tikz-3dplot, tikz-cd, framed}
\usetikzlibrary{calc}
\usepackage[colorlinks]{hyperref}
\usepackage{pdfpages}
\usepackage{fancyhdr} 
\usepackage{authblk}

\numberwithin{equation}{section}%

\usepackage{setspace}
\hypersetup{colorlinks={true},linkcolor={black}}
\DeclareMathAlphabet{\pazocal}{OMS}{zplm}{m}{n}
\newtheorem{theorem}{Theorem}

\numberwithin{theorem}{section}

\newcommand{\COR}[1]{\textcolor{black}{#1}}

\makeatletter
 {\par\unskip\endMakeFramed}
\makeatother

\newcommand*{\eea}{\end{array}}

\newcommand*{\bme}{\begin{multiequations}}
\newcommand*{\eme}{\end{multiequations}}

\newcommand{\E}{\mathbb{E}(}

\providecommand\bcdot{\boldsymbol{\cdot}}

\renewcommand*{\Omega}{\varOmega}
\renewcommand*{\Sigma}{\varSigma}

\newcommand{\R}{\mathbb{R}}

\newcommand{\N}{\mathbb{N}}

\def\squarebox#1{\hbox to #1{\hfill\vbox to #1{\vfill}}}
\newcommand{\reel}{\mathbb R}

\newcommand{\w}{\boldsymbol{v}}
\newcommand{\bu}{\boldsymbol{u}}
\newcommand{\bv}{\boldsymbol{v}}

\newcommand{\bsigma}{\boldsymbol{\sigma}}

\newcommand{\xx}{\boldsymbol{x}}

\newcommand{\XX}{\boldsymbol{X}}

\newcommand{\nab}{\boldsymbol{\nabla}}
\newcommand{\nabh}{\nab_{\!\!\hor}}

\newcommand{\transp}{^{\scriptscriptstyle T}}

\newcommand{\Exp}{\mathbb{E}}

\newcommand{\dif}{{\mathrm{d}}}
\newcommand{\ba}{\boldsymbol{a}}
\newcommand{\mbs}[1]{\ensuremath{\boldsymbol{#1}}}

\renewcommand*{\vec}{\boldsymbol} 
\renewcommand*{\Omega}{\varOmega} 
\renewcommand*{\vec}{\boldsymbol} 
\newcommand{\bdot}{\boldsymbol{\cdot}} 
\newcommand{\sdbt}{\vec{\sigma} \df \vec{B}_t}  
\newcommand{\sodbt}{\vec{\sigma\circ} \df \vec{B}_t}
\newcommand{\sodbs}{\vec{\sigma\circ} \df \vec{B}_s}
  
\renewcommand*{\div}{\boldsymbol{\nabla\cdot}} 

\newcommand{\dt}{\df t}
\newcommand{\grad}{\boldsymbol{\nabla}} 
\newcommand{\tp}{^{\scriptscriptstyle T}} 

\renewcommand*{\div}{\boldsymbol{\nabla\cdot\,}} 
\newcommand{\adv}{\boldsymbol{\cdot\nabla}} 
\newcommand{\dpt}{\df h_t^{\sub\sigma}}

\newcommand{\hor}{{\scriptscriptstyle H}}

\newcommand{\bs}{\boldsymbol} 
\newcommand{\mb}{\mathbb} 
\newcommand{\mc}{\mathcal} 
\newcommand{\mr}{\mathrm}

\newcommand{\sub}{\scriptscriptstyle} 
\newcommand{\beq}{\begin{equation}}
\newcommand{\eeq}{\end{equation}}
\newcommand{\beqs}{\begin{subequations}}
\newcommand{\eeqs}{\end{subequations}}
\newcommand{\benum}{\begin{enumerate}[label=(\roman*)]}
\newcommand{\eenum}{\end{enumerate}}

\newcommand{\dom}{\mc{D}} 

\newcommand{\df}{\mr{d}} 
\newcommand{\sto}{\mr{\mb{D}}} 
\newcommand{\alf}{\frac{1}{2}} 
\newcommand{\salf}{\sfrac{1}{2}}
\newcommand{\kerl}{\breve{\sigma}} 
 
\newcommand{\bdiv}{\grad \bdot} 

\newcommand{\Pb}{\mb{P}} 
\newcommand{\filt}{\mc{F}} 
\newcommand{\ds}{\df s} 
\newcommand{\bsig}{\bs{\sigma}} 
\newcommand{\bkerl}{\breve{\bs{\sigma}}} 
\newcommand{\bxi}{\bs{\xi}} 

\newcommand{\B}{\bs{B}} 
\newcommand{\noi}{\bsig \df \B_t}
\newcommand{\Snoi}{\bsig \circ \df \B_t}

\newcommand{\x}{\bs{x}}
\newcommand{\y}{\bs{y}}

\newcommand{\X}{\bs{X}}
\newcommand{\bk}{\bs{k}}

\newcommand{\us}{\bv^*}

\newcommand{\dx}{\df \x} 
\newcommand{\dy}{\df \y}

\newcommand{\divh}{\nabh\bdot}

\newcommand{\modif}[1]{\color{black} #1 \color{black}}

\begin{document}

\title{Derivation of stochastic models for coastal waves}

\author[1]{Arnaud Debussche}
\author[2]{Étienne Mémin}
\author[,2]{Antoine Moneyron \thanks{Corresponding author: antoine.moneyron@inria.fr}}

\affil[1]{Univ Rennes, CNRS, IRMAR UMR 6625, F35000, Rennes, France.}
\affil[2]{Univ Rennes, INRIA, IRMAR UMR 6625, F35000, Rennes, France.}

\date{\today}

\onehalfspacing

\maketitle
\thispagestyle{empty}

\begin{center}
    \rule{0.5\linewidth}{1pt}
\end{center}

\begin{abstract}

In this paper, we consider a stochastic nonlinear formulation of classical coastal waves models under location uncertainty (LU). In the formal setting investigated here, stochastic versions of the Serre-Green-Naghdi, Boussinesq and classical shallow water wave models are obtained through an asymptotic expansion, which is similar to the one operated in the deterministic setting. However, modified advection terms emerge, together with advection noise terms. These terms are well-known features arising from the LU formalism, based on momentum conservation principle.

\end{abstract}

\singlespacing
\tableofcontents
\onehalfspacing

\section{Introduction}
The ocean surface waves constitute an essential component of ocean dynamics, as they are directly related to a strong energy exchange with the underlying current. However, the mean current and the waves follow different dynamics, operating at different times and spatial scales. Moreover, they are based on significantly different physical assumptions: surface waves rely on potential -- or irrotational -- flow,   whereas the turbulent dynamics of the current is expressed through a non-zero vorticity. As a result, the wave-current coupling is very complex to model numerically. For this purpose, designing simplified stochastic models of surface waves would be advantageous to account for both dynamics. This can be achieved by adding as a specific noise term in the stochastic representation of ocean circulation \cite{LDLM_2023, TML_2023}.

To derive such stochastic wave models, \modif{ways to adapt the classical deterministic derivation to the stochastic setting have been widely investigated.} It was shown for example that deep-ocean stochastic long waves can arise from a linearised stochastic shallow water system \cite{MLLTC_2023} or from a stochastic Hamiltonian formulation \cite{Dinvay-Memin-PRSA-22}. Water waves travelling from deep water areas to shallower regions -- where the water depth is much less than their wavelength -- undergo significant alterations. Swift variations in height, velocity, and direction lead to substantial modifications of the free water surface profiles. Initially resembling almost perfect sine waves, these profiles evolve into an asymmetric shape. Mathematically, coastal waves are described by different dynamical models corresponding to various approximations of the irrotational Euler equations -- which \modif{are averaged along the water column ultimately.}

\modif{More specifically, we use the location uncertainty principle developed in \cite{Memin14}, which is based on the addition of a noise term on the Lagrangian formulation of the displacement. In fact, this approach is a critical aspect in ocean stochastic modelling and has been widely discussed: an SDE based stochastic generalisation of the deterministic Lagrangian expression of the flow was proposed in \cite{Restrepo_2007,RRmWB_2011} for instance, from which the authors derive an Eulerian expression. Such idea also exists in turbulence modelling: for example McWilliam and Berloff proposed stochastic parametrization built upon Langevin models of turbulence in \cite{BmW_2002}, devised in the wake of Kraichnan's seminal work \cite{Kraichnan59}.}

\modif{In this paper,} we aim to derive a location uncertainty interpretation of \modif{nonlinear coastal waves models,} based on the modelling principle introduced in \cite{Memin14}. These models are naturally accompanied by some nonlinear variabilities associated with numerical or modelling uncertainties. More fundamentally, this stochastic formalism provides a way to incorporate the effect of the non-resolved vertical turbulent component. In this work, we focus on the family of models associated with the Serre-Green-Naghdi equations \cite{Serre_1953, GN_1976}, which allow capturing the non-hydrostatic phenomena related to vertical acceleration. In particular, we aim to assess the behavior of the associated numerical simulations. \modif{Our main findings are that the LU interpretations of the Saint-Venant, Boussinesq and Serre-Green-Naghdi wave models allow to break the symmetry of the wave and introduce variability compared to their deterministic counterparts. Nevertheless, they do enjoy the same conservation principles with slight conditions on the noise structure. Hence, the LU setting would allow to explore the influence of a ``conservative randomness'' on the classical wave model dynamics. Plus, we provide some graphical observations of this influence using numerical simulations.}

The article is organised as follows: after a brief summary on the location uncertainty principle (LU), we derive a stochastic representation of the Serre-Green-Nagdhi shallow water waves model following the same physical approximations used in the deterministic setting \cite{Berthelemy-04}. Then, we show how simpler models can be obtained from further approximations, \modif{and discuss the conservation of usual quantities such as mass, momentum and energy.} The last section describes and assesses some numerical results associated with different models.

\section{Stochastic flow and transport}

Let $\dom \subset \R^d$ (with $d = 2$ or $3$) denote a bounded spatial domain. The LU principle relies on a stochastic formulation of the Lagrangian trajectory $(X_t)$ of the form
\begin{equation}
    dX_t(x) = v(X_t(x), t) dt + (\sigma_t dB_t)(x), \quad X_0(x) = x \in \dom. \label{LU-separation}
\end{equation}
This stochastic differential equation (SDE) decomposes the flow in terms of a smooth-
in-time velocity component $v$ that is both spatially and temporally correlated, and a fast unresolved flow component $\sigma dB_t$ (called noise term in the following). This latter component must be understood in the It\={o} sense, and is uncorrelated in time yet correlated in space. Importantly, since the It\={o} integral is of null expectation (as a martingale), the relation \eqref{LU-separation} decomposes the flow unequivocally as a large-scale displacement component and a null mean fluctuation, the mean displacement corresponding to the large-scale component expectation.

Let us now provide a precise definition of this random fluctuation component. The noise term takes values in the Hilbert space $H := (L^2 (\dom))^d$ and is defined from the Wiener process (also called cylindrical Brownian motion) \cite{DaPrato} on a stochastic basis $(\Omega, \filt_t, (\filt_t)_{t\in[0,T]}, \Pb)$. By definition, this stochastic basis is composed of a probability space $(\Omega, \Pb)$ with filtration $(\filt_t)_{t\in[0,T]}$ -- {\it i.e.} a non decreasing family of sigma algebras evolving in time. The noise term involves independent Wiener process components $(B_t^i)_{i=1,\ldots,d}$ defined on an orthogonal basis $({\mbs e_n})_{n\in \N}$ of $H$ as
\modif{
\beq
(B_t^i)(x) = \sum_{n\in\N} \beta_t^n {\mbs e}^i_n(x),
\eeq}
 where \modif{$(\beta_t^i)_{i\in \N}$} is a sequence of independent one dimensional standard Brownian motions. The spatial structure of the unresolved flow component is modelled by the correlation operator, $\bsig_t$, defined as an integral operator on $H$. Let a matrix kernel $\bkerl = (\kerl_{ij})_{i,j=1,\ldots,d}$ that is bounded in space and time, then define
\beq\label{eq:sigma}
    \bsig_t\, \bs{f} (\xx) = \int_{\dom} \bkerl (\x, \y, t) \bs{f} (\y)\, \dy,\ \quad \bs{f} \in H,\ \quad \x \in \dom .
\eeq
Here, we have assumed that this operator is deterministic. However, it is important to note that, if required, it could be defined  as a random operator. \COR{An example of such a random correlation operator can be found in \cite{Li-et-al-2022} where it is defined from the dynamic mode decomposition technique (DMD) \cite{Schmid10} and a Girsanov transform. Girsanov transform enables in particular to introduce a drift term associated to a non-centred noise, through a change of probability measure. This procedure has proven particularly useful for noise calibration in data assimilation \cite{Dufee-QJRMS-22}.} The composition of $\bsig_t  [\bullet]$ and its adjoint operator $\bsig^*_t [\bullet]$ defines a compact self-adjoint positive operator, of which eigenfunctions and eigenvalues are denoted $\bxi_n (\bdot,t)$ and $\lambda_n (t)$ respectively. These eigenvalues fulfil $\sum_{n \in \mb{N}} \lambda_n (t) < + \infty$ and decrease toward zero at infinity, and the eigenfunctions form an orthonormal basis of $H$. As such, the noise can be equivalently defined on this basis as the following spectral expression,
\beq\label{eq:KL}
    \bsig_t\, \df \B_t (\x)= \sum_{n \in \mb{N}} \lambda_n^{1/2} (t) \bxi_n (\x, t)\, \modif{\df \beta_t^n}.
\eeq
Consequently, the noise component is a $H$-valued Gaussian process of null mean and with bounded variance \modif{-- that is $\Exp_{\sub \Pb} [\int_0^t \bsig_s\, \df \B_s] = \bs{0}, \forall t>0$ and $\Exp_{\sub \Pb} \big[ \| \int_0^t \bsig_s\, \df \B_s \|_{\sub H}^2 \big] < +\infty, \forall t>0$ -- under the probability measure $\mb{P}$.} Moreover, the auto-covariance at point $\x \in \dom$ of the unresolved flow component at each instant $t \in [0,T]$ is given by the matrix kernel of the composite operator $\bsig \bsig^*$, and denoted by $\ba (\x, t)$, namely
\beq\label{eq:bracket}
    \ba (\x, t) := \int_{\dom} \bkerl (\x, \y, t) \bkerl\tp (\y, \x, t)\, \dy 
    = \sum_{n \in \mb{N}} \lambda_n (t) \big( \bxi_n \bxi_n\tp \big) (\x, t) .
\eeq
The process $\int_0^t \ba (\x, s)\, \ds$ corresponds to the quadratic variation of $\int_0^t \bsig_s \, \df \B_s (\x)$ \cite{Bauer-et-al-JPO-20}.

The stochastic integral defining the noise  could have been defined in terms of a Stratonovich integral instead of an It\={o} integral \cite{Bauer-et-al-JPO-20}.  For  a deterministic  correlation operator, $\ba$ boils down to the local variance of the noise, due  to the It\={o} isometry \cite{DaPrato}. Hence, we referred to $\ba$ as the variance tensor, although for a random correlation operator this denomination is misleading since it is a random process. Physically, the symmetric non-negative tensor $\ba$ represents the friction coefficients of the unresolved fluid motions, and is physically homogeneous to a viscosity with unit m$^2$/s. 

Consider an extensive random tracer $\Theta$ (e.g. temperature, salinity or buoyancy) transported by the stochastic flow, fulfilling the following conservation property along the trajectories: $\Theta (\X_{t + \delta t}, t + \delta t) = \Theta (\X_t, t)$ with $\delta t$ an infinitesimal time variation. Thus, the evolution law of $\Theta$ is given by the following stochastic partial differential equation (SPDE),
\beq\label{eq:Dt-LU}
    \sto_t \Theta = \df_t \Theta + (\us\, \dt + \noi) \adv \Theta - \alf \bdiv (\ba \grad \Theta)\, \dt = 0 ,
\eeq
where $\sto_t$ is introduced as a stochastic transport operator and $\df_t \Theta (\x) := \Theta (\x, t+\delta t) - \Theta (\x, t)$ stands for the forward time-increment of $\Theta$ at a fixed point $\x \in \dom$. This operator $\sto_t$ was derived in \cite{Memin14} using the generalized It\={o} formula (also called It\={o}-Wentzell formula in the literature) \cite{Kunita} and plays the role of a transport operator in a stochastic setting. Remarkably, it  encompasses physically meaningful terms \cite{Bauer-et-al-JPO-20, Resseguier-GAFD-I-17}: the two first terms correspond to the classical terms of the material derivative, while the third term describes the tracer advection by the unresolved flow component. As shown in \cite{Bauer-et-al-JPO-20, Resseguier-GAFD-I-17}, the resulting (non Gaussian) multiplicative noise $\noi \adv \Theta$ continuously backscatters random energy to the system through the quadratic variation $\salf\, (\grad \Theta) \bdot (\ba \grad \Theta)$ of the random tracer. The last term in \eqref{eq:Dt-LU} represents the tracer diffusion due to the mixing of the unresolved scales. The energy loss by the diffusion term is exactly balanced with the energy brought by the noise. This pathwise balance (i.e for any realisation) leads to tracer energy conservation and highlights the parallel between the classical material derivative and the stochastic transport operator. This will be precised below.

Under specific noise and/or variance tensor definitions, the resulting diffusion \cite{Kadri-CF-17,Memin14} can be connected to the additional eddy viscosity term introduced in many large-scale circulation models \cite{Smagorinsky63, Redi-82}. As an additional feature of interest, this evolution law introduces an \emph{effective} advection velocity $\us$ in \eqref{eq:Dt-LU} defined for an incompressible noise term as
\beq\label{eq:def-ISD}
    \us = \bv -  \alf \bdiv \ba .
\eeq
This statistical eddy-induced velocity drift captures the action of inhomogeneity of the random field on the transported tracer and the possible divergence of the unresolved flow component. It is shown in \cite{Bauer-et-al-JPO-20} that the \emph{turbophoresis} term, $\bv_s = \salf \bdiv \ba$, can be interpreted as a generalization of the Stokes drift associated to surface wave current and that it plays a key role in the triggering of secondary circulations such as  the Langmuir circulation \cite{Craik-Leibovich-76,Leibovich80}. Consequently, this velocity is termed It\={o}-Stokes drift (ISD) in \cite{Bauer-et-al-JPO-20}. Notice that, in the modified advection \eqref{eq:def-ISD}, the ISD cancels out for homogeneous random fields (since then the variance tensor is constant over space).

Many useful properties of the stochastic transport operator $\sto_t$ have been explored by \cite{Resseguier-GAFD-I-17, Resseguier2020arcme}. In particular, if a random tracer is transported by the incompressible stochastic flow under suitable boundary conditions, then the pathwise $p$-th moment $(p \geq 1)$ of the tracer is materially and integrally invariant, namely
\beq\label{eq:p-moment}
    \sto_t \left( \frac{1}{p} \Theta^p \right) = 0,\ \quad \df_t \left(\int_{\dom} \frac{1}{p} \Theta^p\, \dx\right) = 0 .
\eeq

It is worth noting that the transport equation \eqref{eq:Dt-LU} can be written in terms of Stratonovich integral \cite{Bauer-et-al-JPO-20}, as
\beq\label{eq:Dt-LU-Strato}
    D_t \Theta = \df_t \Theta + (\us\, \dt + \Snoi) \adv \Theta  = 0.
\eeq
The Stratonovich integral has the advantage of fulfilling a ``standard" chain rule, so that the notation $D_t$ for the transport operator similar to the material derivative. It\={o} calculus introduces second order terms -- such as the diffusion term -- which becomes implicit in the Stratonovich integral. Because of this,  Stratonovich noise is not anymore a martingale and is not of null expectation. Moreover, it is possible to move safely from one integral to the other under some regularity assumptions. In the following we will use the Stratonovich notation.

\section{Non linear Shallow Water equations under location uncertainty}

In this section, we derive a stochastic representation of ocean waves in the near-shore zone, focusing on the derivation of models going from the Shallow water model to the Serre-Green-Nagdi model. Our stochastic derivation is similar to the scaling procedure described in \cite{Barthelemy-04}. Our derivation starts from the general 3D Euler equations in the LU setting, which read
\begin{align}
    &\dif_t \w + (\w - \alf \div \ba)\bcdot\nab \w \dif t + (\sodbt \bcdot\nab) \w  = -\frac{1}{\rho}\nab(P \dif t + \dpt) - \mathbf{g}z\dif t,\\
    &\nab\bcdot(\w - \alf \div\ba) =0 \label{Cont-eq},
\end{align}
denoting $\bv= (\bu,w)$ the three-dimensional velocity decomposed in terms of horizontal, $\bu$, and vertical, $w$, components. The pressure is denoted $P$, while $\rho$ and $\mbs{g}$ stand for the density and the gravitational vector directed along the vertical direction, respectively.
 
\begin{figure}[h] 
    \centering \includegraphics[width=16cm]{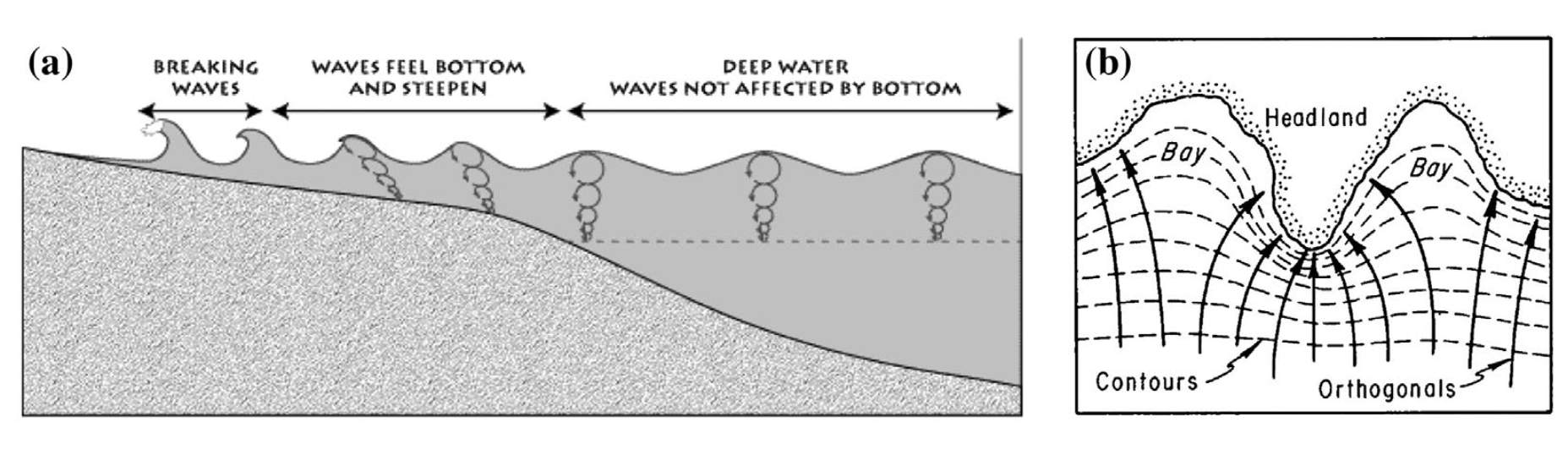}
    \caption{Illustration of shoaling processes, borrowed from \cite{KG_2016}. (a) Wave shoaling diagram. (b) Wave refraction diagram. Source: USACE coastal engineering Manual. As the waves approach the coasts, both the typical wavelength and the water depth decrease, the former much faster than the latter. Hence, $\beta$ decreases and $\epsilon$ increases as the waves approach the coast. Thus, shoaling processes are characterised by moderately large values of the parameters $\epsilon$ and $\beta$ -- that is, away enough from the coast to avoid breaking waves, yet near enough for the bottom topography to influence the wave dynamics.} \label{fig:shoaling}
\end{figure}

Shallow water conditions are characterised by a water depth $h$ being much smaller than the wave length scale \modif{$L \sim 1/|\bk|$}, where $\bk$ denotes the wave number. This condition is often expressed through the quantity $\beta = \bk h_0 \ll 1$, which informs about the predominance of dispersion. In linear theory, the amplitude $A$ of the wave is small and tends to zero when the characteristic scale tends to infinity. Another quantity is usually introduced to measure the non-linear effects: $\epsilon = A/h_0$. Shoaling processes start for $\beta \leq 1$ (i.e. when wavelength and depth have the same order) and ends when the waves break at $\epsilon \geq 1$ (an illustration is provided in figure \ref{fig:shoaling}, which was borrowed from \cite{KG_2016}). Thus, shoaling requires asymptotic models with short wavelength and high wave amplitude. These adimensional numbers are used to build approximated solutions of the shallow water waves system, ranging from weakly nonlinear Boussinesq models (small amplitude regime, $\beta^2 << 1$ and $\epsilon =O(\beta^2)$) to the nonlinear Serre-Green-Naghdi system (large amplitude regime). \modif{As it will be shown} in section \ref{Derivation}, within the modeling under uncertainty setting the Serre-Green-Naghdi model reads
\begin{subequations}
\label{S-SGN-I}
\begin{align}
    &\overline{D}^\hor_t \eta = - h \bigl(\divh (\overline{\bu} -\alf \Upsilon\epsilon \overline{\bu}_s)\dif t + \Upsilon \divh\overline{\sodbt^\hor}\bigr),\label{S-SGN-I-a}\\
\label{S-SGN-I-b}
\begin{split}
    &\overline{D}^\hor_t \overline{\bu}  + \nabh \eta\dif t - \frac{1}{h} \epsilon\beta^2 \nabh \bigl(\frac{h^3}{3} (\dif\overline{G}) \bigr) = {\cal O}(\beta^4,\epsilon \beta^4),
\end{split}
\end{align}
\end{subequations}
with 
\beq
    \dif \overline{G}(\xx,t)= \overline{D}^\hor_t (\divh \overline{\bu}) -  \epsilon\divh\bigl(( \overline{\bu}- \alf \Upsilon \epsilon \overline{\bu}_s)\dif t +  \Upsilon^{\salf} \overline{\sodbt^\hor}\bigr)  \divh \overline{\bu}.
\eeq

\begin{figure}[h] 
    \centering \includegraphics[width=14cm]{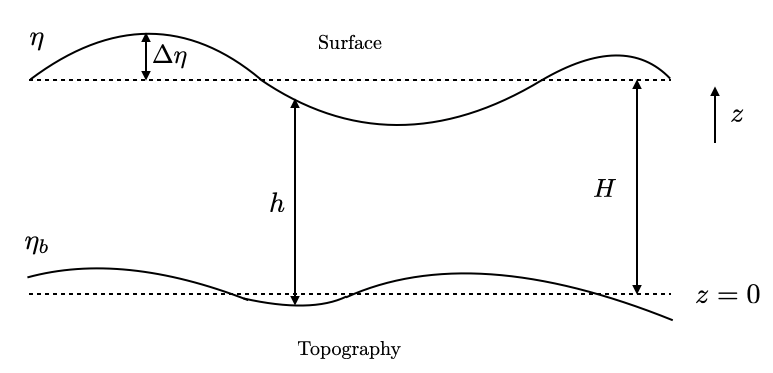}
    \caption{Illustration of the bottom topography $\eta_b$, water depth $h$ and surface deformation $\eta$.} \label{fig:summary}
\end{figure}

The notation ${D}^{\hor}_t f $ stands for the Stratonovich transport operator with respect to vertically averaged stochastic flow components,
\begin{equation}\label{DHf}
    \overline{D}^{\hor}_t f = \dif_t f  +\epsilon\, \nabh f \bcdot \overline{\bu}^*\, \dif t + \Upsilon^{\salf}\epsilon\,\nabh f \bcdot \overline{\sodbt^\hor},
\end{equation}
 with the depth averaged horizontal velocity
\begin{equation}
    \overline{{\bu}}(x,t) = \frac{1}{h}  \int_0^{h(x,t)}  \!\!\!\!\!\! \!\!\!\!\!\! \bu(z) \dif z.
\end{equation}
The last left-hand side term of \eqref{S-SGN-I-b} can be understood as a pressure term associated to the vertical velocity component acceleration corrected by compressibility effects. This term is quite intuitively a function of the average horizontal velocity divergence. Overall, the system constitutes a stochastic version of the Serre-Green-Naghdi equations \cite{Green-Naghdi-76, Serre-53}: compared to the original deterministic model, it includes an additional transport noise term and the contribution of the It\={o}-Stokes drift induced by the inhomogeneity of the small-scale fluctuations. A graphical illustration of the bottom topography, water depth and surface deformation is provided on figure \ref{fig:summary}.

\subsection{Deriving the stochastic Serre-Green-Naghdi system}
\label{Derivation}

In this section, we derive a stochastic representation of the Serre-Green-Naghdi model.

\subsubsection{Evolution equations}

\paragraph{Scaling relations}
First, we proceed to the following adimentionalisation to define the asymptoptic models. Consider $x = L \tilde x$, $z= h_0 \tilde z$, $\w = \epsilon\, C_0 \tilde \w$ with $C_0= \sqrt{gh_0}$ the characteristic long wave velocity. Also, time is scaled as $T\sim L/C_0$. The notation $\tilde\bullet$ stands for adimensioned variables. We will further assume the following incompressible assumption for the noise and the ISD,
\begin{equation}
    \div\sodbt =0 \text{ and } \div \bv_s =0.
\end{equation}
These assumptions on the noise lead to the global energy conservation -- expressed as the sum of the potential and kinetic energies -- of the shallow water system \cite{Brecht-et-al-2021}, as well as the energy conservation of transported scalar \cite{Resseguier-GAFD-I-17, Bauer2020jpo}. The incompressibility condition on the flow boils down to the classical divergence-free condition
\begin{equation}
    \div \w =0.
\end{equation}
This incompressibility condition leads to the scaling $w = \epsilon \beta C_0 \tilde{w}$ for the vertical velocity. In addition, we scale the noise term and the variance tensor $\ba$ by  $\Upsilon^{\salf}$ and $\Upsilon$, respectively. The incompressibility condition on the noise term writes
\beq
    \noi^h = \Upsilon ^{\salf}  L \widetilde{\noi^h},
\eeq
with a similar scaling on the horizontal noise displacement. This provides a scaling on the vertical noise component as 
\beq
    {\sodbt^z} = \Upsilon^{\salf} \epsilon \beta L \,\widetilde{\sodbt^z}.
\eeq
The variance tensor has the  dimension of a viscosity in $m^2/s$. For a deterministic correlation tensor, $\bsigma$, it is related to the variance of the flow fluctuation, $\bv' = (\bu',w')^{\transp}$ through 
\begin{equation}
    \ba = \Exp (\bv'{\bv'}\transp)\tau,
\end{equation}
where $\tau$ is a decorrelation time. Now denote $\ba^{\hor \hor}$, $\ba^{\hor z}$ and $a^{zz}$ the horizontal (2D) matrix variance tensor, the horizontal-vertical cross vector, and the vertical variance, respectively. Thus, we get the following scaling relations, 
\begin{subequations}\label{scaling-a-sys}
\begin{align}
    &\ba^{\hor \hor} = \Upsilon \epsilon^2 C_0 L \,\tilde{\ba}^{\hor \hor}\\
    &\ba^{\hor z} = \Upsilon \epsilon^2\beta C_0 L \,\tilde{\ba}^{\hor z} \label{scalling-a-h}\\
    &a^{zz} = \Upsilon \epsilon^2\beta^2 C_0 L \,\tilde{a}^{zz}.
\end{align}
\end{subequations}
Note that these relations could have been inferred directly from the noise scalings and relation 
\beq
\ba \dif t = \Exp (\noi (\noi)\transp).
\eeq
Additionally, the ISD component scales as
\begin{equation}\label{ISD-scaling}
    \salf\div\ba = (\bu_s, w_s) = \Upsilon \epsilon^2 C_0 \bigl( \nabla_{\tilde{H}} \bcdot \tilde \ba^{\hor \hor} +\partial_{\tilde{z}} \tilde \ba^{\hor z}, \beta \nabla_{\tilde{H}} \bcdot \tilde \ba^{\hor z} + \beta \partial_{\tilde{z}} \tilde a^{zz}\bigr),
\end{equation}
where $\nabh$ refers to the gradient with respect to the horizontal coordinates. 
As in the classical setting, we will assume that the large-scale component of the flow is irrotational and thus
\begin{equation}
    \partial_{\tilde{y}} \tilde u = \partial_{\tilde{x}} \tilde v \text{ and } \partial_{\tilde{z}} \tilde \bu = \beta^2 \nabla_{\tilde{H}} \bcdot \tilde w,
\end{equation}
while the adimensional continuity equation is 
\beq
    \nabla_{\tilde{H}} \tilde \bu = -\partial_{\tilde{z}} \tilde w.
\eeq

In the long wave assumption, the pressure is scaled by the static pressure as follows
\modif{
\beq
    P\sim \rho g h_0 \text{ and thus } P=\rho g h_0\tilde P.
\eeq}
The same scaling is assumed for the turbulent pressure, with the additional noise variance scaling $\Upsilon^{\salf}$,
\modif{
\beq
    \dpt = \Upsilon^{\salf}\rho g h_0 \tilde{\dpt}.
\eeq}
\modif{For simplicity,} from now on we will drop the tilde accentuation for the adimensional variables in the following.

\paragraph{Incompressibility condition}
Integrating over the z axis, the incompressibility condition \eqref{Cont-eq} gives 
\modif{\beq
    \int_0^{h(\boldsymbol{x},t)} \divh ( \bu -\alf \Upsilon\epsilon\, \bu_s)\dif z = - \biggl( w\bigl(h(\boldsymbol{x})\bigr) -\alf \Upsilon\epsilon\, w_s\bigl(h(\boldsymbol{x})\bigr) \biggr).
\eeq}
Using Leibniz formula and introducing the depth averaged horizontal effective velocity $\overline{{\bu}^*}$,
\begin{equation}
    \overline{{\bu}^*}(\boldsymbol{x},t) = \frac{1}{h}  \int_0^{h(\boldsymbol{x},t)}  \!\!\!\!\!\! \!\!\!\!\!\! \bu^*(z) \dif z,
\end{equation}
we have
\beq \label{averaged-continuity}
    \divh \int_0^{h(x,t)} \!\!\!\!\!\! \!\!\!\!\!\! \bu^* \dif z - \bu^*(h)\bcdot \nabh h= \divh (\overline{\bu^*}\, h) -  \bu^*(h)\bcdot \nabh h = \COR{-}  w^* (h).
\eeq
In the same way,  we have for the noise term
\beq \label{averaged-continuity-noise1}
    \divh (h \,\overline{\noi^\hor}) - \divh \bigl(h \noi^\hor(h)\bigr) = -\noi^z(h).
\eeq

\paragraph{Adimensionned momentum equations}
Gathering the previous scaled relations, we obtain the following equations for the horizontal and vertical velocity,
\begin{multline}
    \epsilon \dif_t \bu +  \COR{\epsilon^2 \divh \bigl((\bu -\alf \Upsilon \epsilon \bu_s) \bu\bigr)\dif t +  \epsilon^2 \partial_z\bigl((
    w-\alf \Upsilon \epsilon w_s) \bu\bigr)\dif t
    + \Upsilon^{\salf}\epsilon^2\div \bigl({\sodbt^\hor}\,  \bu\bigr) + \Upsilon^{\salf}\epsilon^2  \partial_z\bigl({\sodbt^z} \,\bu \bigr)}\\
   = \COR{-}\nabh P\dif t \COR{-} \Upsilon^{\salf} \partial_x \dpt, 
\end{multline}
and
\begin{multline}
    \epsilon \beta^2  \dif_t w +  \epsilon^2 \beta^2 \bigl( (\bu - \alf \Upsilon \epsilon \bu_s) \bcdot \nabh w\bigr) + \epsilon^2 \beta^2 \bigl( (w - \alf \Upsilon \epsilon w_s)\partial_z w \bigr)\dif t+
    \epsilon^2 \beta^2 \Upsilon^{\salf}(\sodbt^\hor\bcdot\nabh w) +  \epsilon^2 \beta^2 \Upsilon^{\salf}(\sodbt^z\partial_z w)\\
     = -(\partial_z P + 1)\dif t - \COR{\Upsilon^{\salf}} \partial_z \dpt.
\end{multline}
The vertical momentum equation can be written in a more compact form, that is
\begin{equation}\label{vert-acc}
    \epsilon \beta^2 D_t w  = -\partial_z \dif P -1\dif t,
\end{equation}
where $P$ represents the total pressure, sum of the finite variation  and  martingale pressures. In this formula, $1$ stands for the rescaled gravity term.

\subsubsection{Boundary conditions}

\paragraph{Bottom boundary conditions} For the boundary conditions we will assume that $w(\boldsymbol{x},z,t) = \sodbt^z = 0$ on the bottom $z=0$. Due to this last condition we have also $\ba^{Hz}= a^{zz} =0$ on the bottom. The constraint of the vertical ISD from \eqref{ISD-scaling} implies that $\partial_z a^{zz}(\boldsymbol{x},0,t)=0$ on the bottom. As a consequence, 
\[
    \div\ba (\boldsymbol{x},0,t) = \Upsilon \epsilon^2 C_0 \bigl( \divh \ba^{\hor \hor}(\boldsymbol{x},0,t) +\partial_z  \ba^{\hor z}(\boldsymbol{x},0,t), 0 \bigr).
\]
The divergence-free constraint of the ISD leads to following condition on the bottom, 
\begin{equation}
    \divh \divh \ba^{\hor \hor}(\boldsymbol{x},0,t) =0.
\end{equation}
We will also consider that the noise and the ISD have the same characteristics as the large scale velocity near the bottom: due to the large-scale component being irrotational, the following hold in the vicinity of the bottom
\begin{equation}
    \partial_y \sodbt^\hor = \partial_x \sodbt^\hor, \quad   \beta^2 \nabh \sodbt^z =   \partial_z \sodbt^\hor  \text{ at } z= \delta z,
\end{equation}
and 
\begin{equation}
    \partial_x v_s = \partial_y u_s, \quad  \beta^2 \nabh w_s=   \partial_z \bu_s  \text{ at } z=\delta z.
\end{equation}

\paragraph{Free surface boundary conditions}
At the free surface $z= h(\boldsymbol{x},t)$, we denote by $\eta$ the surface elevation variation and consider the scaling
\[
    z= h(\boldsymbol{x},t) = h_0\bigl(1 + \epsilon \eta (\boldsymbol{x},t) \bigr),
\]
which leads to 
\begin{equation}
    \tilde z = \frac{1}{h_0}h(\boldsymbol{x},t) =\bigl(1+\epsilon \eta (\boldsymbol{x},t)\bigr).
\end{equation}
The evolution of the \emph{dimensional} surface elevation is given by 
\begin{equation}
    (w - w_s)\dif t  + \sodbt^z = D_t \eta,
\end{equation}
where an effective vertical velocity driving the dynamics of $\eta$ appears on the left-hand side. This velocity is composed of the vertical components of the velocity, the ISD, and an additive noise variable linked with rapid vertical oscillations. 
As usual in long wave approximations, we assume that the \COR{whole} pressure is constant at the interface, and for sake of simplicity we consider that this constant is null, that is
\begin{equation}
    \dif P_t(\boldsymbol{x} ,h(\boldsymbol{x} ,t) ,t) =0.
\end{equation}
The adimensional form of the free surface evolution (dropping the tilde accentuation) reads,
\modif{
\begin{equation}
    \bigl(w(h) - \alf\Upsilon \epsilon w_s(h)\bigr)\dif t + \Upsilon^{\salf}\, \sodbt^z(h) = \dif_t\eta +\epsilon \nabh\eta\bcdot (\bu(h) -\alf \Upsilon \epsilon \bu_s(h))\, \dif t + \Upsilon^{\salf}\epsilon\, \nabh \eta \bcdot \sodbt^\hor(h).
\end{equation}}
Considering the averaged continuity equations (\ref{averaged-continuity}) we obtain
\begin{equation}
\label{S-free-surf}
    \dif_t\eta +  \epsilon\, \nabh \eta \bcdot (\overline{\bu} -\alf \Upsilon\epsilon \overline{\bu}_s)\, \dif t + \Upsilon^{\salf}\epsilon\, \nabh \eta \bcdot \overline{\sodbt^\hor} = \\- h \bigl(\divh (\overline{\bu} -\alf \Upsilon\epsilon \overline{\bu}_s)\dif t + \Upsilon^{\salf}  \divh\overline{\sodbt^h}\bigr),
\end{equation}
which can be written more compactly as
\begin{equation}\label{dyn-surf}
    \overline{D}^{\hor}_t \eta = - h \bigl(\divh (\overline{\bu} -\alf \Upsilon\epsilon \overline{\bu}_s)\dif t + \Upsilon^{\salf}  \divh\overline{\sodbt^\hor}\bigr).
\end{equation}
\COR{Notice that we have introduced an operator $\overline{D}^{\hor}_t$, that involves only depth average velocity components.}
Moreover, remark that for homogeneous noise (i.e. with no statistical dependence on space location), the variance tensor is constant over space and the ISD cancels out. In such a case, the surface elevation dynamics simplifies in the more intuitive form, 
\begin{equation}\label{homo-surf-dyn}
    \dif_t \eta +\epsilon\, \nabh \eta \bcdot \overline{\bu}\, \dif t + \Upsilon^{\salf}\epsilon\,\nabh \eta \bcdot \overline{\sodbt^\hor} = \\- h \bigl(\divh \overline{u} \dif t + \Upsilon^{\salf}  \divh \overline{\sodbt^\hor}\bigr).
\end{equation}

Additionally, in the general case, the integrated ISD can be written in terms of horizontal quantities,  
\begin{align}
    h\overline{\bu}_s &=  \bigl( \int^{h(x,t)}_0 \divh \ba^{\hor \hor} + \ba^{\hor z}(h)\bigr)\noindent =  \divh (h \overline{\ba}^{\hor}) - \ba^{\hor \hor}(h)\nabh h + {\ba}^{\hor z}(h) \noindent \nonumber\\
    &=  \divh (h \overline{a}^{\hor})  - \ba^{\hor \hor}(h)\nabh h + \langle \sodbt^z(h),\sodbt^\hor(h) \rangle  \nonumber\\
    &=  \divh (h \overline{\ba}^{\hor})  - \ba^{\hor \hor}(h)\nabh h - \alf h \divh \ba^{\hor \hor}(h).
    \label{aver-ISD}
\end{align}

\subsubsection{The wave model}

\paragraph{Depth averaged horizontal momentum equation}

By taking the average of the horizontal momentum equation, and using Leibniz rules in time and space for stochastic processes (see appendix \ref{sec:Leibniz}), the continuity relation and the elevation evolution together with the velocity boundary conditions, we obtain
\begin{multline}
     \epsilon h\, \dif_t \overline{\bu}+ \epsilon^2 h\, \bigl(\overline{\bu}-\alf \Upsilon\epsilon\, \overline{\bu}_s)\bcdot\nabh\bigr)  \overline{\bu}\,\dif t + \Upsilon^{\salf}\epsilon^2 h \, (\overline{\sodbt^\hor} \bcdot\nabh) \overline{\bu}\\ 
     - \COR{h\epsilon^2} \divh  \int_0^{h(x,t)}\bigg( (\bu-\alf \Upsilon\epsilon \bu_s)\bu\dif t + \Upsilon^{\salf} (\sodbt^\hor \bu) - (\overline{\bu}-\alf \Upsilon\epsilon \overline{\bu}_s)\overline{\bu}\dif t + \Upsilon^{\salf} (\overline{\sodbt^\hor} \overline{\bu}) \biggr) \dif z = - \int_0^h \nabh \dif P \dif z.
\end{multline}
Now we express below the depth-averaged pressure term involved in the right-hand side of this equation. 

\paragraph{Dynamic pressure contribution}
Decomposing $\dif P= p\dif t + \dpt$, we obtain from Leibniz formula   
\begin{equation}
    \COR{\int_0^{h(x,t)} \nabh(\dif P)\dif z  =  \nabh (h \overline{\dif P}) - {\dif P}(h) \nabh h = \nabh (h \overline{\dif P})}. \label{eq-pressure}
\end{equation}
Moreover, remind that $\dif P$ is given by integrating the vertical momentum equation \COR{\eqref{vert-acc}}, i.e
\modif{
\beq
    - \partial_z \, \dif P = \dif t + \epsilon \beta^2 D_t w.
\eeq}
Thus, integrating \eqref{eq-pressure} \COR{vertically from depth $z$ up to the surface,} one has
\modif{
\beq
    - p(x,z,t) \dif t - \dpt(x,z,t) = (z-h)\dif t  -\epsilon \beta^2 \int_{z}^{h} D_t  w(x,z',t) \dif z', 
\eeq}
which implies
\begin{equation}
    -h(\overline{p} \dif t \COR{+} \overline{\dpt}) = -\alf h^2\dif t  -\epsilon\beta^2 \int^h_0 \biggl( \int^h_z D_t w(x,z',t)\dif z'  \biggr)\dif z.
\end{equation}

In the previous relation, the second right-hand side term can be simplified through integration by part as follows
\begin{align}
    \COR{-}\epsilon\beta^2 \int_0^h \biggl( \int_z^h D_t w \dif z' \biggr)\COR{\dif z} &= - \epsilon\beta^2 \int_0^h z \,\partial_z \biggl( \int^z_0 D_t w \dif z'   - \int_0^h D_t w \dif z' \biggr)\dif z + \bigg[ z  \biggl( \int_0^z D_t w \dif z' - \int_0^h D_t w \dif z' \biggr) \biggr]_0^h\nonumber\\
    &= - \epsilon \beta^2 \int_0^h z D_t w \dif z. 
\end{align}

Then, the horizontal momentum equation reads
\begin{multline}
    \dif_t \overline{\bu}+ \epsilon \bigl( (\overline{\bu}- \alf \Upsilon\epsilon \overline{\bu}_s) \bcdot\nabh\bigr)  \overline{\bu}\dif t + \Upsilon^{\salf} \epsilon (\sodbt^\hor \bcdot\nabh) \overline{bu} + \nabh \eta \dif t + \frac{1}{h}\epsilon \beta^2 \,\nabh \!\!\int_0^h z D_t w \dif z\\
    = - \frac{\epsilon}{h} \divh \int_0^h \biggl[ (\bu -\alf \Upsilon\epsilon \bu_s) \bu -  (\overline{\bu}-\alf \Upsilon\epsilon \overline{\bu}_s) \overline{u}\biggr]\dif z - \frac{\epsilon}{h} \Upsilon^{\salf} \divh \int_0^h \bigg[ (\sodbt^\hor \bu) -  (\overline{\sodbt^\hor} \overline{\bu})\biggr] \dif z.
    \label{eq-u-exact}
\end{multline}

So far, no approximation has been introduced in the derivation of this equation. In order to compute the different terms of this averaged equation and build simplified systems, we introduce now some asymptotic approximations. In particular, the integral terms on the right-hand side of \eqref{eq-u-exact} involve both $\bu$ and its vertical average value $\overline{\bu}$. This means the system is not closed since one cannot deduce $\bu$ from $\overline{\bu}$ solely. Therefore, one needs to investigate reliable approximations of these integral terms, as well as of the vertical acceleration term on the left-hand side of \eqref{eq-u-exact}.

\paragraph{Expansion of the velocities at the bottom}
Since the velocity potential is harmonic and irrotational, expanding the velocity component through a Taylor expansion at $z=0$ gives
\modif{
\begin{align}
    \bu(\xx,z,t) &= \bu(x,0, t) + \partial_z \bu(\xx,0,t) z + \alf \partial_{zz}^2 \bu(\xx,0,t) z^2 + \frac{1}{6} \partial_{zzz}^3 \bu(\xx,0,t) z^3 + {\cal O}(\beta^4 z^4) \nonumber\\
    &= \bu(\xx,0, t) + \beta \underbrace{\nabh w (\xx,0,t)}_{=0} z +  \alf \beta \partial_{z}\nabh w(\xx,0,t) z^2 + \frac{1}{6} \beta \partial_{zz}^2 \nabh w (\xx,0,t) z^3 + {\cal O}(\beta^4 z^4) \nonumber\\
    &= \bu(\xx,0, t) - \alf \beta^2 \nabh \divh \bu(\xx,0,t) z^2 - \frac{1}{6} \beta^3 \underbrace{\Delta_H \nabh w (\xx,0,t)}_{=0} z^3 + {\cal O}(\beta^4 z^4).
\end{align}}
Notice that, due to both harmonicity and irrotationality conditions, the expansion involves only even orders terms. Hence, since $z$ at most of the order of $h_0$, we obtain
\begin{equation}
    \label{eq-lin-u-adim}
    \bu(\xx,z,t)  = \bu(\xx,0,t) - \alf \beta^2  \nabh \divh\bu(\xx,0,t)  z^2 + {\cal O}(\beta^4).
\end{equation}

Averaging this equation, the bottom horizontal velocity $\bu^b (\xx,t)= u(\xx,0,t) $ is expressed as
\begin{equation}
    \bu^b(\xx,t) = \overline{\bu}(\xx,t)  + \frac{1}{6} \beta^2 h^2 \nabh \divh \bu^b(\xx,t) + {\cal O}(\beta^4). \label{eq-lin-u-adim-avg}
\end{equation}
Using the following expansion for ${u}^b$ in \eqref{eq-lin-u-adim-avg},
\modif{
\beq
    {\bu}^b = \overline{\bu} (\xx,t) + \epsilon \, \mbs x' \bcdot\nabh \overline{u},
\eeq}
we obtain
\begin{equation}
    \bu(\xx,z,t)  = \overline{\bu}(\xx,t) + \frac{1}{6} h^2 \beta^2 \,\nabh\divh  \overline{\bu}(\xx,t) - \alf  \beta^2 \,\nabh\divh \overline{\bu}(\xx,t) z^2 + {\cal O}(\beta^4,\epsilon \beta^4).
\end{equation}
Remark that this expansion brings an additional error term of order ${\cal O}(\epsilon \beta^4)$.

For the ISD, proceeding in the exact same way, and assuming it is irrotational at the bottom, we obtain as well
\modif{
\begin{align}
    \bu_s(\xx,z,t)  &=  \bu_s(\xx,0, t) + \partial_z \bu_s(\xx,0,t) z + \alf \partial_{z^2}^2 \bu_s (\xx,0,t) z^2 + \frac{1}{6} \partial^3_{z^3} \bu_s (\xx,0,t) z^3 + {\cal O}(\beta^4 z^4) \nonumber\\
    &= \bu_s(\xx,0, t) - \alf \beta^2 \nabh \divh \bu_s(\xx,0,t) z^2 + {\cal O}(\beta^4),
\end{align}}
where we used irrotationality at the bottom and incompressibility of the ISD. Averaging along depth, and replacing in the above expression the bottom ISD, we get
\begin{equation}\label{ISD-lin}
    \bu_s(\xx,z,t)  =  \overline{\bu}_s(\xx, t) + {\cal O}(\beta^2, \epsilon \beta^2).
\end{equation}

Furthermore, for the vertical velocity component $w$ we have
\modif{
\begin{align}
    w(\xx,y,t) &= \partial_z  w (\xx, 0,t) z + \alf \partial ^2_{z^2} w (\xx,0,t) z^2 + \frac{1}{6} \partial_{z^3}^3 w (\xx,0,t) z^3  + {\cal O}(\beta^4 z^4), \nonumber\\
    & = -  \divh \bu(\xx,0,t) z -\beta^2 \alf \underbrace{\Delta_H w(\xx,0,t)}_{=0} z^2 + \beta^2 \frac{1}{6} \Delta_H \divh \bu(\xx,0,t) z^3 + {\cal O}(\beta^4 z^4), \nonumber\\
    &= - \divh u(\xx,0,t) z + \beta^2 \frac{1}{6} \Delta_H \divh \bu(\xx,0,t) z^3 + {\cal O}(\beta^4 z^4).
\end{align}}
Making use of the expression of the bottom velocity, we get 
\begin{equation}
    w (\xx,y,t) = - \divh \overline{\bu}(\xx,t)  z + {\cal O}(\beta^2). \label{approx-w}
\end{equation}

Similarly, for the vertical component of the ISD, one has
\modif{
\begin{align}
    w_s(\xx,y,t) &= \partial_z  w_s (\xx, 0,t) z + \alf \partial ^2_{z^2} w_s (\xx,0,t) z^2 + \frac{1}{6} \partial_{z^3}^3 w_s (\xx,0,t) z^3  + {\cal O}(\beta^4 z^4), \nonumber\\
    &= - \divh \bu_s(\xx,0,t) z + \beta^2 \frac{1}{6} \Delta_H \divh \bu_s(\xx,0,t) z^3 + {\cal O}(\beta^4 z^4),
\end{align}}
then 
\begin{equation}
    w_s (\xx,y,t) = - \divh \overline \bu_s(\xx,t)  z + {\cal O}(\beta^2).
\end{equation}

For the noise term, following the same procedure, the horizontal component is expressed as
\begin{equation}\label{line-noise-u}
    \sodbt^\hor (\xx,z,t)  =  \overline{\sodbt^\hor}(\xx,t) + \frac{1}{6} \beta^2 h^2 \nabh \divh \overline{\sodbt^\hor}(\xx,t) - \alf  \beta^2 \nabh\divh  \overline{\sodbt^\hor}(\xx,t) z^2 + {\cal O}(\beta^4,\epsilon \beta^4).
\end{equation}
and the vertical one as
\begin{equation}
    \sodbt^z (\xx,y,t) = -  \divh \overline{\sodbt^\hor}(\xx,t)  z + {\cal O}(\beta^2).
\end{equation}

From these formulae, we can approximate the variance tensor components $a^H$, $a^{Hz}$ and $a^{zz}$:
\modif{
\begin{subequations}
    \begin{align}
        a^{\hor}\dif t   &= \Exp(\sdbt^\hor \otimes \sdbt^\hor)=   \Exp (\overline{\sdbt^\hor} \otimes \overline{\sdbt^\hor}) +  {\cal O}(\beta^2,\epsilon \beta^2)\nonumber\\
        &= \overline{\ba}^{\hor} dt +  {\cal O}(\beta^2,\epsilon \beta^2), \label{line-a-hh}
    \end{align}
    \begin{align}
        \ba^{\hor z}\dif t   &= \Exp(\sdbt^\hor \: \sdbt^z)=  -  \beta \Exp (\overline{\sdbt^\hor} \;  \divh \overline{\sdbt^\hor})\, z  +  {\cal O}(\beta^2,\epsilon \beta^2)\nonumber\\
        &= - \alf z \,\divh \overline{\ba}^{\hor} dt +  {\cal O}(\beta^2,\epsilon \beta^2), \label{line-a-hz}
    \end{align}
    \begin{align}
        a^{zz}\dif t   &= \Exp(\sdbt^z \: \sdbt^z)=    \Exp (\divh \overline{\sdbt^\hor} \; \divh \overline{\sodbt^\hor})\, z^2 +  {\cal O}(\beta^4,\epsilon \beta^4)\nonumber\\
        &= z^2\, \overline{\ba}^{\rm div} dt +  {\cal O}(\beta^4), \label{line-a-zz}
    \end{align}
\end{subequations}}
    
With the above approximations, we observe that all the quadratic integrals scale as
\begin{align}
    & \int_0^h \biggl[ (\bu -\alf \Upsilon\epsilon \bu_s) \bu -  (\overline{\bu}-\alf \Upsilon\epsilon \overline{\bu}_s) \overline{u}\biggr]\dif z \sim {\cal O}(\beta^4,\epsilon \beta^4),\\ 
    &\int_0^h \bigg[ (\sodbt^\hor \bu) -  (\overline{\sodbt^\hor} \overline{\bu})\biggr] \dif z \sim  {\cal O}(\beta^4,\epsilon \beta^4)
\end{align}

Besides, the ISD \eqref{aver-ISD} reads now
\modif{
\begin{align}
    \bu_s  &= \frac{1}{h}\bigl(\divh (h \overline{\ba}^{\hor})  - \ba^{\hor \hor}(h)\nabh h - \alf h \divh \ba^{\hor \hor}(h)\bigr)\nonumber\\
    &= \frac{1}{2}(\divh \overline{\ba}^{\hor}) + {\cal O}(\beta^2,\epsilon \beta^2),\\
    w_s  &= \divh \ba^{\hor z} + \partial_z a^{zz} \nonumber \\
    &= - z (\alf \divh\divh \overline{\ba}^{\hor} - 2 \overline{\ba}^{\rm{div}} ) +  {\cal O}(\beta^2,\epsilon \beta^2).
\end{align}}

Using the approximation of the vertical velocity \eqref{approx-w}, the vertical acceleration relation \eqref{vert-acc} reads
\begin{multline}
    D_t w =    - z \biggl( \dif_t  \divh \overline{\bu} +  \epsilon \bigl( (\overline{\bu} - \alf \Upsilon \epsilon \overline{\bu}_s) \bcdot\nabh\bigr) \divh  \overline{\bu}  \dif t -\epsilon\bigl( \divh( \overline{\bu}- \alf \Upsilon \epsilon \overline{\bu}_s) \divh \overline{\bu} \bigr) \dif t + \\
    \epsilon\Upsilon^{\salf}(\overline{\sodbt^\hor}\bcdot\nabh) \divh \overline{\bu} -    \epsilon \Upsilon^{\salf} ( \divh\overline{\sodbt^\hor}   \divh \overline{\bu}) \biggr) +  {\cal O}(\beta^2,\epsilon \beta^2).
\end{multline}

Eventually, we find the evolution equations for the surface elevation and the depth averaged velocity,
\begin{subequations}\label{S-SGN}
\begin{align}
    &\overline{D}^\hor_t \eta = - h \bigl(\divh (\overline{\bu} -\alf \Upsilon\epsilon \overline{\bu}_s)\dif t + \Upsilon \divh\overline{\sodbt^\hor}\bigr),\\
\begin{split}
    &\dif_t \overline{\bu}+ \epsilon \bigl( (\overline{u}- \alf \Upsilon\epsilon \overline{u}_s) \bcdot\nabh \bigr) \overline{\bu}\dif t + \Upsilon^{\salf} \epsilon (\overline{\sodbt^\hor }\bcdot\nabh)  \overline{\bu}+ \nabh \eta\dif t
    - \frac{1}{h} \epsilon\beta^2 \nabh \bigl(\frac{h^3}{3} (\dif\overline{G}) \bigr) = {\cal O}(\beta^4,\epsilon \beta^4)
\end{split}
\end{align}
\end{subequations}
with 
\begin{align}
\begin{split}
    &\dif \overline{G}(\xx,t)= \biggl( \dif_t  \divh \overline{\bu} +  \epsilon \bigl( (\overline{\bu} - \alf \Upsilon \epsilon \overline{\bu}_s) \bcdot \nabh \bigr) \divh \overline{u} \dif t - \epsilon\bigl( \divh( \overline{\bu}- \alf \Upsilon \epsilon \overline{\bu}_s) \divh \overline{\bu} \bigr) \dif t  \\ 
    & \hspace*{2cm}+ \epsilon\Upsilon^{\salf}(\overline{\sodbt^\hor}\bcdot\nabh) \divh \overline{\bu} \dif t - \epsilon \Upsilon^{\salf} ( \divh \overline{\sodbt^\hor}   \divh \overline{\bu} \bigr)\biggr),
\end{split}
\end{align}
which can be more compactly written as 
\beq
    \dif \overline{G}(\xx,t)= \overline{D}^\hor_t (\divh \overline{\bu}) -  \epsilon\divh\bigl(( \overline{\bu}- \alf \Upsilon \epsilon \overline{\bu}_s)\dif t +  \Upsilon^{\salf} \overline{\sodbt^\hor}\bigr)  \divh \overline{\bu}.
\eeq

This last expression can readily be understood as the acceleration of the average horizontal velocity divergence corrected by compressibility effects. This system constitutes a stochastic version of the Serre-Green-Naghdi equations \cite{Green-Naghdi-76, Serre-53}. Compared to the original \emph{deterministic} model, additional noise terms are involved. Those terms correspond to small scale advection and are accompanied with a corresponding ISD correction term in the large scale advection. 

Let us now exhibit some simplified models that arise by neglecting higher order terms. 

\subsection{Shallow water waves approximated models}

In this section we present two stochastic representations of classical approximations of the Serre-Green-Naghdi equations, namely the Shallow Water and the Boussinesq wave models. In addition, we briefly mention a stochastic version of the Kordeveg-De Vries equation. Through the section, we will make use of the Stokes number (also called Ursel number), defined as $S=\epsilon/\beta^2$.

\paragraph{Shallow water \modif{(or Saint-Venant)} long waves approximation}

For long waves regimes such as tidal waves, we have $\beta\ll \epsilon =A/h \ll 1$, which corresponds to a Stokes number $S\gg 1$. Neglecting in system \eqref{S-SGN} the terms of order higher than $\epsilon$ gives
\begin{subequations}\label{S-SW}
\begin{align}
    &\overline{D}^\hor_t \eta = - h \bigl(\divh (\overline{\bu} -\alf \Upsilon\epsilon \overline{\bu}_s)\dif t + \Upsilon^{\salf}  \divh \overline{\sodbt^\hor}\bigr),\\
    \begin{split}
    &\dif_t \overline{\bu}+ \epsilon \bigl( (\overline{\bu}- \alf \Upsilon\epsilon \overline{\bu}_s)\bcdot\nabh  \bigr)  \overline{\bu}\dif t + \Upsilon^{\salf} \epsilon (\overline{\sodbt^\hor }\bcdot\nabh) \overline{\bu}+ \nabh \eta\dif t=0 .
    \end{split}
\end{align}
\end{subequations}
This system corresponds to a stochastic version of the 2D Shallow water model \cite{Brecht-et-al-2021}. We note that the noise term is kept assuming $\Upsilon$ is of order 1 or higher. The diffusion term is in balance with the energy brought by the noise and must be kept to ensure energy conservation. For lower noise amplitude, the system boils down to the classical deterministic system. Notice that this system results from neglecting the vertical acceleration, which corresponds to the usual hydrostatic assumption. However, differently from the deterministic (linear) shallow water system, this corresponding linear stochastic system admits dispersive waves as solutions due to the noise term \cite{MLLTC_2023}.

\paragraph{Boussinesq approximation}

Assuming that $S\approx 1$ and $\beta \ll 1$, we retain only the terms of order $\epsilon$ and $\beta^2$ in the system \eqref{S-SGN}. Thus, we obtain a stochastic interpretation of Boussinesq wave equation. This system reads
\begin{subequations}\label{S-B}
\begin{align}
    &D^\hor_t \eta = - h \bigl(\divh (\overline{\bu} -\alf \Upsilon\epsilon \overline{\bu}_s)\dif t + \Upsilon^{\salf}  \divh\overline{\sodbt^\hor}\bigr),\\
    &\dif_t \overline{\bu}+ \epsilon \bigl( (\overline{\bu}- \alf \Upsilon\epsilon \overline{\bu}_s) \bcdot\nabh\bigr) \overline{\bu}\dif t + \Upsilon^{\salf} \epsilon (\overline{\sodbt^\hor }\bcdot\nabh)  \overline{\bu}+ \nabh \eta \dif t 
    -  \epsilon\beta^2 \bigl(\frac{h^2}{3}  \nabh\divh \dif_t  \overline{\bu}  \bigr) = {\cal O}(\beta^4,\epsilon \beta^4)
\end{align}
\end{subequations}

Remark that an additional dispersive term appears compared to the previous system. From the LU Boussinesq system, one can also derive a stochastic version of the Kordeveg-De Vries (KdV) equation, as developed in appendix \ref{sec:KdV}. 

\modif{
\subsection{Discussion on conserved quantities}

In this section, we discuss the conservation of the following quantities: mass, momentum and (mechanical) energy. We consider a bounded horizontal domain $S_H$ and assume that the noise term $\overline{\sodbt^H}$ is zero on the boundary $\partial S_H$. In particular, $\overline{\sodbt^H}$ is periodic.

\paragraph{Mass:}

We regard mass conservation first. Notice that the total mass $m$ fulfils $m = \int_{S_H}\int_0^{h}\rho \: \dif z \: \dif S_H$, where $S_H$ is the horizontal domain -- which can be 1D or 2D depending on the considered problem -- and $h= (1 + \epsilon \eta)$. Assuming that $\rho = \rho_0$ is constant, we get $m= \rho_0 h(|S_H|+\epsilon \int_{S_H}\eta \: \dif S_H)$. Hence, $\dif_t m \propto \dif_t \Big(\int_{S_H}\eta \: \dif S_H\Big)$.

Moreover, in the three models studied in this paper, the evolution equation of the surface elevation $\eta$ remains unchanged and reads
\begin{equation}
    \dif_t \eta  + \epsilon\, (\overline{\bu}^* \bcdot \nabh) \eta \, \dif t + \Upsilon^{\salf}\epsilon\, (\overline{\sodbt^\hor} \bcdot \nabh) \eta = - h \bigl(\divh \overline{\bu}^*\dif t + \Upsilon^{\salf}  \divh\overline{\sodbt^\hor}\bigr), \label{mass-nonConservative}
\end{equation}
that is to say, in conservative form,
\begin{equation}
    \dif_t \eta + \divh (\overline{\bu}^* h) \dif t + \Upsilon^{\salf} \divh (\overline{\sodbt^\hor} h) =0. \label{mass-conservative}
\end{equation}
Thus, integrating over the horizontal domain $S_H$ and using the divergence theorem, equation \eqref{mass-conservative} yields
\begin{equation}
    \dif_t \Big(\int_{S_H}\eta \: \dif S_H\Big) + \int_{\partial S_H} \Big(\overline{\bu}^* h \: \dif t + \Upsilon^{\salf}\overline{\sodbt^\hor} h \Big) \bcdot \dif n_{S_H} =0
\end{equation}

Therefore, under suitable boundary conditions -- which are periodic with a 1D domain $S_H$ in our study -- the horizontal integral of the surface elevation $\int_{S_H} \eta \: \dif S_H$ is conserved, and consequently so is the total mass.

\paragraph{Momentum:}

To investigate momentum conservation, we define the total momentum $\bold{p} = \int_{S_H}\int_0^{h} \rho \overline{\bu} \: \dif z \: \dif S_H = \rho_0 \int_{S_H} h  \overline{\bu} \: \dif S_H$. Then, we derive its evolution equation: starting from the LU Serre-Green-Naghdi model, we get
\begin{align}
    \dif_t (h \overline{\bu}) = &h \dif_t \overline{\bu} + \overline{\bu} \dif_t h = h \dif_t \overline{\bu} + \epsilon \overline{\bu} \dif_t \eta \nonumber \\
    = &- \epsilon \bigl( h \overline{\bu}^* \bcdot\nabh \bigr) \overline{\bu}\dif t - \Upsilon^{\salf} \epsilon (h \overline{\sodbt^\hor }\bcdot\nabh)  \overline{\bu}- h\nabh \eta\dif t
    + \epsilon\beta^2 \nabh \bigl(\frac{h^3}{3} (\dif\overline{G}) \bigr) \nonumber \\
    &- \epsilon \overline{\bu} \Big(\divh (\overline{\bu}^* h) \dif t + \Upsilon^{\salf} \divh (\overline{\sodbt^\hor} h)\Big) \nonumber\\
    = &- \epsilon \divh (h \overline{\bu} \otimes \overline{\bu}^*)  \dif t - \epsilon \Upsilon^{\salf} \divh (h \overline{\bu} \otimes \overline{\sodbt^\hor } ) - \underbrace{\frac{1}{2\epsilon} \nabh h^2}_{=\nabh\Big(\eta + \frac{\epsilon \eta^2}{2}\Big)} \dif t + \epsilon\beta^2 \nabh \bigl(\frac{h^3}{3} (\dif\overline{G}) \bigr). \label{momentum-conservative}
\end{align}
Similarly as for the mass, for $i \in \{x,y\}$ in 2D, and for $i=x$ in 1D, equation \eqref{momentum-conservative} yields
\begin{equation}
    \dif_t \Big(\int_{S_H} h \overline{\bu}_i \: \dif S_H\Big) + \epsilon \int_{\partial S_H} \Big(h \overline{\bu}_i \overline{\bu}^* \: \dif t + \Upsilon^{\salf}h \overline{\bu}_i \overline{\sodbt^\hor} \Big) \bcdot \dif \bold{n}_{S_H} + \int_{\partial S_H}\Big((\eta + \frac{\epsilon \eta^2}{2})\dif t - \frac{\epsilon\beta^2 h^3}{3} \dif\overline{G}\Big)\dif (\bold{n}_{S_H})_i  =0 \label{momentum-boundary},
\end{equation}
using the gradient and the divergence theorems. Remind that $\dif \overline{G}$ is defined as
\begin{equation}
    \dif \overline{G}(\xx,t)= \overline{D}^\hor_t (\divh \overline{\bu}) -  \epsilon\divh\bigl( \overline{\bu}^* \dif t +  \Upsilon^{\salf} \overline{\sodbt^\hor}\bigr)  \divh \overline{\bu}. \label{def-dG}
\end{equation}
Again, under suitable boundary conditions -- i.e periodic in our work -- it is immediate that
\begin{equation}
    \epsilon \int_{\partial S_H} \Big(h \overline{\bu}_i \overline{\bu}^* \: \dif t + \Upsilon^{\salf}h \overline{\bu}_i \overline{\sodbt^\hor} \Big) \bcdot \dif \bold{n}_{S_H} + \int_{\partial S_H} (\eta + \frac{\epsilon \eta^2}{2})\dif t \: \dif (\bold{n}_{S_H})_i = 0.
\end{equation}
Consequently, one has
\begin{equation}
    \dif_t \Big(\int_{S_H} h \overline{\bu}_i \: \dif S_H\Big) = \frac{\epsilon\beta^2 }{3} \int_{\partial S_H} h^3 \dif\overline{G} \: \dif (\bold{n}_{S_H})_i \label{momentum-dG}.
\end{equation}
The water height $h$ being periodic by assumption, it is enough to show that $\dif\overline{G}$ is periodic. In the LU Saint-Venant model, the RHS of equation \eqref{momentum-dG} is completely neglected, which is equivalent to assuming $\dif\overline{G} = 0$. In such case, momentum conservation is immediate. In the LU Boussinesq model, $\dif\overline{G}$ is approximated as
$$\dif\overline{G} = d_t (\divh \overline{\bu}) =: \dif\overline{G}_{B}.$$
Since $\overline{\bu}$ is periodic, $\divh \overline{\bu}$ also is, as long as the 1st order space derivatives of $\overline{\bu}$ are well defined. Then $\dif_t(\divh \overline{\bu})$ is periodic as well -- as long as this term is well-defined -- which proves momentum conservation. Regarding the LU Serre-Green-Naghdi model, we use the ``full'' equation \eqref{def-dG} on $\dif\overline{G}$, namely
\begin{align*}
    \dif \overline{G} = \dif \overline{G}_{SGN} :&= \overline{D}^\hor_t (\divh \overline{\bu}) -  \epsilon\divh\bigl( \overline{\bu}^* \dif t +  \Upsilon^{\salf} \overline{\sodbt^\hor}\bigr)  \divh \overline{\bu} \\
    &= \dif\overline{G}_{B} + \epsilon\, (\overline{\bu}^* \bcdot \nabh) ( \divh \overline{\bu}) \, \dif t + \Upsilon^{\salf}\epsilon\, (\overline{\sodbt^\hor} \bcdot \nabh) ( \divh \overline{\bu}) \\
    &-  \epsilon\divh\bigl( \overline{\bu}^* \dif t +  \Upsilon^{\salf} \overline{\sodbt^\hor}\bigr)  \divh \overline{\bu}.
\end{align*}
By similar arguments, the new terms on the RHS are periodic since the 2nd order space derivatives of $\overline{\bu}$ are, as long as they are well-defined. Hence, $\dif \overline{G}_{SGN}$ is periodic as well, that is momentum is conserved.

\paragraph{Energy:}

For shallow water models -- in particular, the LU Saint-Venant model -- the total energy $E_{SW}$ is defined as follows (using dimensioned velocities $\overline{\bu}$ and water height $h)$,
\begin{align}
    E_{SW} &= \int_{S_H} \int_0^h \frac{1}{2} \rho_0 \|\bu\|^2 \: dz \: dS_H+ \int_{S_H} \int_0^h \rho_0 g z \: dz \:  dS_h \nonumber \\
    &= \frac{\rho_0}{2} \int_{S_H} h \|\overline{\bu}\|^2 dS_H + \frac{\rho_0 g}{2} \int_{S_H} h^2 dS_H,
\end{align}
where the two terms on the RHS respectively correspond to the kinetic and potential energies. Scaling $\overline{\bu}$ and $h$ as before, we find the equation on the following rescaled energy equation \cite{DJ_2019, Vallis-17}
\begin{align}
    E_{SW} = \frac{\epsilon^2}{2} \int_{S_H} h \|\overline{\bu}\|^2 dS_H + \frac{1}{2} \int_{S_H} h^2 dS_H = \frac{\epsilon^2}{2} (h\overline{\bu},\overline{\bu})_{L^2(S_H, \R^2)} + \frac{1}{2}\|h\|_{L^2(S_H, \R)}^2 := E_c + E_p,
\end{align}
denoting $E_c$ and $E_p$ the scaled total kinetic and potential energies. We also denote $e_c = \frac{\epsilon^2}{2} h \|\overline{\bu}\|^2$ and $e_p = \frac{1}{2} h^2$. Now, we derive the evolution equation of this energy in the LU Saint-Venant model, using Einstein's notation on $i$,
\begin{align*}
    \dif_t e_c = &\frac{\epsilon^2}{2} \overline{\bu} \bcdot d_t (h\overline{\bu})  + \frac{\epsilon^2}{2} h \overline{\bu} \bcdot d_t\overline{\bu} \\
    = &- \frac{\epsilon^3}{2} \sum_{j \in J} \overline{\bu}_i \partial_j (h \overline{\bu}_i \overline{\bu}_j^*)  \dif t - \frac{\epsilon^3 \Upsilon^{\salf}}{2} \sum_{j \in J} \overline{\bu}_i \partial_j (h \overline{\bu}_i (\overline{\sodbt^\hor })_j ) - \frac{\epsilon}{2} h \overline{\bu}_i \partial_i h \: \dif t \bigr)\\
    &- \frac{\epsilon^3 }{2} \sum_{j \in J} h \overline{\bu}_i \overline{\bu}_j^* \partial_j \overline{\bu}_i \dif t - \frac{\epsilon^3 \Upsilon^{\salf}}{2} \sum_{j \in J} h \overline{\bu}_i (\overline{\sodbt^\hor })_j \partial_j \overline{\bu}_i - \frac{\epsilon}{2} h\overline{\bu}_i \partial_i h \: \dif t \bigr)\\
    = &-\frac{\epsilon^3 }{2} \divh (h \overline{\bu}_i^2 \overline{\bu}^*) \dif t - \frac{\epsilon^3 }{2} \Upsilon^{\salf} \divh (h \overline{\bu}_i^2 (\overline{\sodbt^\hor })_j) - \epsilon h \overline{\bu}_i \partial_i h \: \dif t \bigr),
\end{align*}
and
\begin{align*}
    \dif_t e_p = h \bcdot d_t h = - \epsilon h \partial_i (h \overline{\bu}_i^*) \dif t - \epsilon \Upsilon^{\salf} h \partial_i( h (\overline{\sodbt^\hor})_i),
\end{align*}
where $J = \{x\}$ is the problem is 1D and $J = \{x, y\}$ if it is 2D. Thus, the quantity $e = e_c + e_p$ fulfils
\begin{align}
    \dif_t e = &-\frac{\epsilon^3 }{2} \divh (h \|\overline{\bu}\|^2 \overline{\bu}^*) \dif t - \frac{\epsilon^3 }{2} \Upsilon^{\salf} \divh (h \|\overline{\bu}\|^2 \overline{\sodbt^\hor }) - \epsilon \divh ( h^2 \overline{\bu}) \dif t \\
    &+ \frac{\epsilon \Upsilon}{4} h \divh (h \overline{\bu}_s) \: dt - \epsilon \Upsilon^{\salf} h \divh (h \overline{\sodbt^\hor }). \nonumber
\end{align}
Integrating over the domain $S_H$, using the divergence theorem and periodic boundary conditions, we get by integration by parts the evolution equation of the total Saint-Venant energy $E_{SV}$,
\begin{align}
    \dif_t E_{SW} &= \int_{S_H} \Big[-  \frac{\epsilon \Upsilon}{8} \overline{\bu}_s \bcdot \nabh h^2 \: dt + \frac{\epsilon}{2} \Upsilon^{\salf} \overline{\sodbt^\hor } \bcdot \nabh h^2 \Big] dS_H\\
    &= \int_{S_H} \Big[ \frac{\epsilon \Upsilon}{8} h^2 \divh \overline{\bu}_s \: dt - \frac{\epsilon}{2} \Upsilon^{\salf} h^2 \divh \overline{\sodbt^\hor } \Big] dS_H. \nonumber
\end{align}

Consequently, for the LU Saint-Venant model -- that is assuming $\dif\overline{G} =0$ -- energy conservation is enforced by choosing the noise term $\sodbt$ such that $\divh \overline{\sodbt^\hor } = \divh \overline{\bu}_s =0$. We denote this assumption \textbf{(DF-BHNISD)}, standing for ``divergence free barotropic horizontal noise and It\={o}-Stokes drift''. However, in 1D problems, this condition does not make much physical sense since it is equivalent to considering a constant horizontal noise over space. For this reason, energy conservation is not ensured in our numerical simulations, since they were performed with more general noises which do not fulfil this assumption. Nevertheless, we anticipate that performing 2D test simulations of this stochastic Saint-Venant equation with non trivial divergence free noise would lead to numerical energy conservation results.

For the LU Serre-Green-Naghdi model, the energy is rather defined as \cite{DJ_2019, MS_1985}
\begin{align}
    E_{SGN} = \frac{\epsilon^2}{2} (h\overline{\bu},\overline{\bu})_{L^2(S_H, \R^2)} + \frac{1}{2}\|h\|_{L^2(S_H, \R)}^2 + \frac{\epsilon^3 \beta^2}{6}(h^3 \divh \overline{\bu}, \divh \overline{\bu})_{L^2(S_H, \R)}.
\end{align}
Using the same notations $e_c$ and $e_p$ as before, and defining $e_{pv} = \frac{1}{6} h^3 (\divh \overline{\bu})^2$ and $e = e_c + e_p + e_{pv}$, one has similarly
\begin{align}
    d_t E_{SGN} &= \int_{S_H} \Big[ - \frac{\epsilon^3 \beta^2 }{3} h^3 (\divh \overline{\bu}) \dif\overline{G} + \frac{\epsilon^3 \beta^2}{2} h^2 \dif_t h (\divh \overline{\bu})^2 + \frac{\epsilon^3 \beta^2}{3} h^3 (\divh \overline{\bu}) \dif_t(\divh \overline{\bu}) \Big] dS_H,
\end{align}
using that $\divh \overline{\sodbt^\hor } = \divh \overline{\bu}_s =0$. Now, computing the first term in the integrand yields
\begin{align}
     - \frac{\epsilon^3 \beta^2 }{3} h^3 (\divh \overline{\bu}) \dif\overline{G} = &- \frac{\epsilon^3 \beta^2}{3} h^3 (\divh \overline{\bu}) \dif_t(\divh \overline{\bu}) - \frac{\epsilon^4 \beta^2}{3} h^3 (\divh \overline{\bu}) (\overline{\bu}^* dt + \Upsilon^{\salf} \overline{\sodbt^\hor}) \bcdot \nabh (\divh \overline{\bu}) \nonumber\\
     &+ \frac{\epsilon^4 \beta^2}{3} h^3 (\divh \overline{\bu})^2 \divh (\overline{\bu}^* \: dt + \Upsilon^{\salf}\overline{\sodbt^\hor}),
\end{align}
that is
\begin{align}
     - \frac{\epsilon^3 \beta^2 }{3} h^3 (\divh \overline{\bu}) \dif\overline{G} + \frac{\epsilon^3 \beta^2}{3} h^3 (\divh \overline{\bu}) \dif_t(\divh \overline{\bu}) = &- \frac{\epsilon^4 \beta^2}{6} h^3 (\overline{\bu}^* dt + \Upsilon^{\salf} \overline{\sodbt^\hor}) \bcdot \nabh (\divh \overline{\bu})^2 \\
     &+ \frac{\epsilon^4 \beta^2}{3} h^3 (\divh \overline{\bu})^2 \divh (\overline{\bu}^* \: dt + \Upsilon^{\salf} \overline{\sodbt^\hor}). \nonumber
\end{align}
In addition, the second term in the integrand is
\begin{align}
     \frac{\epsilon^3 \beta^2}{2} h^2 &\dif_t h (\divh \overline{\bu})^2 = - \frac{\epsilon^4 \beta^2}{2} h^2 (\divh \overline{\bu})^2 \divh (h \overline{\bu}^* dt + \Upsilon^{\salf} h \overline{\sodbt^\hor}) \\
     &= - \frac{\epsilon^4 \beta^2}{2} h^3 (\divh \overline{\bu})^2 \divh (\overline{\bu}^* dt + \Upsilon^{\salf} \overline{\sodbt^\hor}) - \frac{\epsilon^4 \beta^2}{6} (\divh \overline{\bu})^2 (\overline{\bu}^* dt + \Upsilon^{\salf} \overline{\sodbt^\hor}) \bcdot \nabh h^3. \nonumber
\end{align}
Consequently,
\begin{align}
    \dif_t E_{SGN} &= \frac{\epsilon^4 \beta^2}{6} \int_{S_H} \Big[ h^3 (\overline{\bu}^* dt + \Upsilon^{\salf} \overline{\sodbt^\hor}) \bcdot \nabh (\divh \overline{\bu})^2 \Big] dS_H \nonumber\\
    &+ \frac{\epsilon^4 \beta^2}{6} \int_{S_H} \Big[ h^3 (\divh \overline{\bu})^2 \divh (\overline{\bu}^* \: dt + \Upsilon^{\salf} \overline{\sodbt^\hor})\Big] dS_H \nonumber\\
    &+ \frac{\epsilon^4 \beta^2}{6} \int_{S_H} \Big[ (\divh \overline{\bu})^2 (\overline{\bu}^* dt + \Upsilon^{\salf} \overline{\sodbt^\hor}) \bcdot \nabh h^3 \Big] dS_H \nonumber \\
    &= \frac{\epsilon^4 \beta^2}{6} \int_{S_H} \divh \Big[ h^3 (\divh \overline{\bu})^2 (\overline{\bu}^* dt + \Upsilon^{\salf} \overline{\sodbt^\hor}) \Big] dS_H = 0,
\end{align}
using the divergence theorem and the periodic boundary conditions again. Notice that no assumptions were made in addition to the one for the Saint-Venant model, that is \textbf{(DF-BHNISD)}. As before, this assumption leads to a space constant noise in the 1D case, therefore it is anticipated that the energy is not conserved in our simulations. Moreover, the previous calculations show that \textbf{(DF-BHNISD)} \emph{is not enough to enforce energy conservation} in the LU Boussinesq model for both the energies $E_{SW}$ and $E_{SGN}$, which is coherent with the deterministic Boussinesq model.}

\section{Numerical simulations}

In this section, we present some numerical simulations we made to test the three models derived. The Julia code that we produced is based on the work of Vincent Duchêne and Pierre Navarro, who proposed a variety of wave models implementations in the deterministic setting -- see the documentation in the following link: \href{http://waterwavesmodels.github.io/WaterWaves1D.jl/dev/}{http://waterwavesmodels.github.io/WaterWaves1D.jl/dev/} \cite{website_DN_2024}. \modif{These models are essentially based on pseudo-spectral resolution methods, which justifies the use of periodic boundary conditions.} We adapted their \modif{numerical framework} to the stochastic case, introducing implementations of the noise terms and the ISDs for this purpose.

Regarding the purely stochastic aspects, we consider noises with wave spatial structure. This is justified by the shape of solutions found in \cite{MLLTC_2023}: considering a constant noise a first, the authors showed that the system admits progressive wave solutions. Then, they extend the analysis to systems where the noise is itself a progressive wave. In the end, the \modif{1D} noise we consider is the following,
\begin{equation}
    \sigma(x) dW_t = A \cos(kx) d\beta_t^1 + A \sin(kx) d\beta_t^2,
\end{equation}
where $A$ denotes the amplitude of the noise and $d\beta_t^1, d\beta_t^2$ denote Brownian motions. Using a Box-Muller argument, this shape is equivalent to
\begin{equation}
    \sigma(x) dW_t = A \cos(kx + \phi_t) d\beta_t,
\end{equation}
where $\phi_t$ is a uniformly distributed random phase on $(-\pi,\pi)$, such that for all $t \neq s$, $\phi_t$ and $\phi_s$ are independent. Additionally, $d\beta_t$ is a Brownian motion. In our simulations, the noise wave number is set to \modif{$k=2\pi/100$}, and the noise amplitude may take the following values: $A=0.001$, $A=0.005$ or $A=0.01$. \modif{Notice that the dimensioned noise scales like $A \epsilon \sqrt{gh_0}$ as a consequence. The value of wave number $k$ is chosen to be at least one order of magnitude smaller than the typical wave number of the deterministic wave. This is because our simulations showed that noise terms with too small space scale oscillations lead to numerical instability, and enables us to further discuss the presence of an additive noise term in the water elevation dynamics. The values of amplitude where chosen to be much smaller than the typical height of the wave. Numerically, we have observed that $A=0.001$ yields slight perturbations of the deterministic waves -- that is typical realisations of each LU model is similar to its deterministic counterpart -- while $A=0.01$ induces a more ``noisy'' dynamics -- that is typical realisations are essentially noise driven. In addition, we chose $A=0.005$ as an enlightening intermediate case.}

\begin{figure}[h] 
    \centering
    \includegraphics[width=7.5cm]{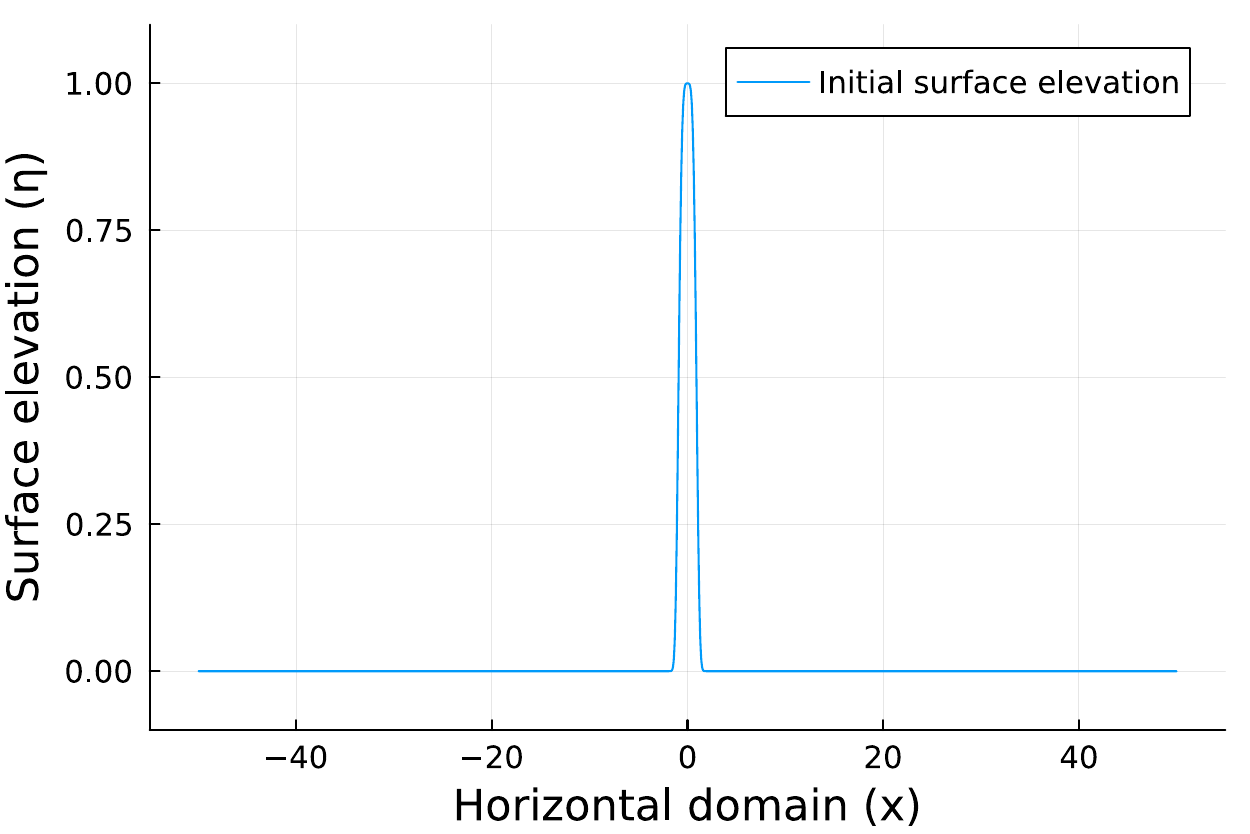}
    \includegraphics[width=7.5cm]{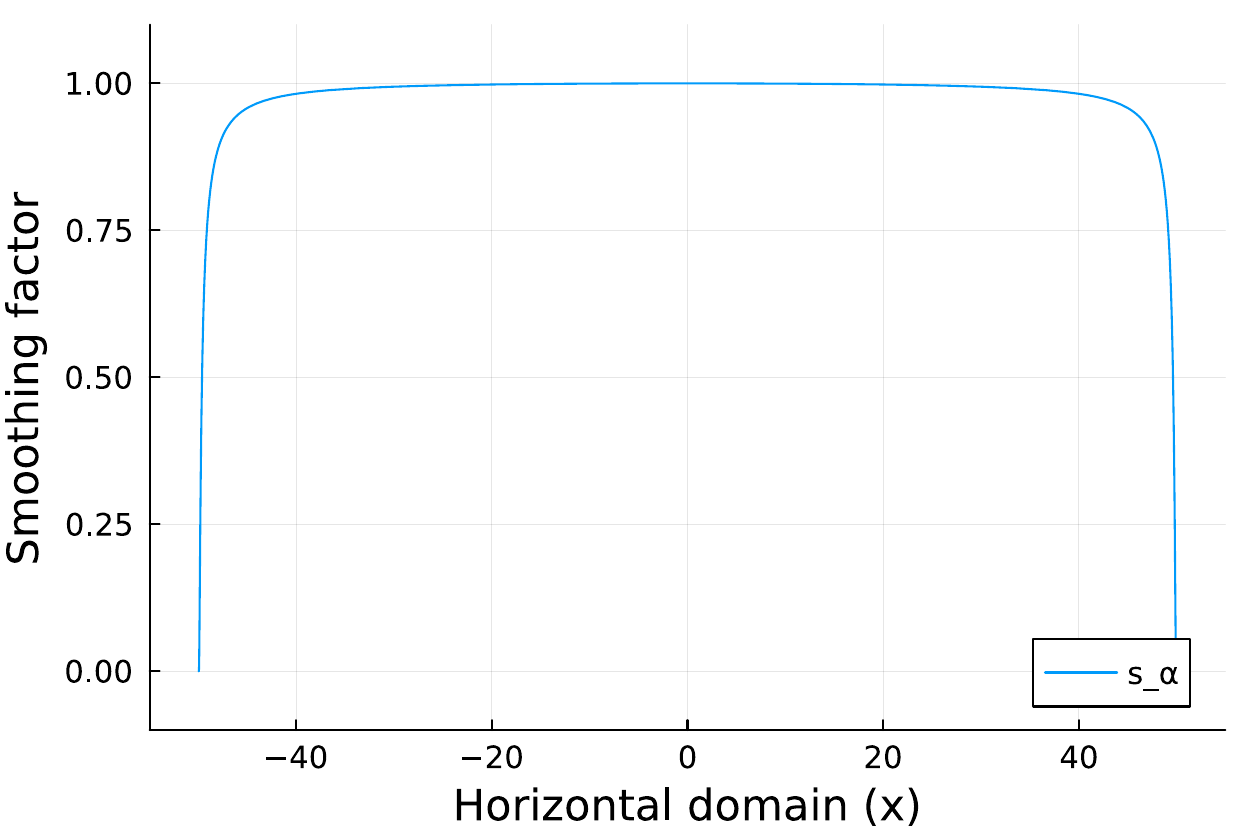}
    \caption{Initial surface deformation (left) and boundary conditions enforcement function $s_\alpha$ (right, $\alpha = 10$).} \label{fig:initCond-BCfactor}
\end{figure}

Our tests are based on computing the evolution of the deformation surface $\eta$, with a ``heap of water" type initial condition. Namely, the initial surface deformation is set to be $\eta(x,t=0) = \exp( - x^4)$, while the initial velocity is set to $u(x,t=0) = 0$. All of our tests were performed on a numerical 1D tank $[-L,L]$ with $L=50$, which is discretised with $N=2^{11}$ spatial points. The timestep is chosen to be $dt = 0.005s$, and we assume periodic boundary conditions. To enforce these conditions on the noise terms as well, we multiply them the function $s_\alpha(x) = \exp\Big(\frac{1}{\alpha^2}\big(1-\frac{1}{1-(x/L)^2}\big)\Big)$, with $\alpha=10$, in order to make them vanish on the boundary. The initial condition on $\eta$ and the profile of $s_\alpha$ are given on figure \ref{fig:initCond-BCfactor}.

Moreover, we propose two sets of parameters for testing our models
\begin{itemize}
    \item $\beta = 0.01$ and $\epsilon = 0.1$ $(\mathcal{P}_1)$, so that the scaling conditions associated with the Shallow Water model hold (and in particular, they hold for the Boussinesq model as well). Therefore, the three models should give qualitatively similar results. \modif{
    \item $\beta=0.1$ and $\epsilon= 0.1$ $(\mathcal{P}_2)$, so that the scaling conditions of the Shallow Water and the Boussinesq models do not hold. Thus, these models should give qualitatively distinct -- yet not dramatically different -- results.}
\end{itemize}
Additionally, we will compare our stochastic models to the deterministic one dimensional Water Waves model.

Furthermore, two options are available for simulating our models, depending on whether we consider a ``true" Stratonovitch noise or if we rewrite it in It\={o} form. In the former case, we may adopt a stochastic Euler-Heun approach as developed in \cite{XT_2016}, where the computation of the stochastic diffusion becomes implicit. In the latter case however, one needs to give an analytical expression for the stochastic diffusion. Even though one can compute them in the LU Saint-Venant and the LU Boussinesq systems, this correction term becomes extremely complex in the LU Serre-Green-Naghdi model due to the term $\overline{dG}$. Therefore, we rather adopted the first approach using the stochastic Euler-Heun method.

Moreover, for numerical stability reasons, we used a stochastic version of the order 4 Runge-Kutta algorithm (RK4), rather than the (simpler) Euler-Maruyama algorithm. For instance, the (deterministic) Saint-Venant equations are known to be dramatically more stable when solved with the order 4 Runge-Kutta rather than the Euler method. In summary, the solving algorithm we used is essentially the following: we treat the bounded variation term as in the classical RK4 method, and the martingale term as in the stochastic Euler-Heun method. \modif{Such approach has been studied in more details in \cite{GF_2019, PS_2008, FBLM_2023}.}

\begin{figure}[h] 
    \centering
    \includegraphics[width=7.5cm]{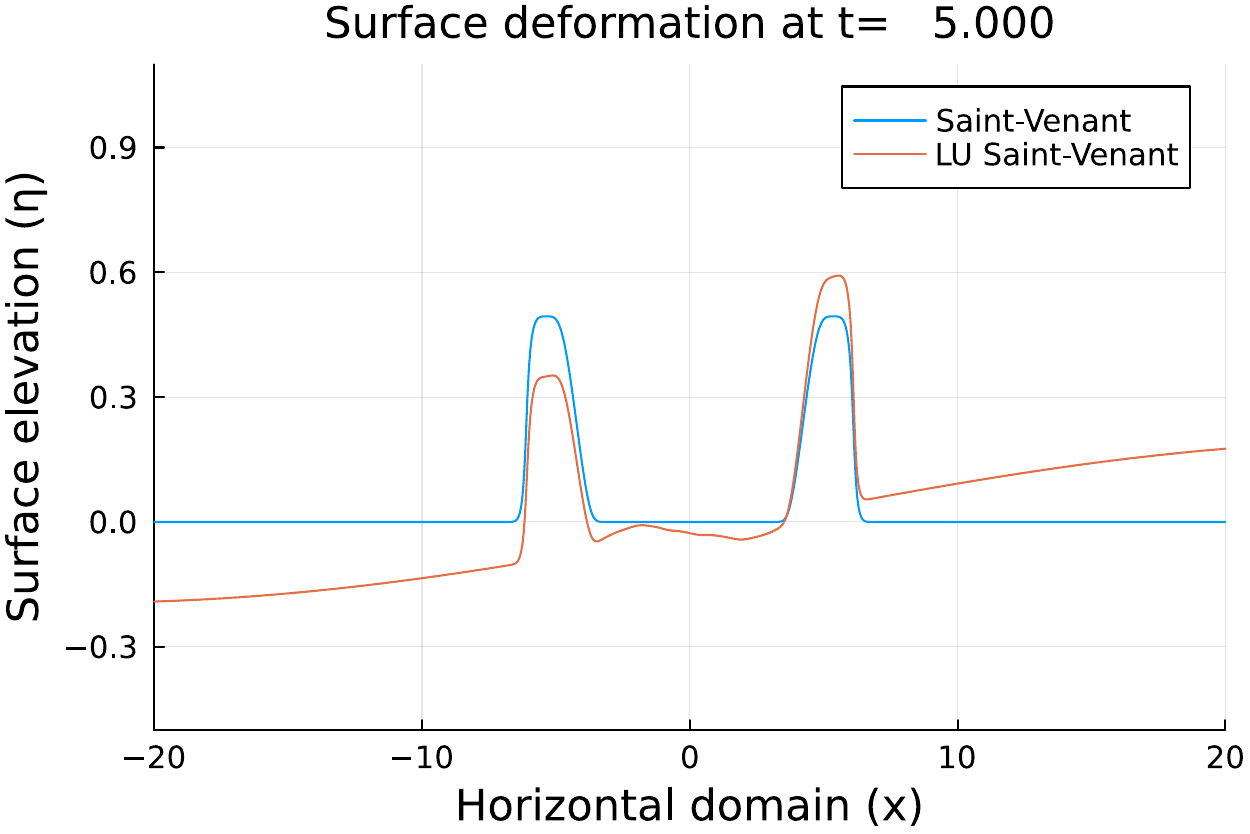}
    \includegraphics[width=7.5cm]{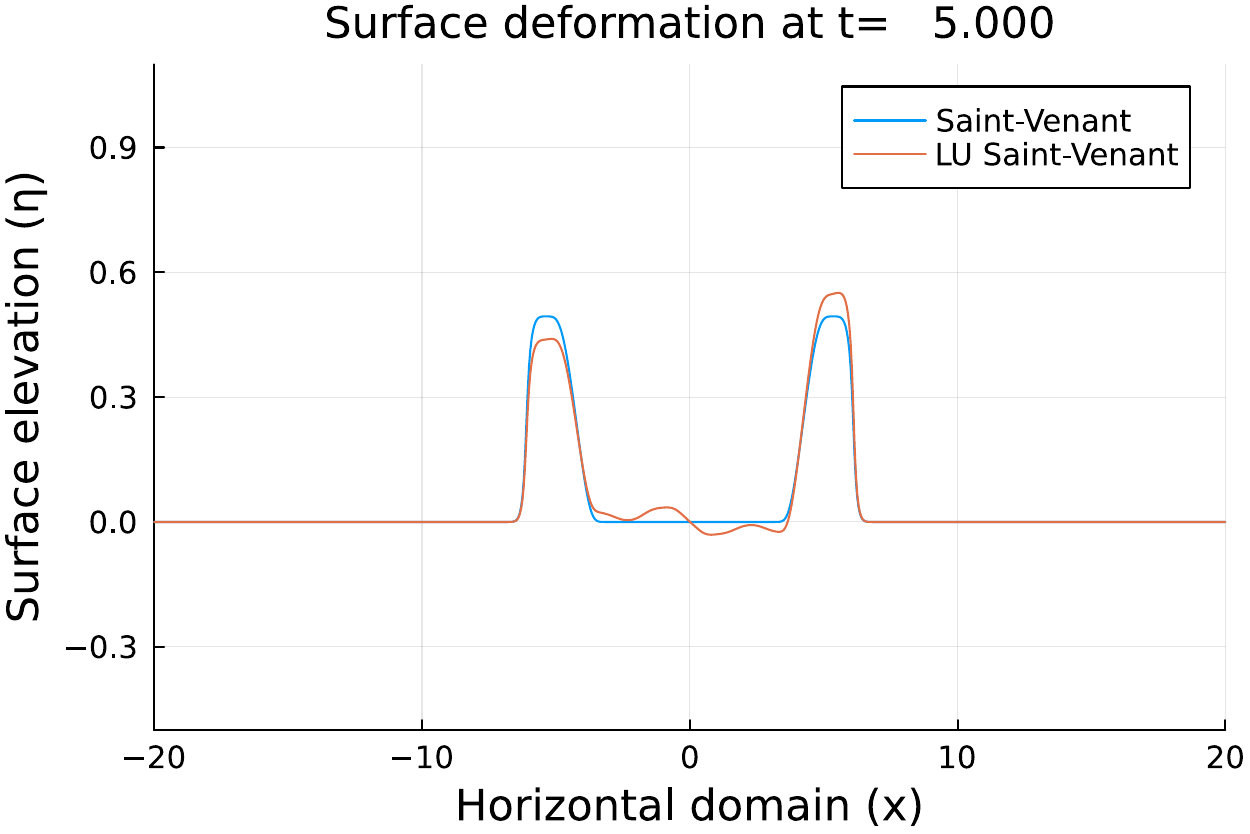}
    \caption{\modif{Realisations of the LU Saint-Venant model, with additive noise (left) and without (right). We only plot the solution over the domain $[-20,20]$, at $t = 5s$. Parameter set: $(\mathcal{P}_1)$ -- Wave number: $k = 2\pi/100$ -- Amplitude: $0.001$.}} \label{fig:additiveNoise}
\end{figure}

\modif{Furthermore, since water elevation equations on $\eta$ are the same in the three models we study -- regarding for example equation \eqref{mass-nonConservative} -- one may notice the presence of the term
$$\Upsilon^{\salf} h \divh (\overline{\sodbt^\hor}) = \Upsilon^{\salf} \divh (\overline{\sodbt^\hor}) + \epsilon \Upsilon^{\salf} \eta \divh (\overline{\sodbt^\hor}).$$
This shows the existence of an additive noise term $\Upsilon^{\salf} \divh (\overline{\sodbt^\hor})$ in the water elevation dynamics. Its effect is illustrated by Figure \ref{fig:additiveNoise}: due to this term, the surface elevation is \emph{not flat} ``away from the wave'', which is a strong difference compared to the deterministic setting. To facilitate the comparison between our models and their deterministic counterparts, we chose to \emph{disregard} this additive noise term in our simulations. This can be interpreted as a filtering of the lower wave numbers -- i.e. large scale dynamics.}

\subsection{\modif{Qualitative analysis} on the effect of the noise}

In this subsection, we give some qualitative insights about the effect of the noise on the dynamics of the system. For this purpose, we compare each deterministic solution to a solution of the associated LU interpretation, \modif{ disregarding the effect of the additive noise term previously mentioned.}

\subsubsection{Parameters $(\mathcal{P}_1)$}

\begin{figure}[htbp]
    \centering
    \begin{tabular}{c|c|c}
        \includegraphics[width=5cm]{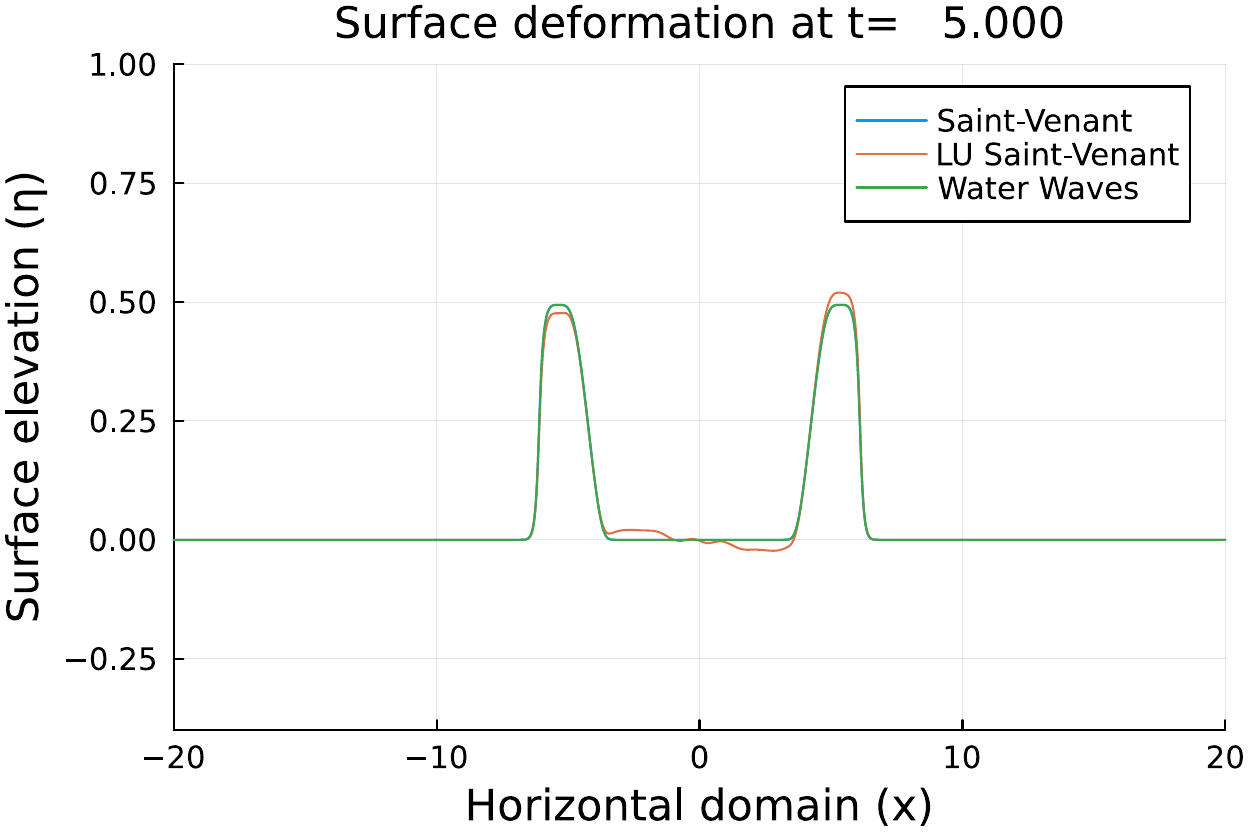} &
        \includegraphics[width=5cm]{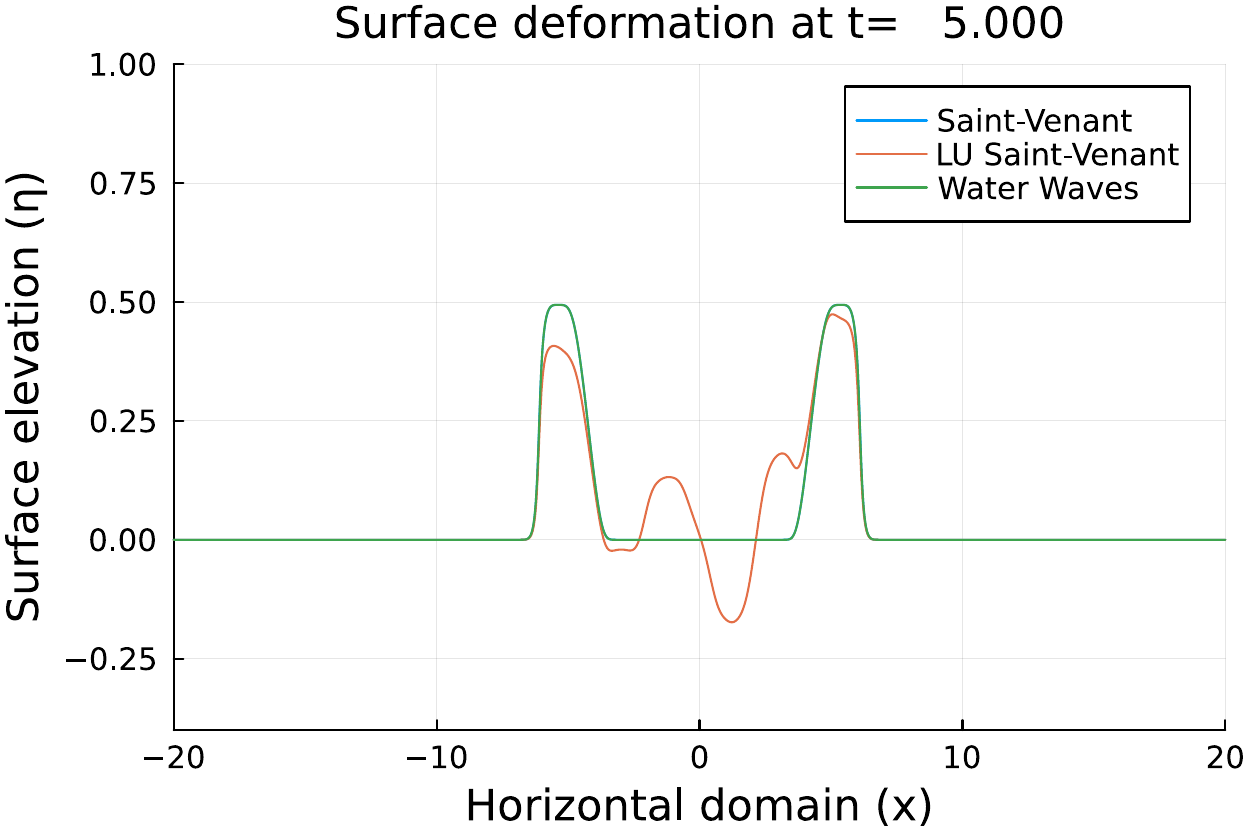} &
        \includegraphics[width=5cm]{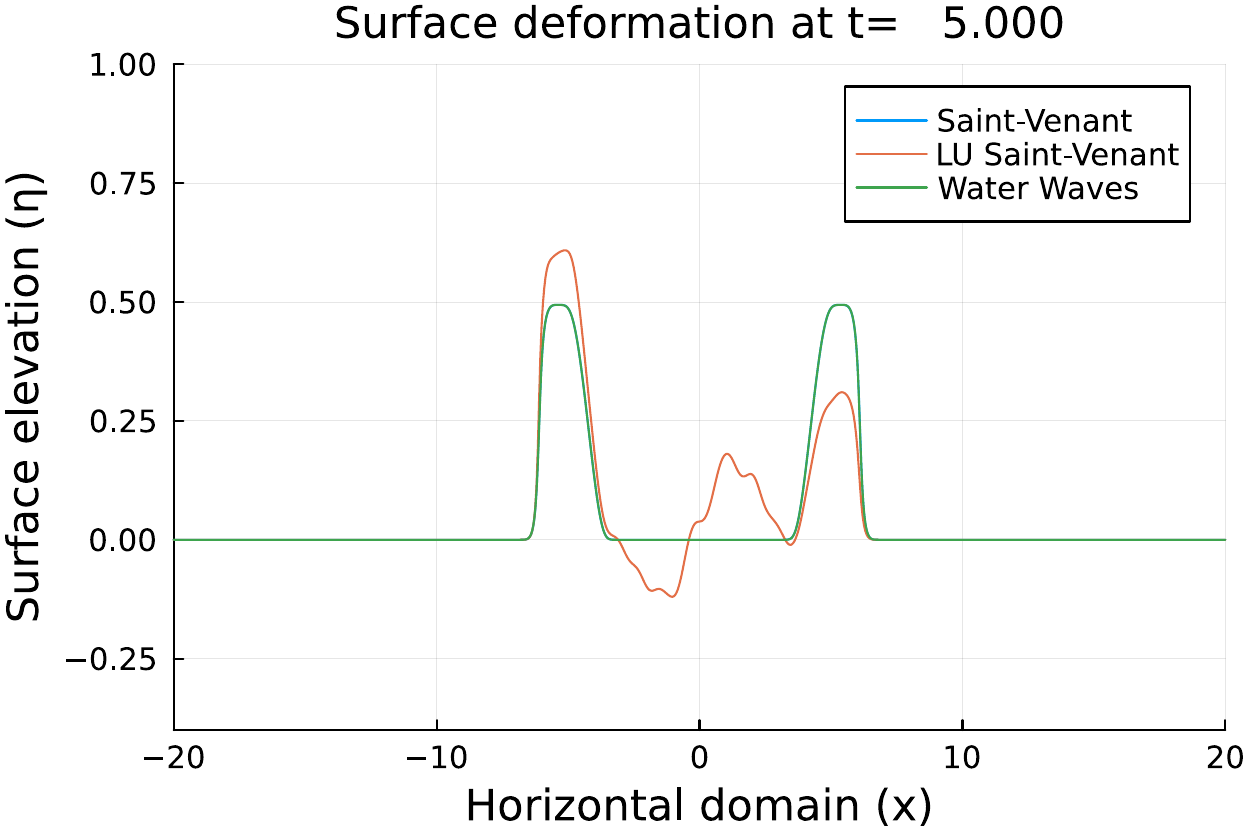}\\
        \hline
        \includegraphics[width=5cm]{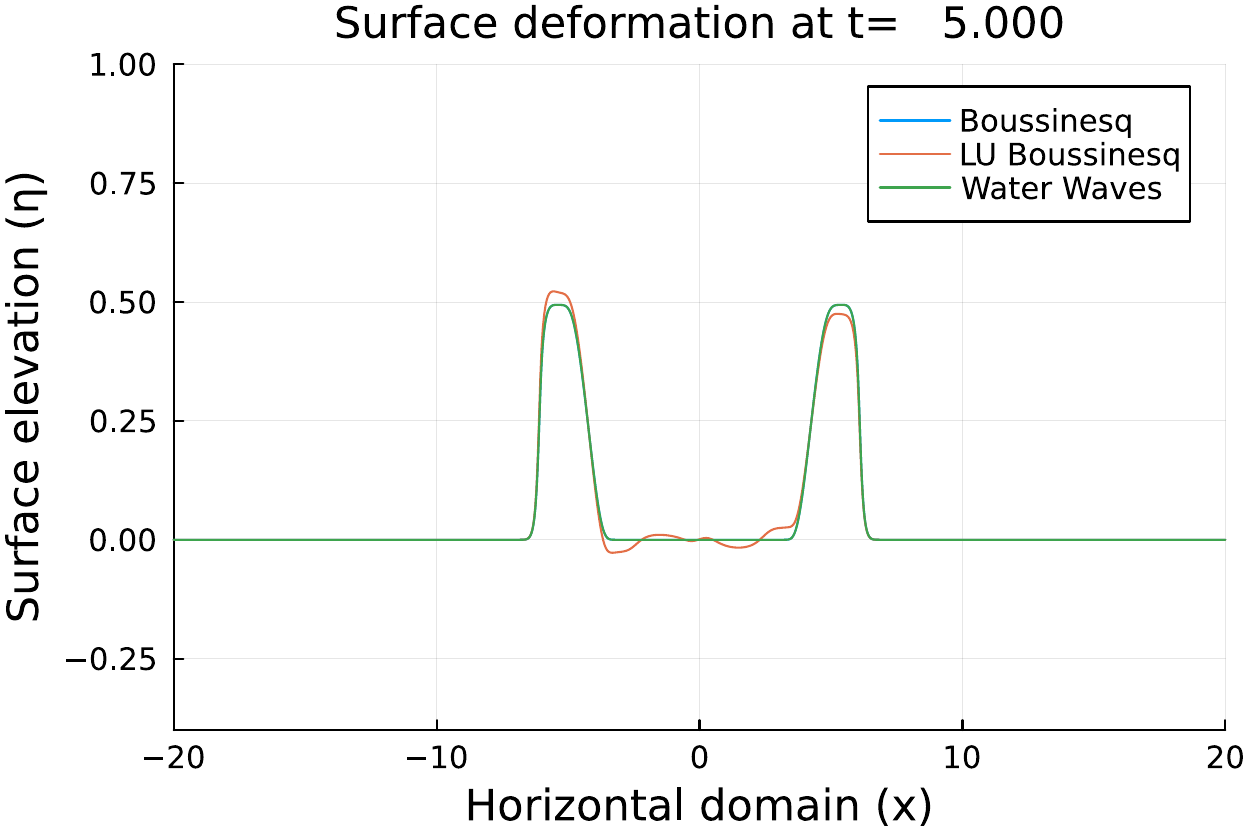} &
        \includegraphics[width=5cm]{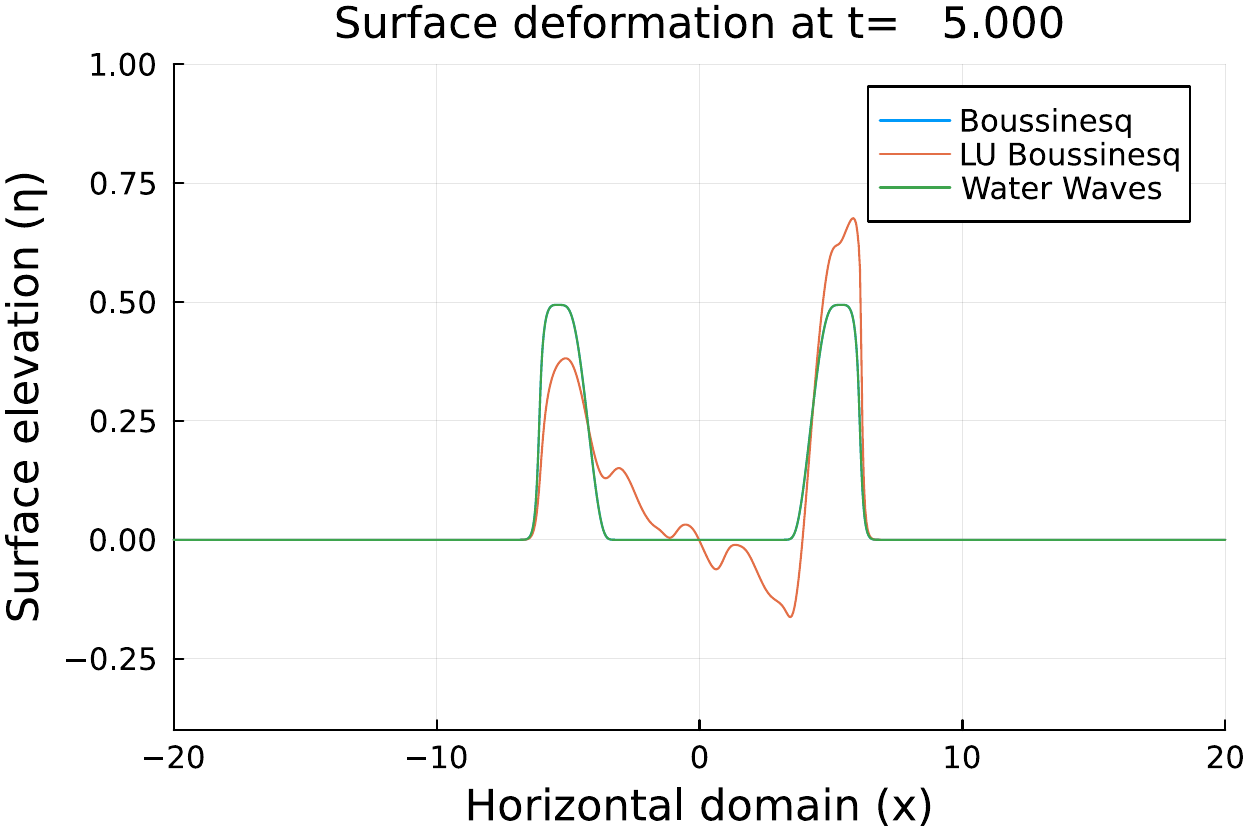} &
        \includegraphics[width=5cm]{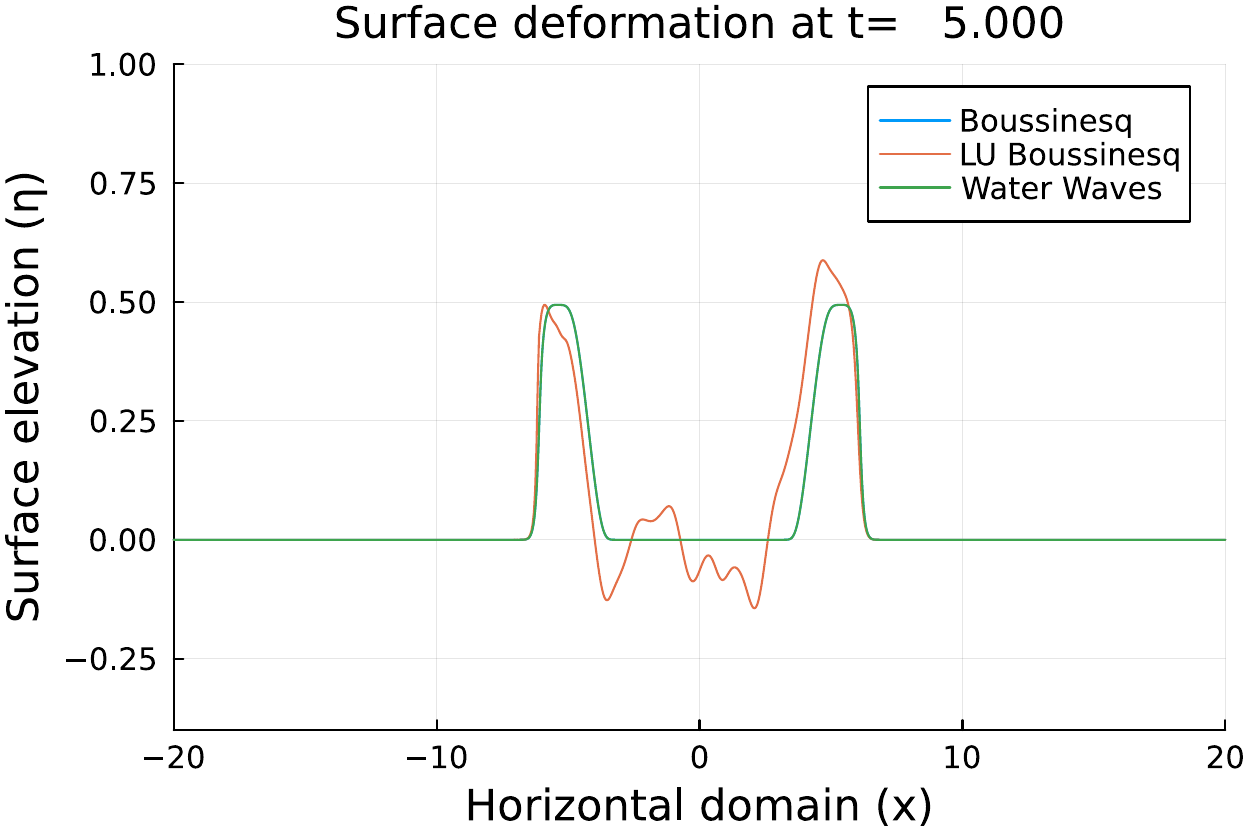}\\
        \hline
        \includegraphics[width=5cm]{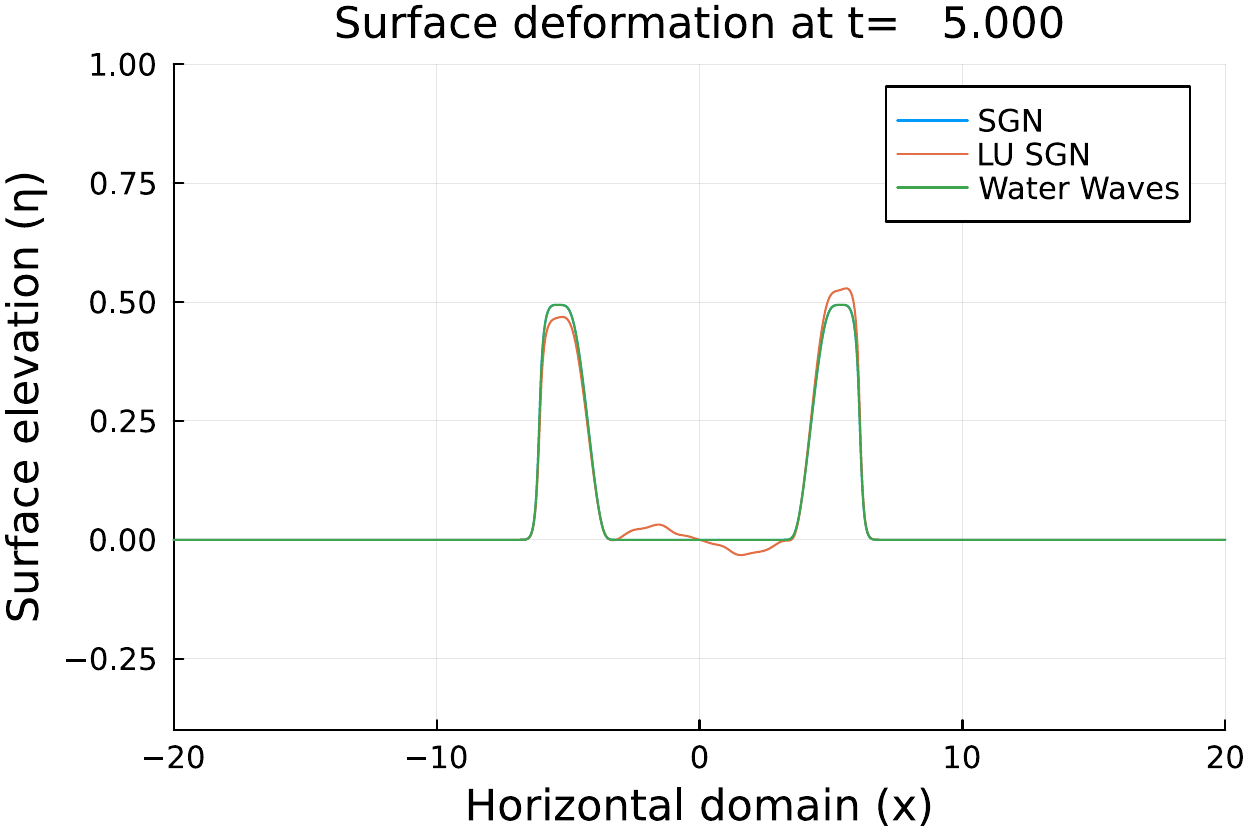} &
        \includegraphics[width=5cm]{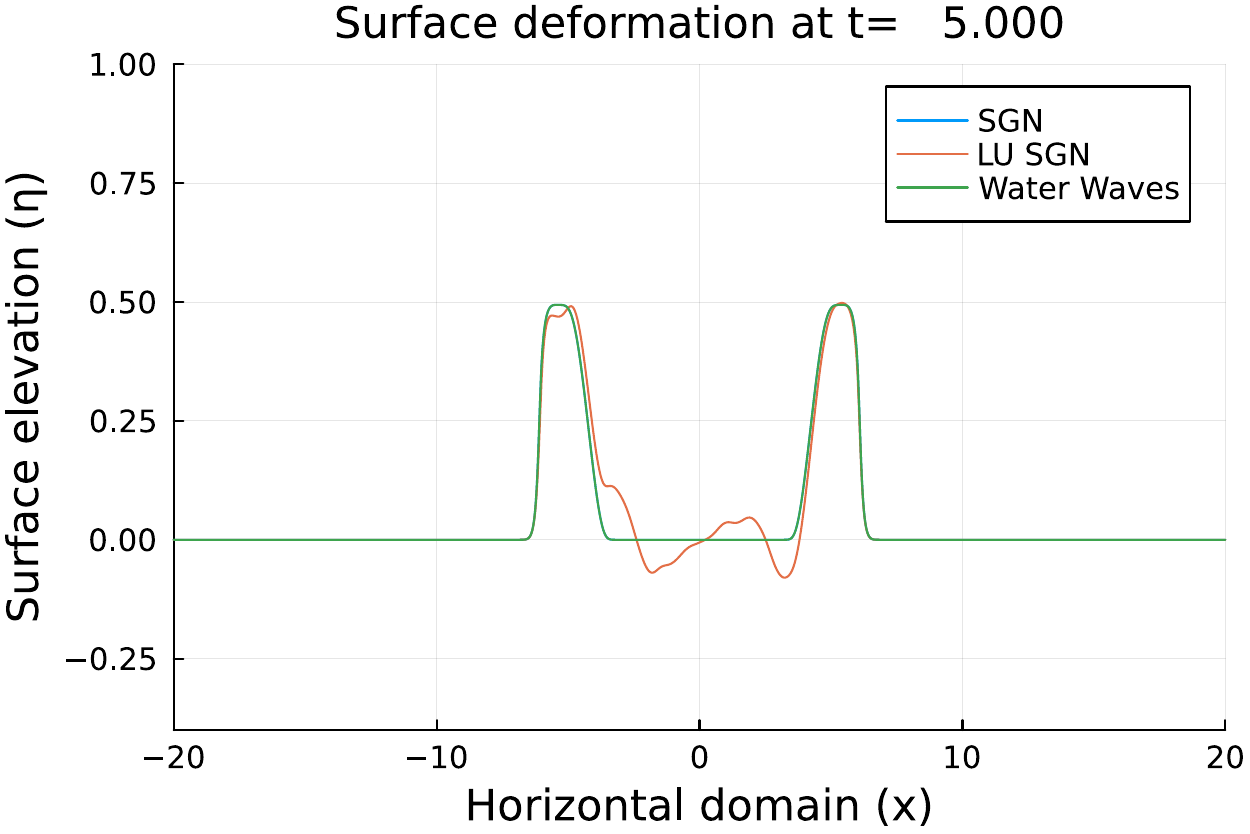} &
        \includegraphics[width=5cm]{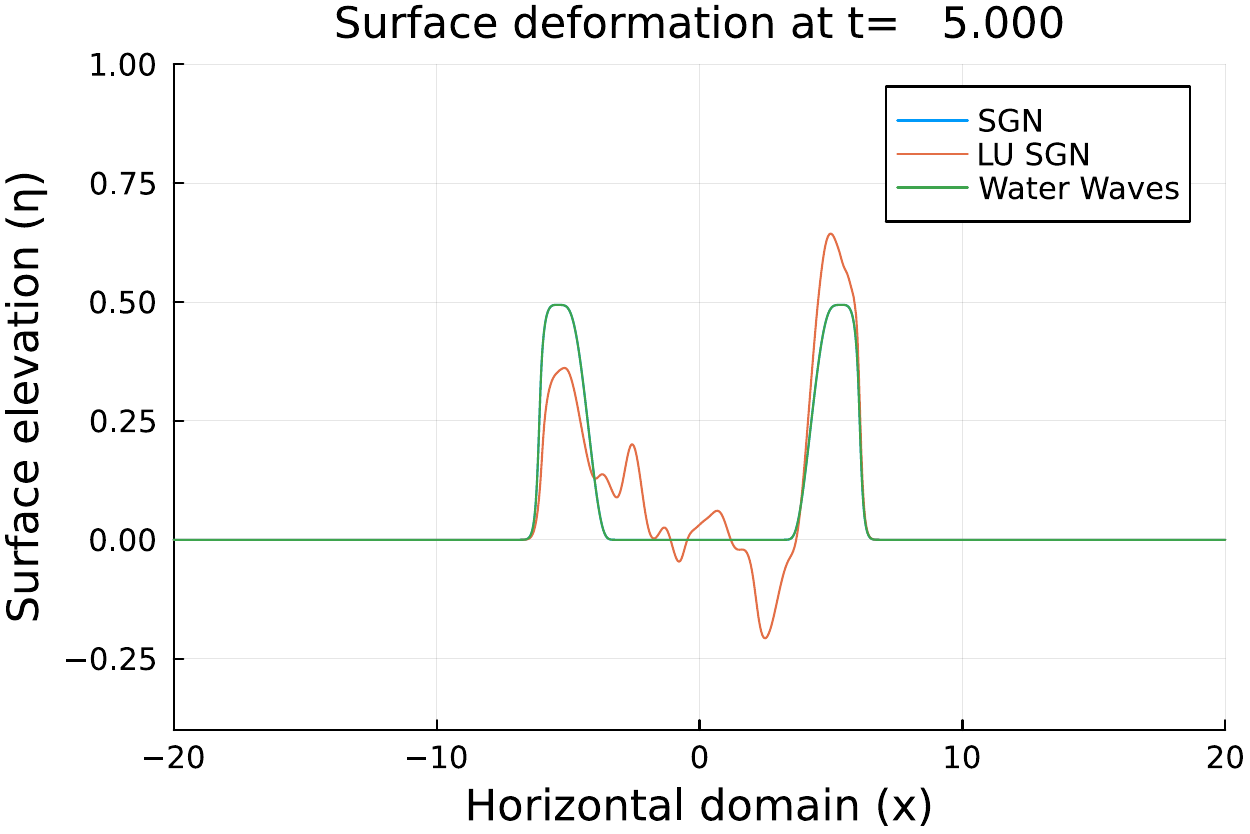}
    \end{tabular}
    \caption{\modif{Comparison between the surface deformation of the deterministic Saint-Venant (1st row), Boussinesq (2nd row) and Serre-Green-Naghdi (3rd row) models and realisations of their LU interpretations. Parameter set: $(\mathcal{P}_1)$ -- Wave number: $k = 2\pi/100$ -- Amplitude, from left to right: $0.001$, $0.005$ and $0.01$.}}\label{fig:plots-p1}
\end{figure}

Considering the parameters $(\mathcal{P}_1)$, we compared the deterministic Saint-Venant, Boussinesq and Serre-Green-Naghdi models (blue curves) to realisations of their associated LU models (orange curves), using the 1D Water Waves as a reference (green curves). The noise we chose has the shape of a stationary wave, with wave number \modif{$k = 2\pi/10$}. Figure \ref{fig:plots-p1} show \modif{realisations} of these models at $t = 5s$, for different noise amplitudes.

We observe that deterministic solutions are close to the Water Waves solution. This is expected regarding the scaling of $\beta = 0.01 \ll 1$. Moreover, the LU systems seem to converge to their respectively associated deterministic solution for a vanishing noise, thus giving consistency to the LU interpretation of the wave equations.
In both cases, the noise tends to break the spatial symmetry of the wave. This stems from the shape of the noise, which is a sum of a symmetric function -- $A\cos(kx)d\beta_t^1$ -- and an antisymmetric function -- $A\sin(kx)d\beta_t^2$.

\subsubsection{Parameters $(\mathcal{P}_2)$}

\begin{figure}[htbp]
    \centering
    \begin{tabular}{c|c|c}
        \includegraphics[width=5cm]{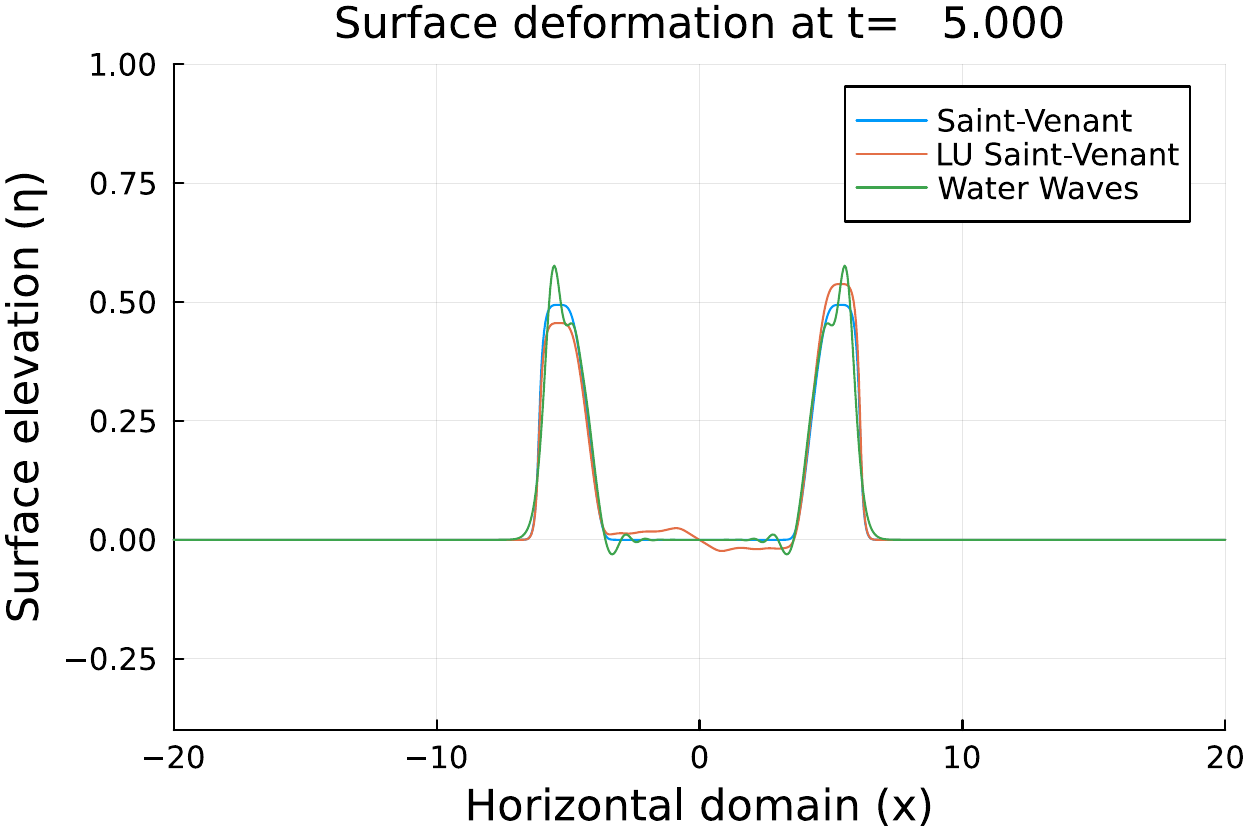}&
        \includegraphics[width=5cm]{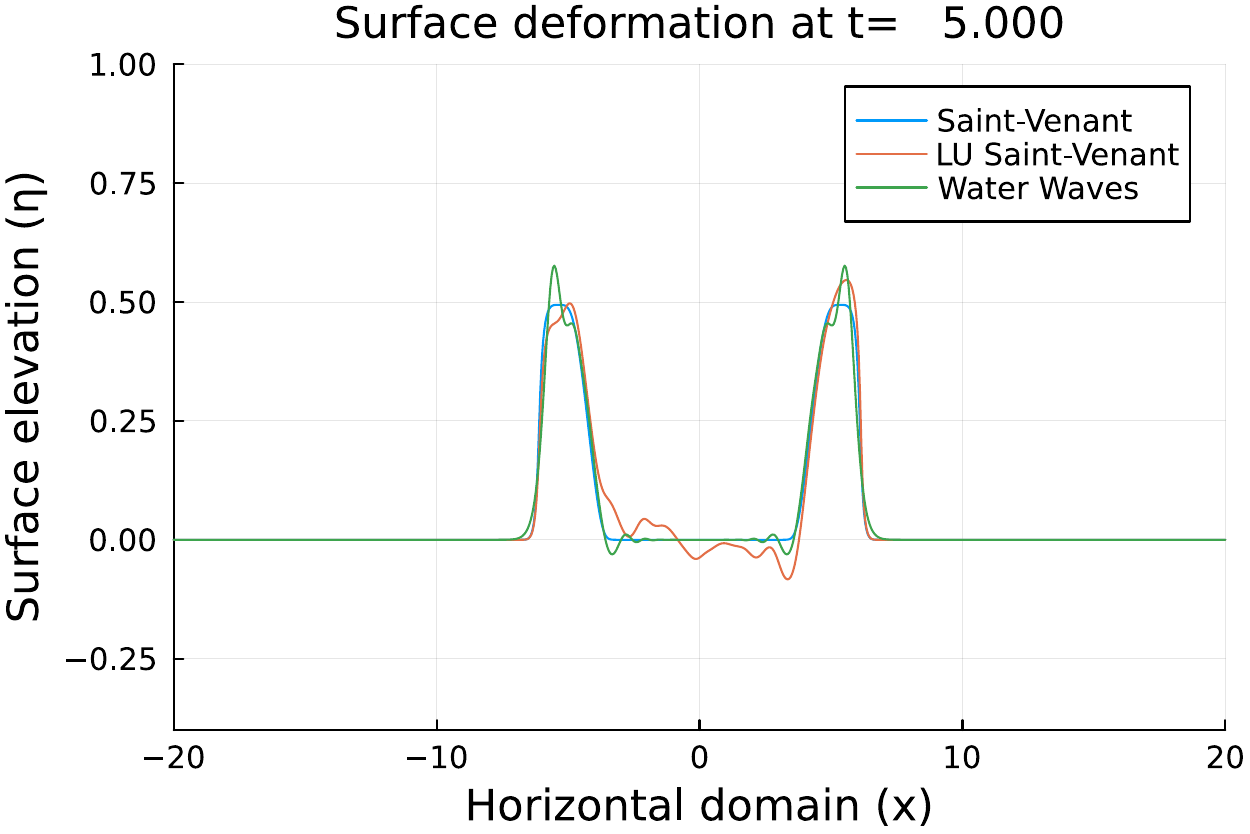}&
        \includegraphics[width=5cm]{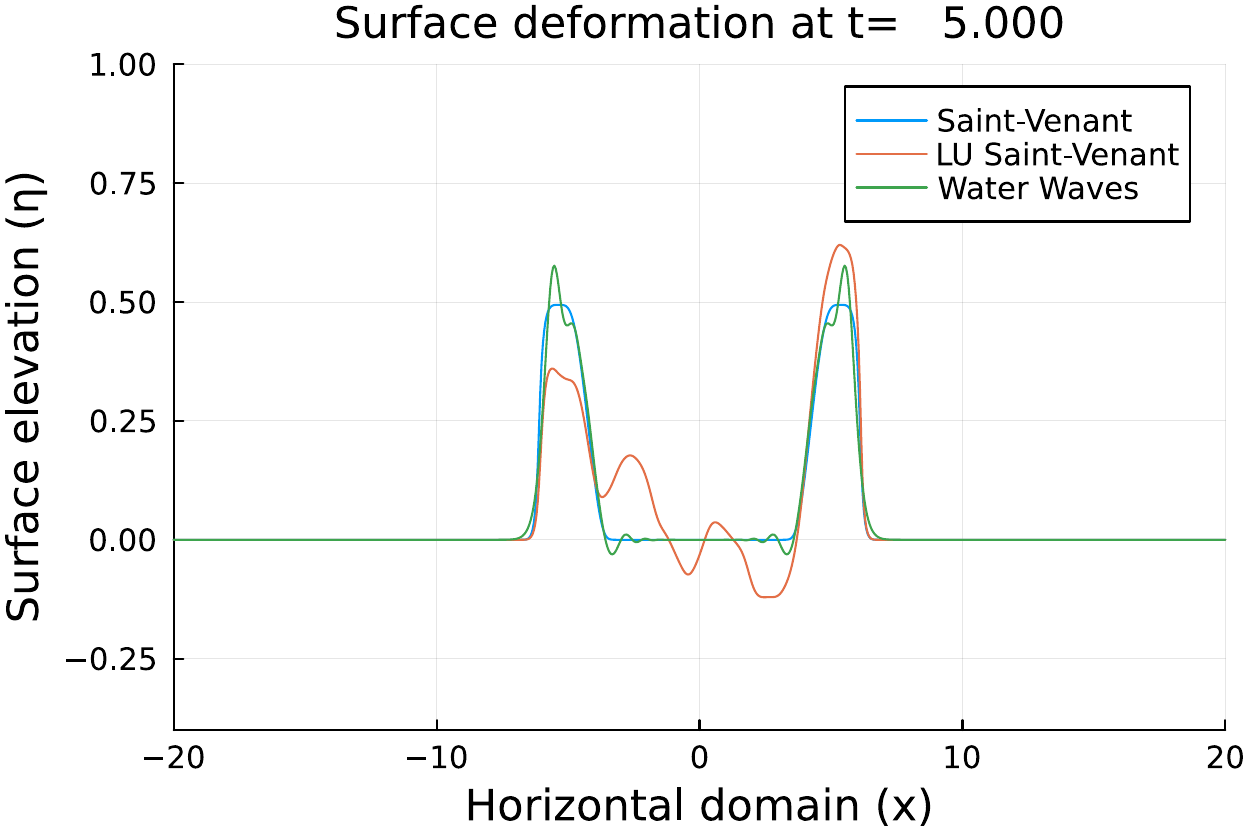}\\
        \hline
        \includegraphics[width=5cm]{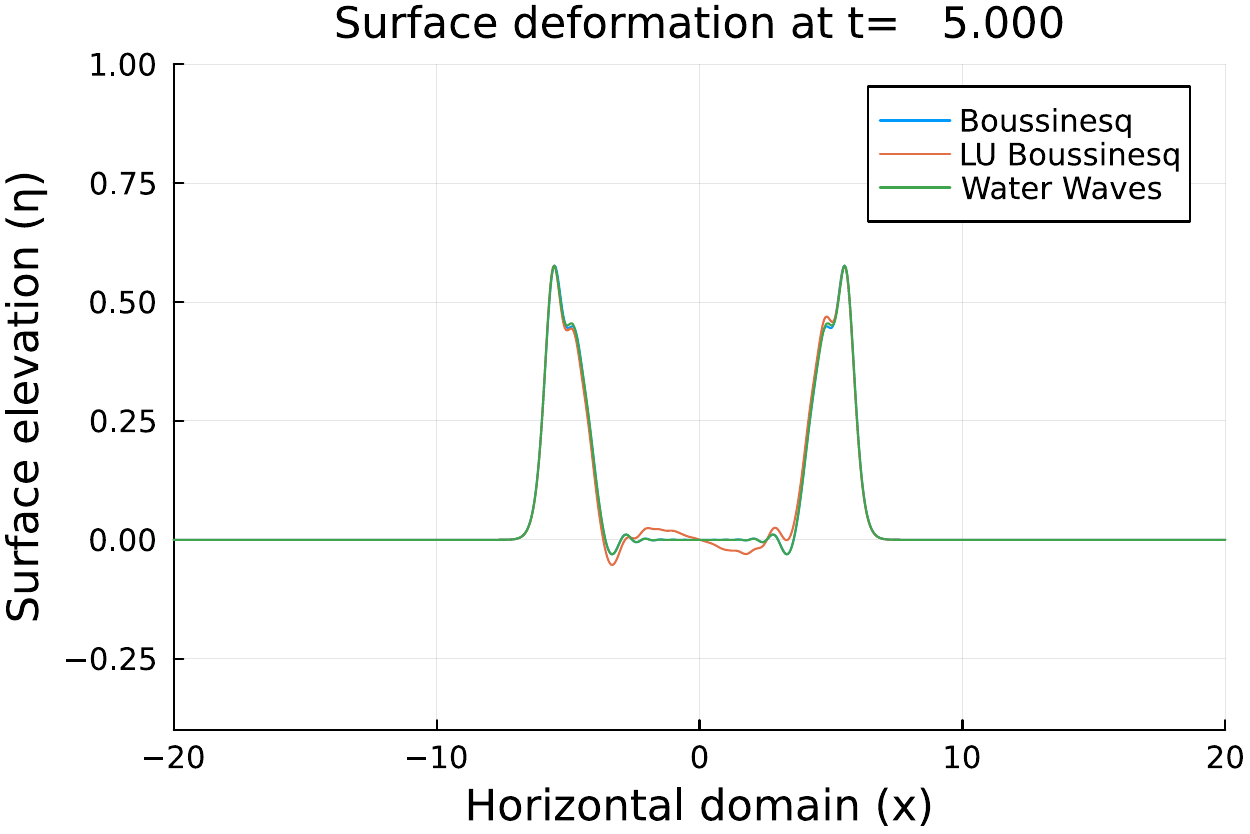}&
        \includegraphics[width=5cm]{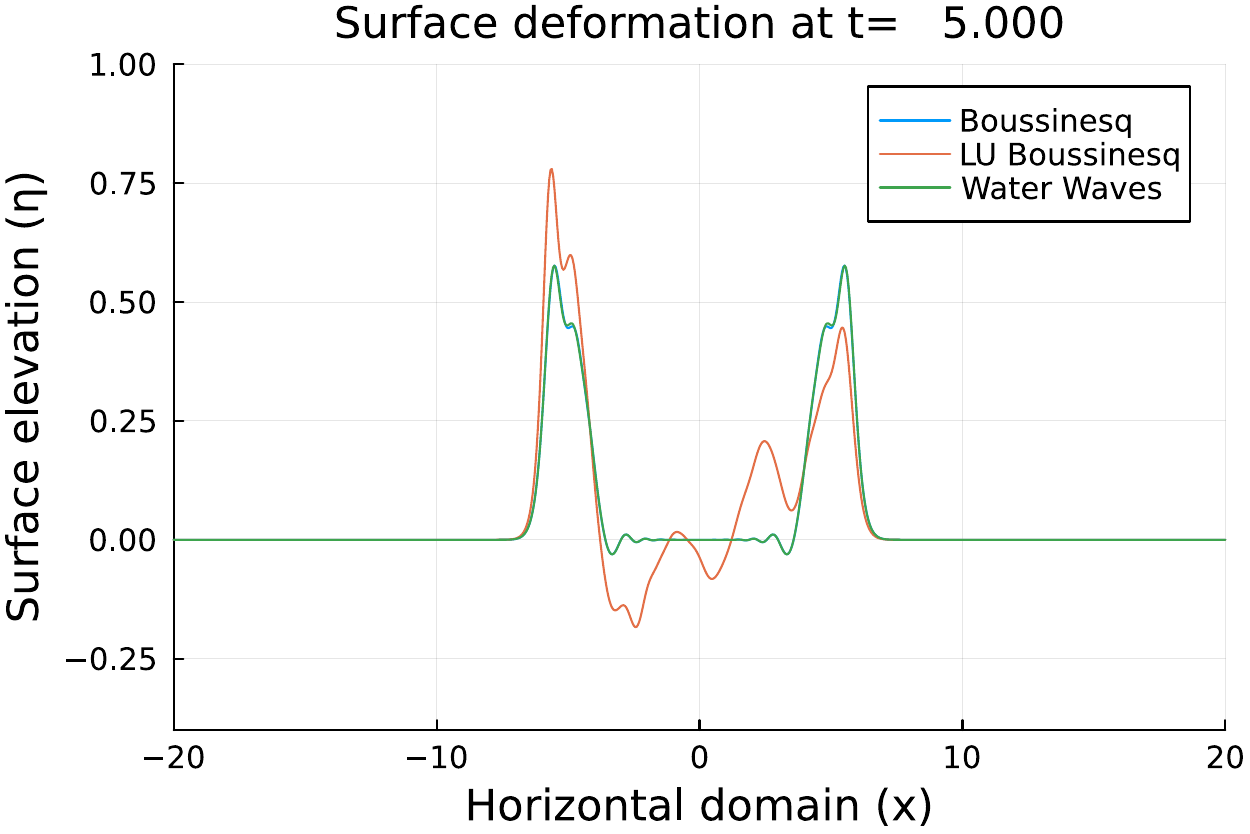}&
        \includegraphics[width=5cm]{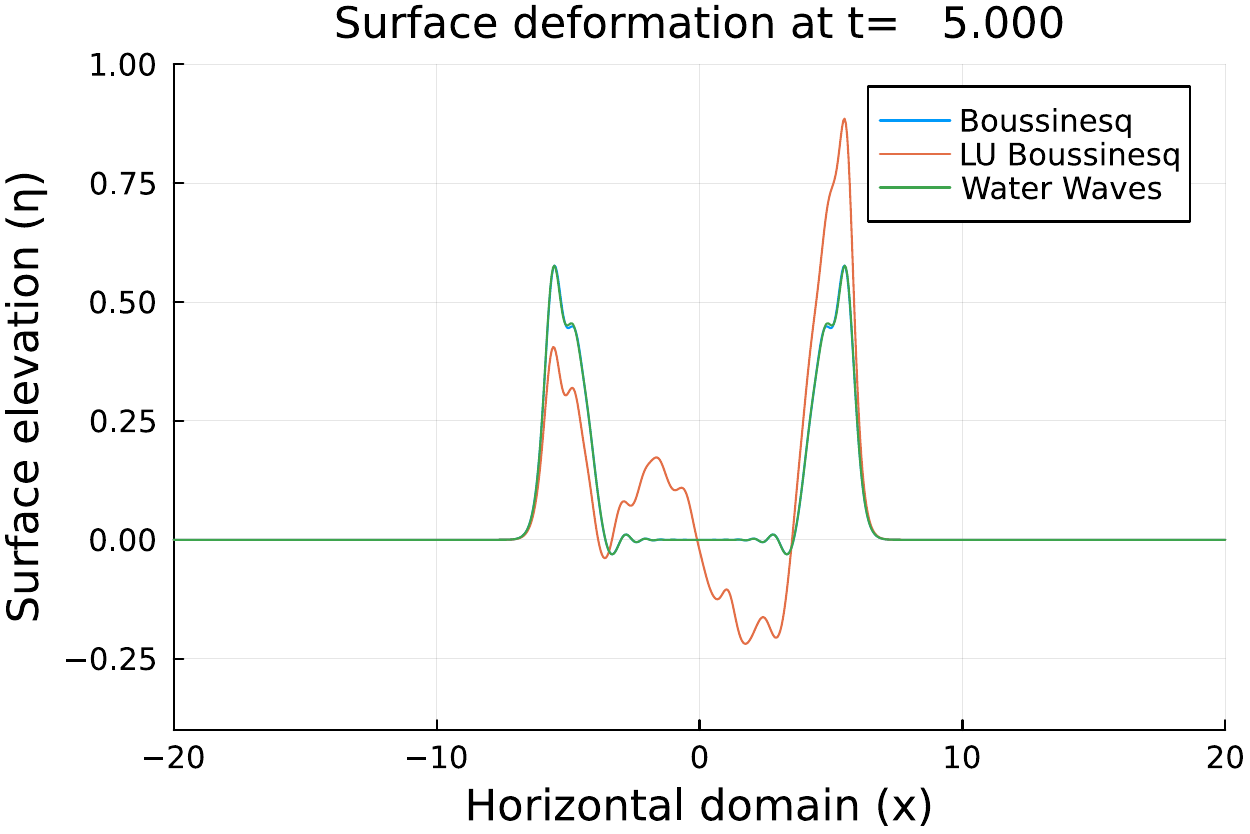}\\
        \hline
        \includegraphics[width=5cm]{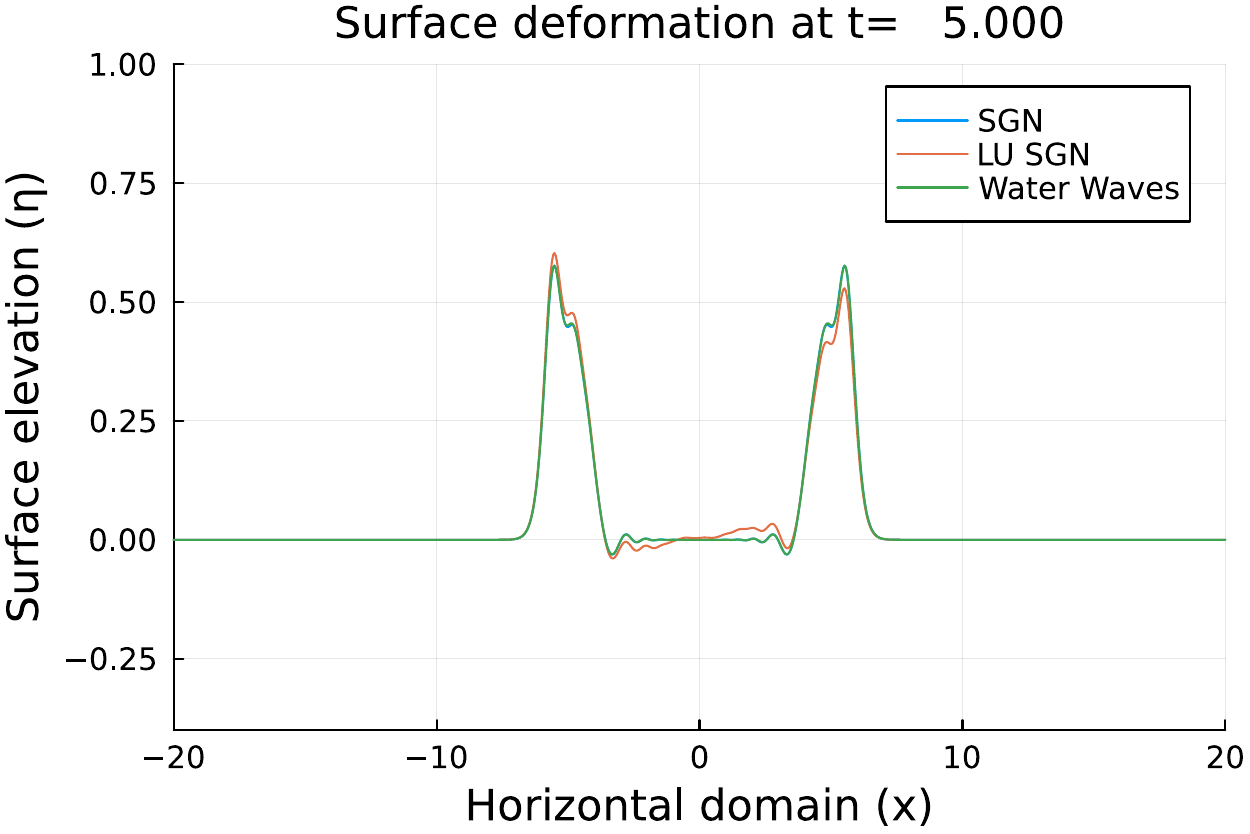}&
        \includegraphics[width=5cm]{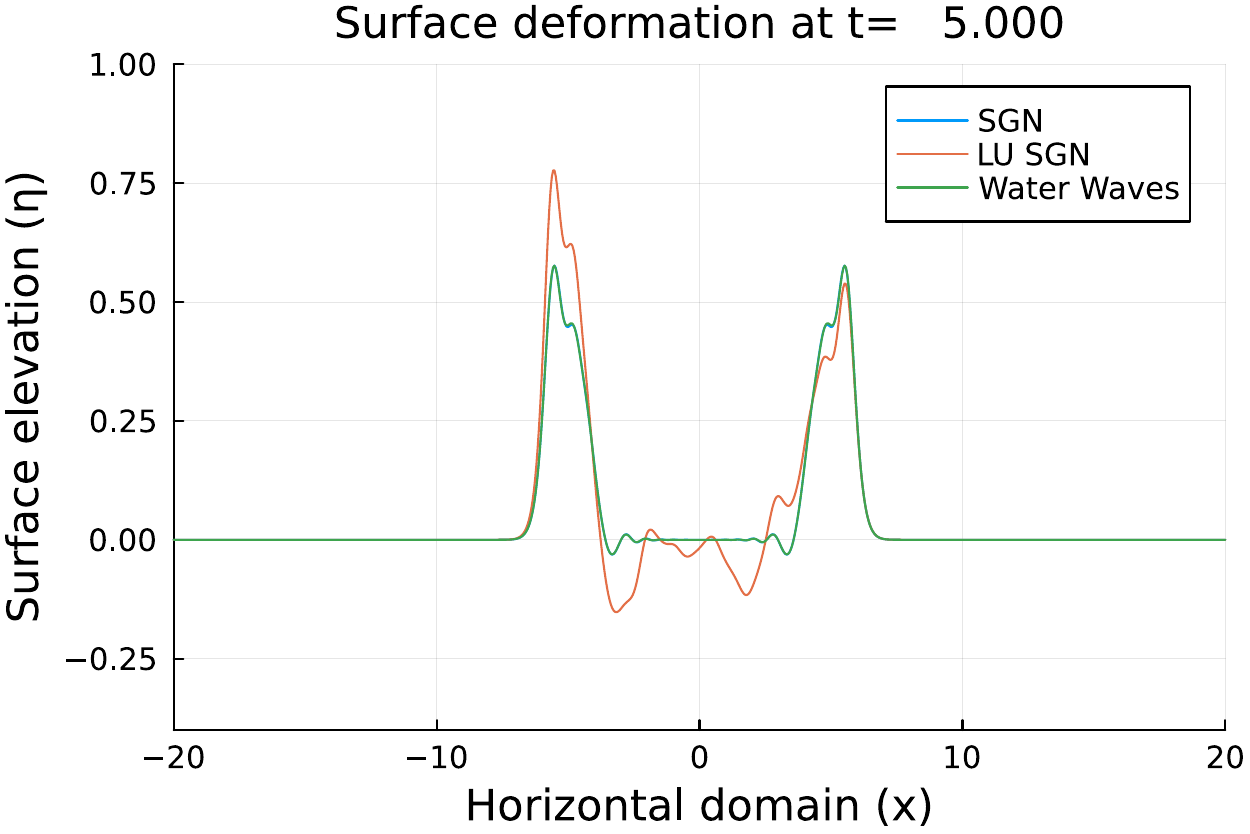}&
        \includegraphics[width=5cm]{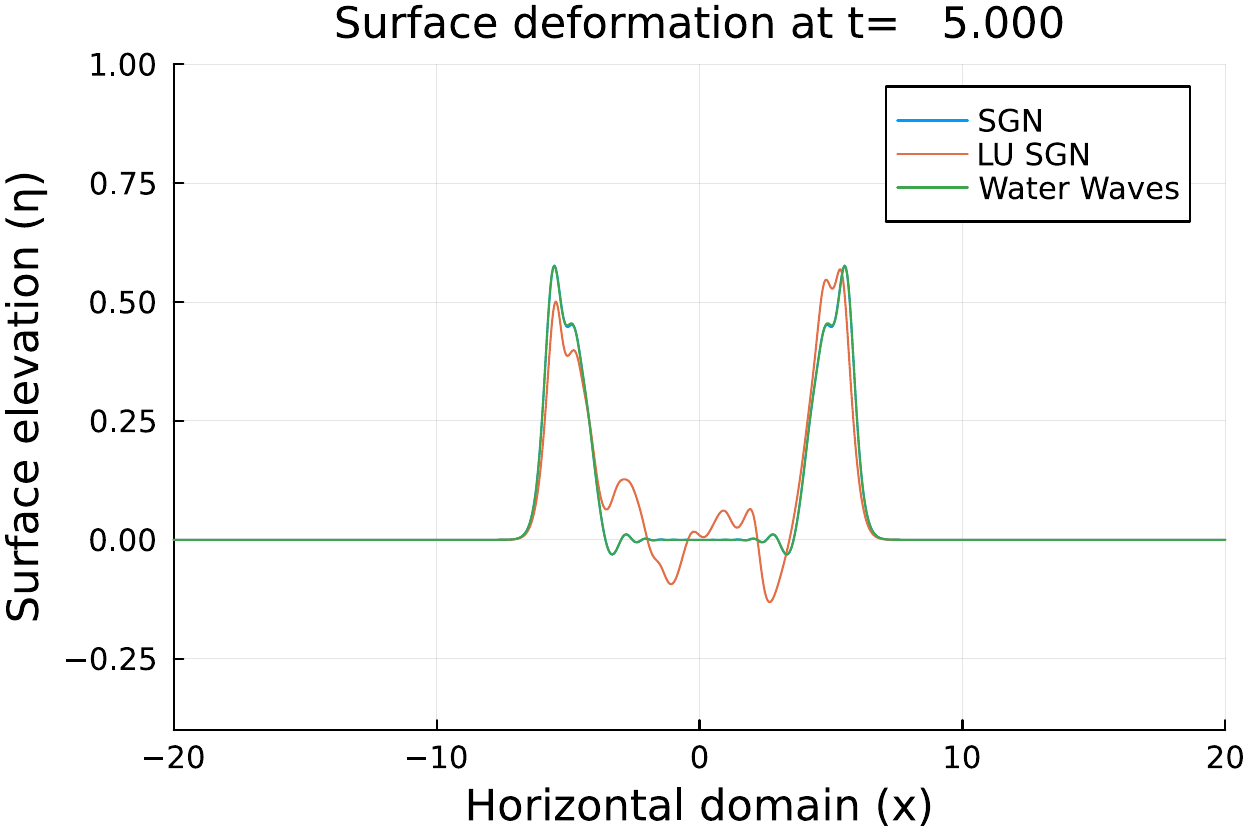}
    \end{tabular}
    \caption{\modif{Comparison between the surface deformation of the deterministic Saint-Venant (1st row), Boussinesq (2nd row) and Serre-Green-Naghdi (3rd row) models and their LU interpretations. Parameter set: $(\mathcal{P}_2)$ -- Wave number: $k = 2\pi/100$ -- Amplitude, from left to right: $0.001$, $0.005$ and $0.01$.}}\label{fig:plots-p2}
\end{figure}

The symmetry breaking mentioned in the previous subsection is also observed in \ref{fig:plots-p2}, where we use parameters $(\mathcal{P}_2)$. However, in this case, the parameter \modif{$\beta=0.1$} does not match the validity conditions of the Saint-Venant equations \modif{since $\beta \sim \epsilon$.} Therefore, the Water Waves solution is expected to \modif{differ} from the Saint-Venant ones, in both deterministic and LU forms. This is observed in figure \ref{fig:plots-p2}, indeed. In addition, the LU Boussinesq and LU Serre-Green-Naghdi appear to give an interesting variability to their deterministic versions. \modif{However, the  LU Serre-Green-Naghdi model yields numerical instabilities when choosing larger values of $\beta$ -- e.g. $\beta = 1$. This is due to a periodicity default caused by the symmetry breaking of the wave, which does not exist in the deterministic setting. To tackle this issue, we think of investigating numerical models in conservative form -- that is computing the water elevation/momentum variables $(\eta, h\overline{u})$ rather than the water elevation/momentum variables $(\eta,\overline{u})$. For similar reasons, our implementation does not enjoy conservation of physical quantities discussed earlier. Although it is relatively simple to translate our LU Saint-Venant algorithm into conservative form -- which then conserves mass, momentum and energy up to machine accuracy -- our LU Boussinesq and LU Serre-Green-Naghdi algorithms are more challenging to adapt due to the presence of the term $\overline{dG}$. We expect such numerical method to be more stable than the currently used one, and to allow investigating how multi-scale location uncertainty affects the waves dynamics. This is subject to further work.} 

\subsection{Numerical estimation of the noise-induced spreading}

In this section, we analyse the first and second order statistics of each LU wave model in the setting described above, at time $t=5s$. Again, we will study these models in the sets of parameters $(\mathcal{P}_1)$ for the LU Saint-Venant model, and $(\mathcal{P}_2)$ for the LU Boussinesq and Serre-Green-Naghdi models. The LU models statistics we analyse are computed with 130 realisations of each stochastic model. 

\begin{figure}[h]
    \centering
    \begin{tabular}{c|c|c}
        \includegraphics[width=5.25cm]{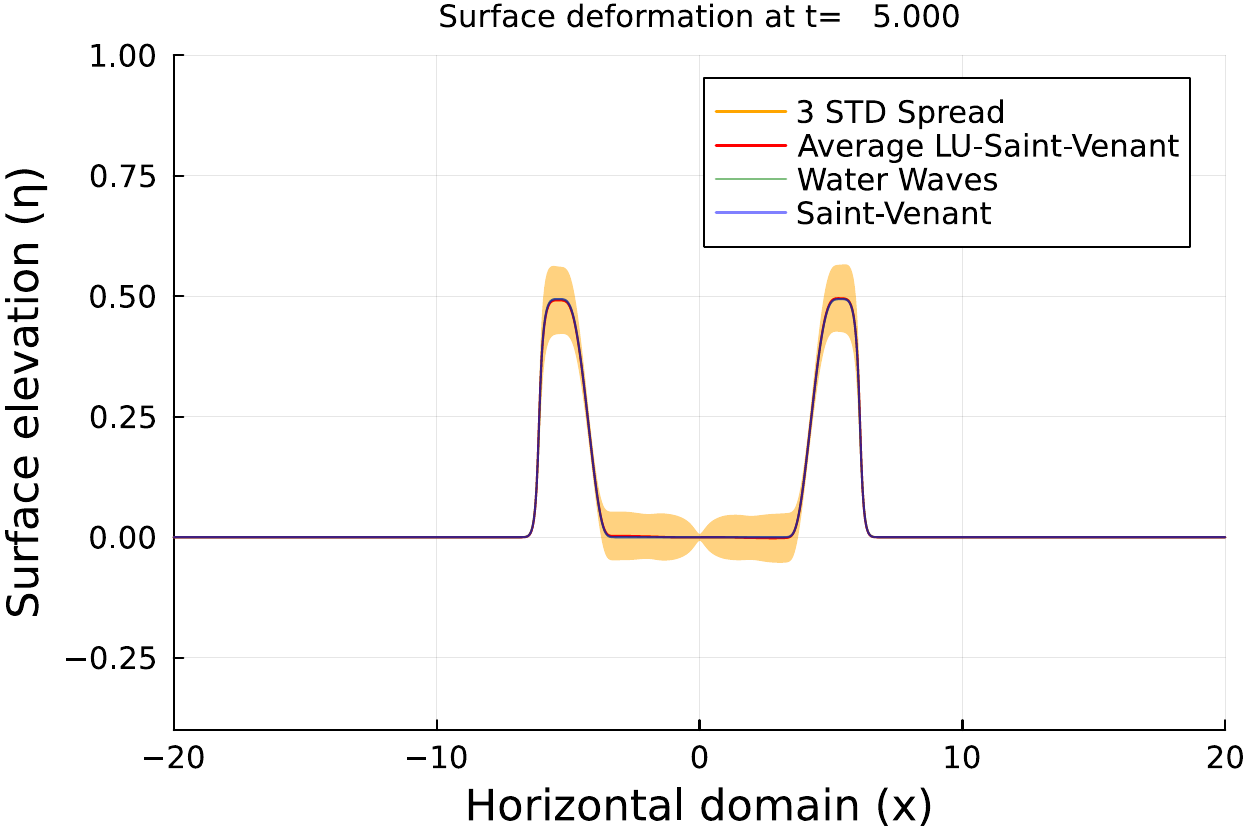} &
        \includegraphics[width=5.25cm]{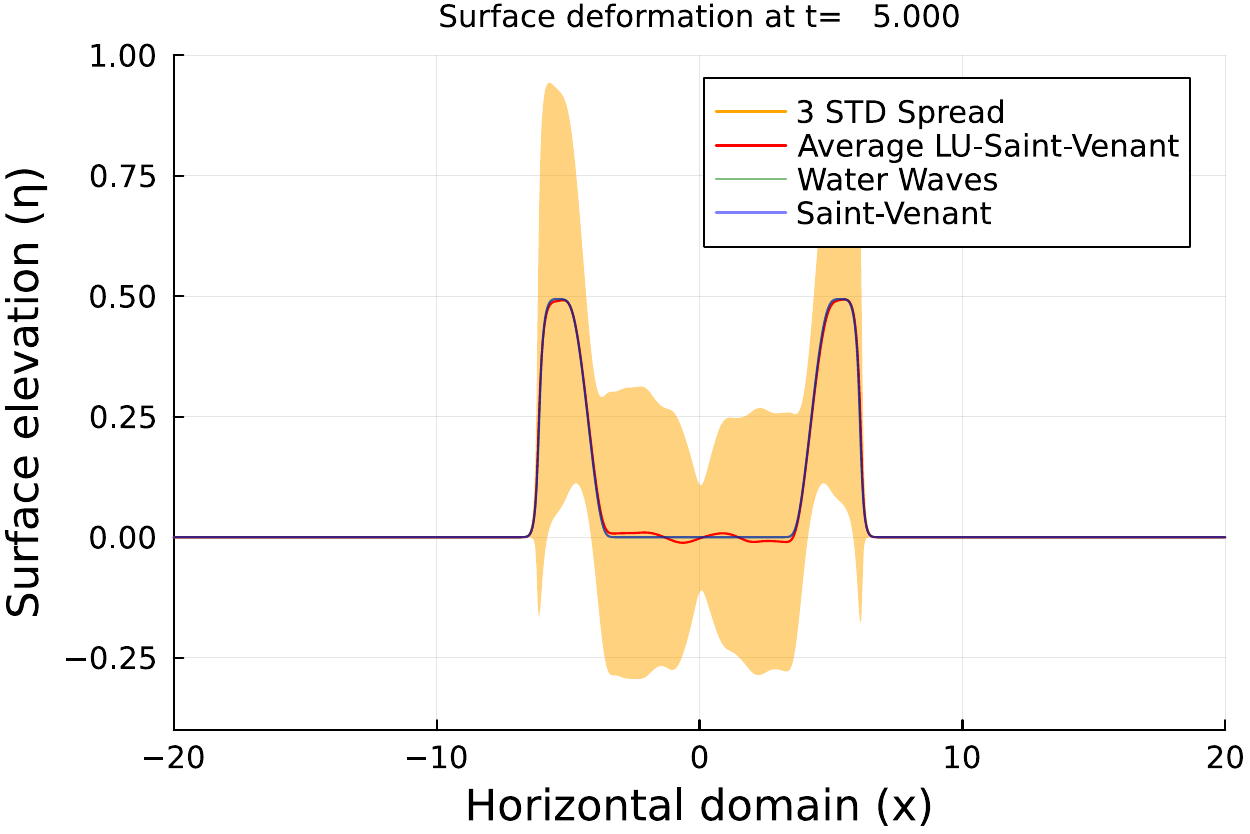} &
        \includegraphics[width=5.25cm]{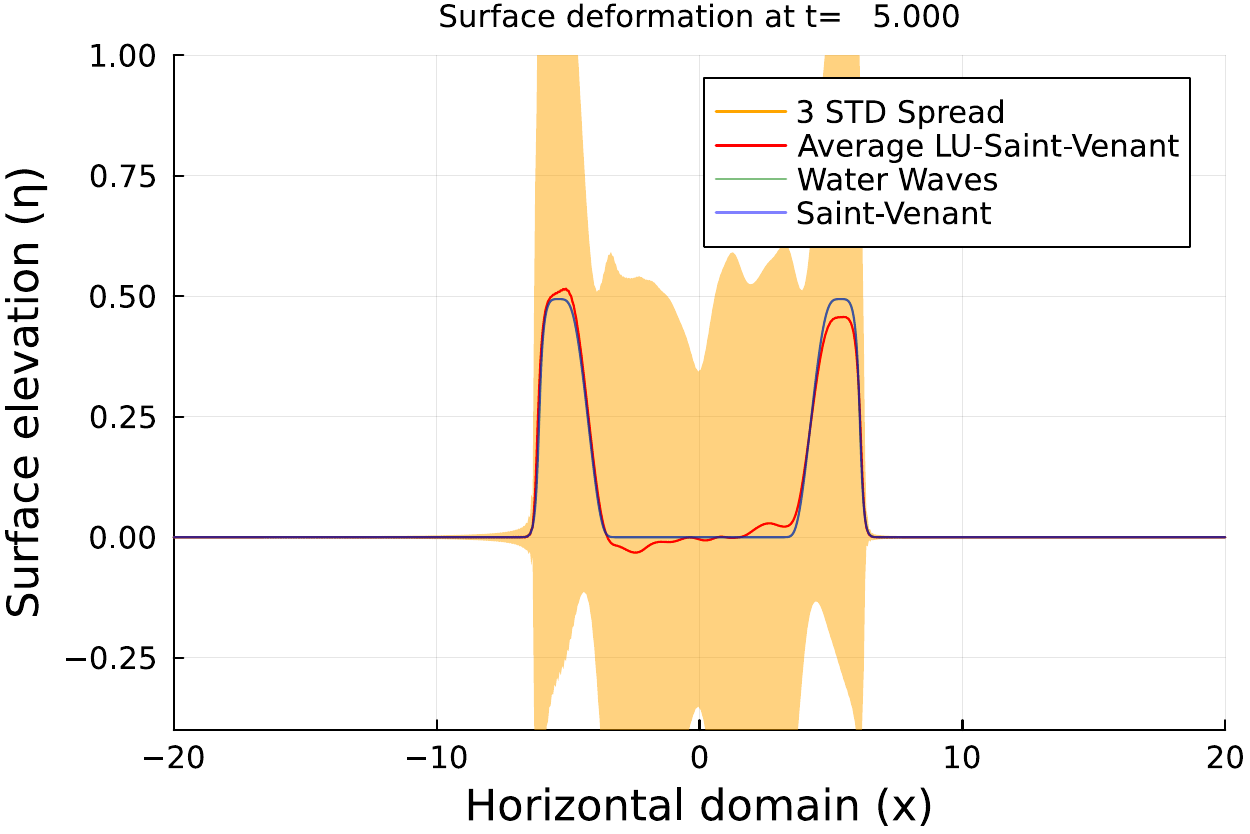}
    \end{tabular}
    \caption{\modif{First and second order statistics of the LU Saint-Venant model at $t=5s$, compared to the associated deterministic model and the water waves one. In addition, an evaluation of the spreading -- defined here as 3 times the empirical standard deviation -- is given for different values of noise amplitude (orange area). We only plot the solution over the domain $[-20,20]$. Parameter set: $(\mathcal{P}_1)$ -- Wave number: $k = 2\pi/100$ -- Amplitude, from left to right: $0.001$, $0.005$ and $0.01$. }} \label{fig:plots-stats-p1}
\end{figure}

A spreading arises from LU Saint-Venant model, which appears to grow linearly with the amplitude of the noise \modif{-- see Fig~\ref{fig:plots-stats-p1}}. Also, is concentrated in the ``upstream" of the wave, which is expected since it would be physically irrelevant for the wave to affect a region where it has not passed. In addition, \modif{the means and the standard deviations appear to be space symmetric -- up to statistical error -- which suggests the distribution is also space symmetric. Moreover, there exists a significant spreading at the peak of the wave for strong enough noise amplitudes.} 

\begin{figure}[h]
    \centering
    \begin{tabular}{c|c|c}
        \includegraphics[width=5.25cm]{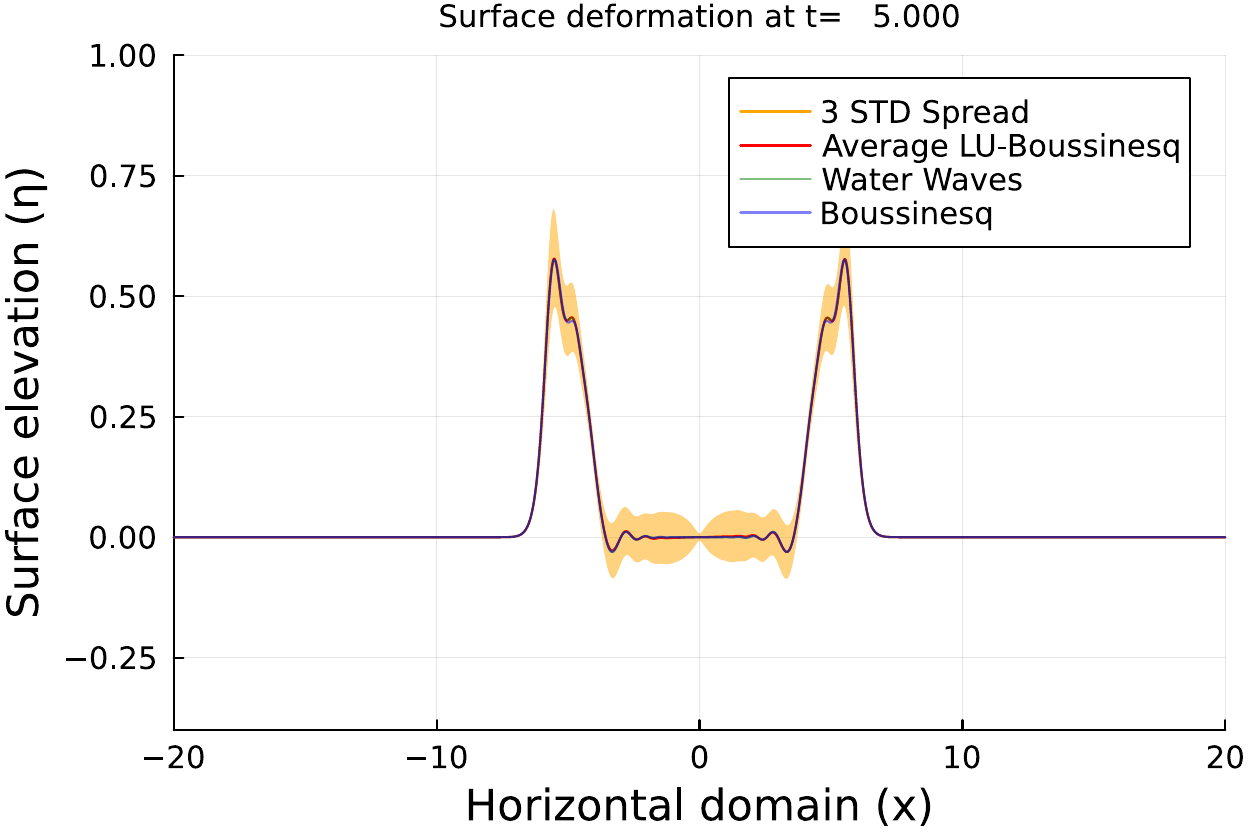} &     \includegraphics[width=5.25cm]{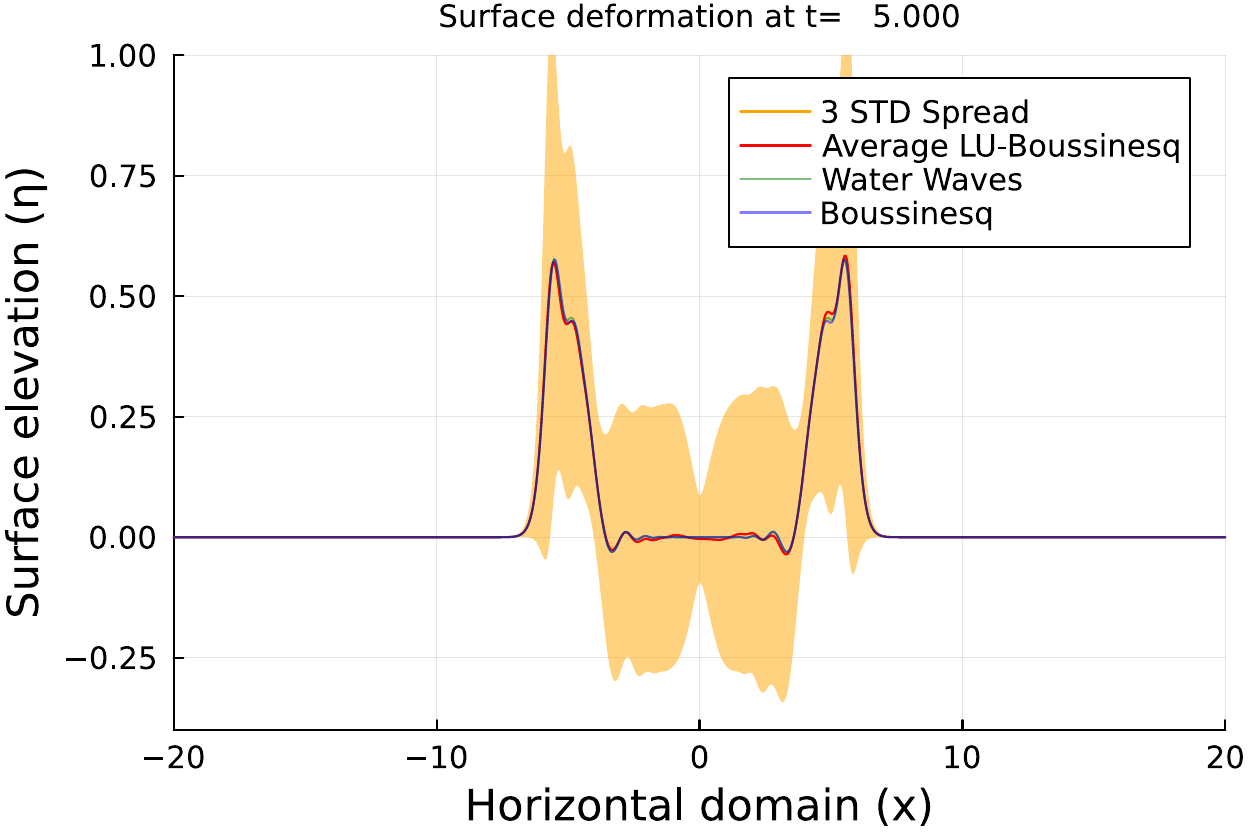} &    \includegraphics[width=5.25cm]{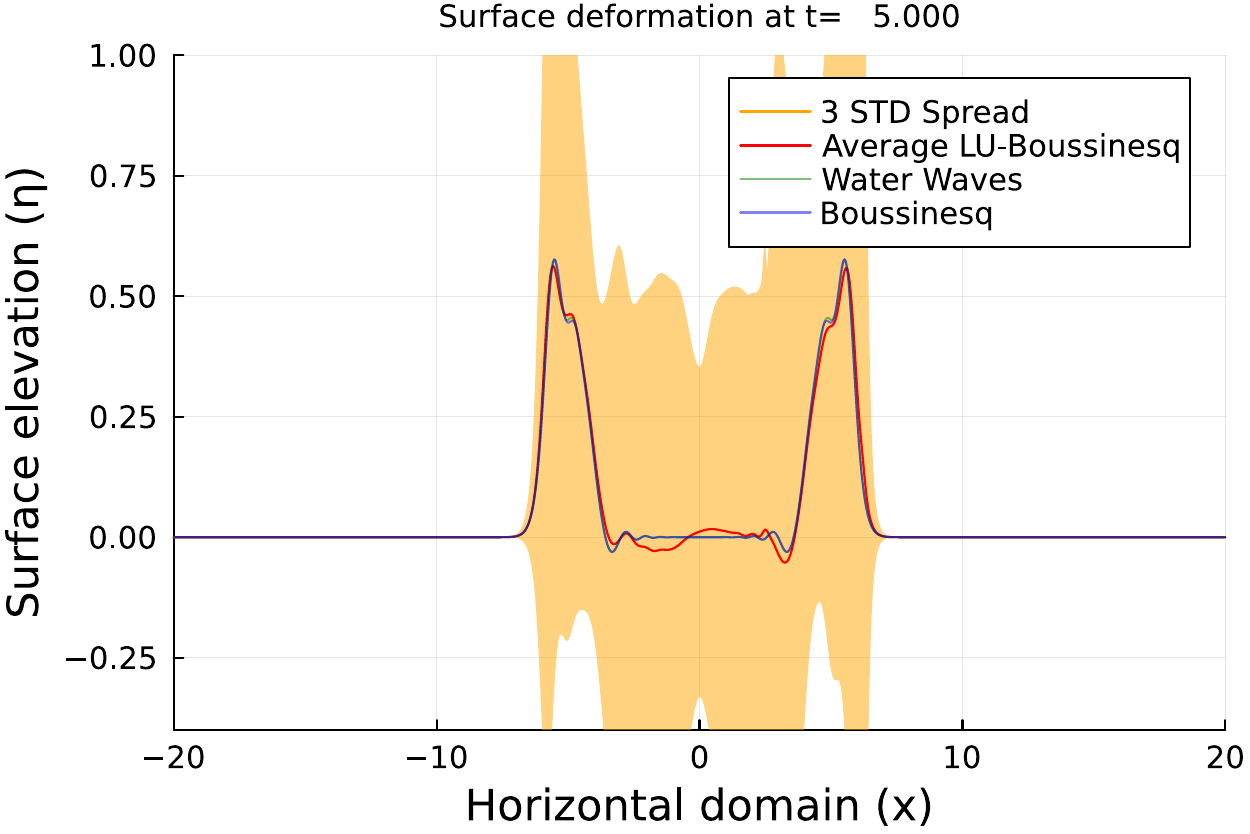}\\
        \hline
        \includegraphics[width=5.25cm]{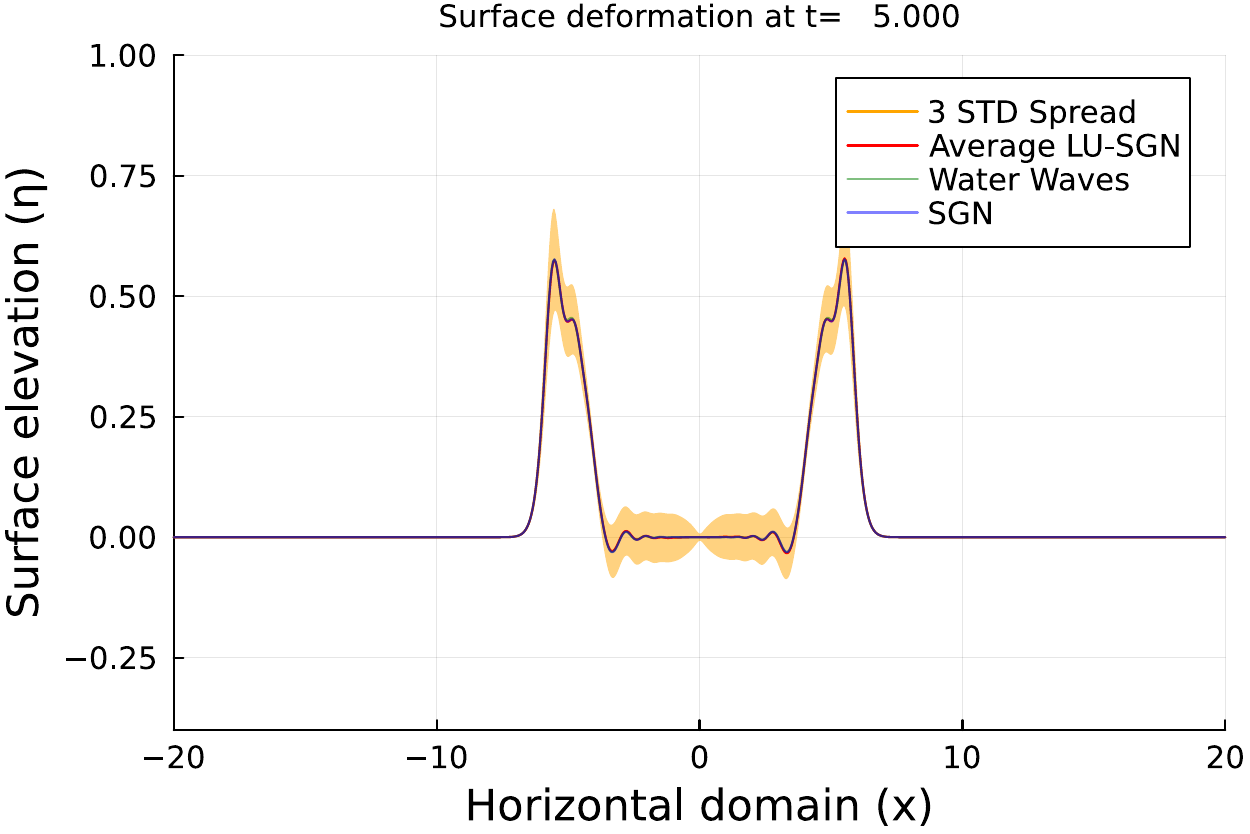}&
        \includegraphics[width=5.25cm]{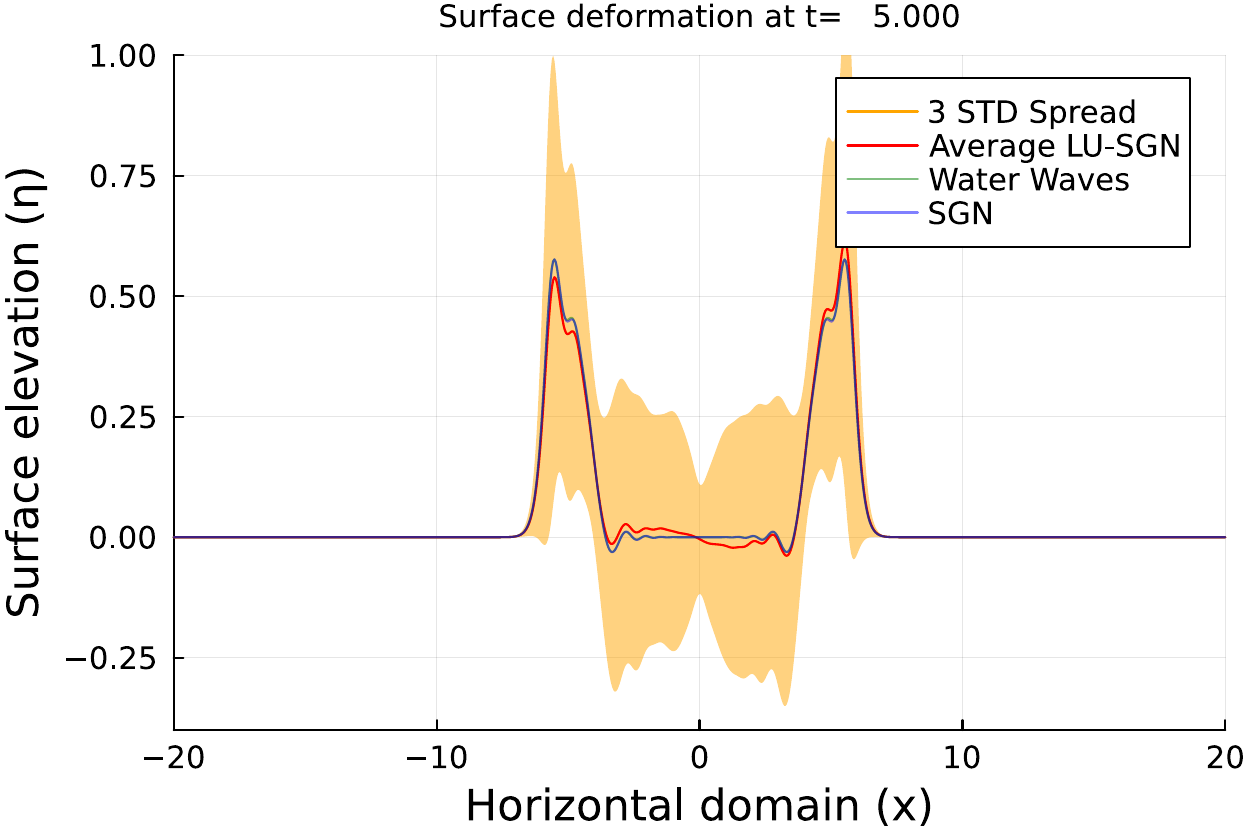}&
        \includegraphics[width=5.25cm]{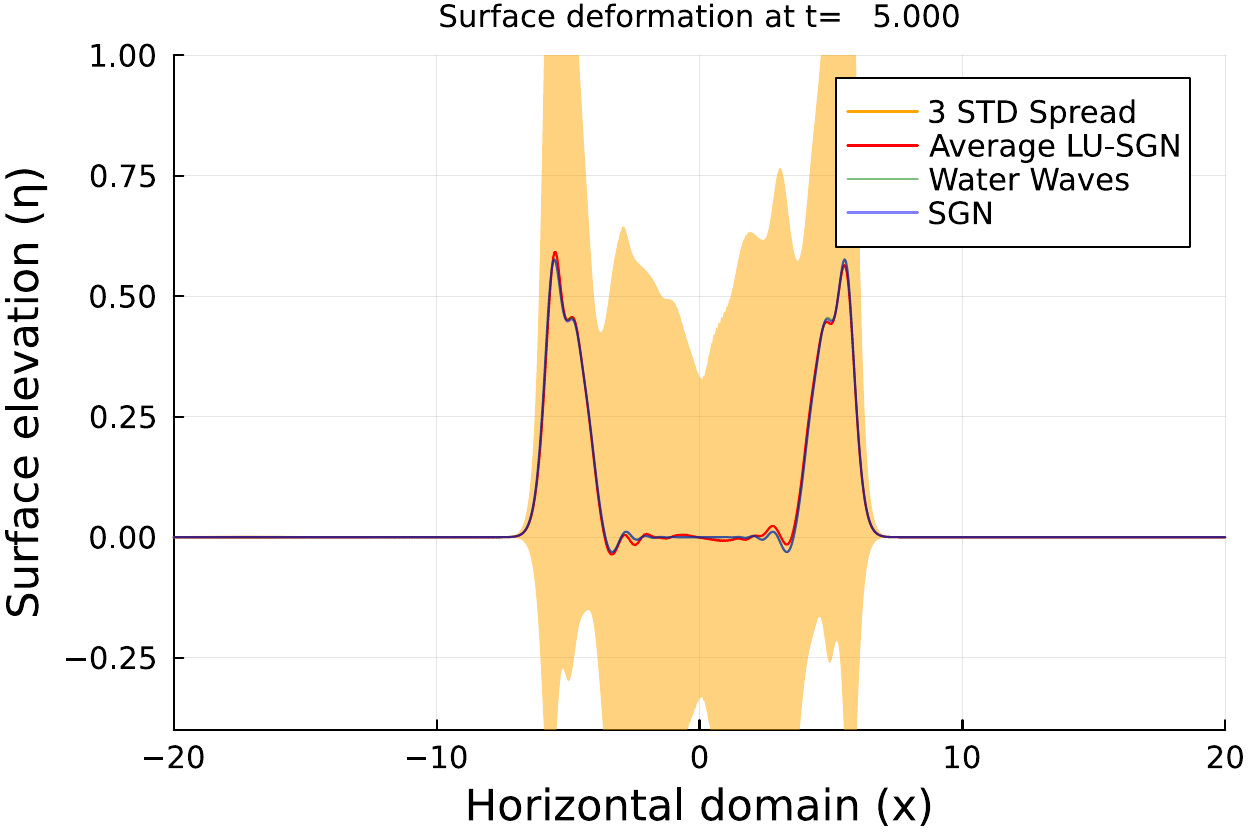}
    \end{tabular}
    
    \caption{\modif{First {and second} order statistics of the LU Boussinesq (1st row) and LU Serre-Green-Naghdi (2nd row) models at $t=5s$, compared to the associated deterministic models and the water waves one. In addition, an evaluation of the spreading -- defined here as 3 times the empirical standard deviation -- is given for different values of noise amplitude (orange area). We only plot the solution over the domain $[-20,20]$. Parameter set: $(\mathcal{P}_2)$ -- Wave number: $k = 2\pi/100$ -- Amplitude, from left to right: $0.001$, $0.005$ and $0.01$. }} \label{fig:plots-stats-p2}
\end{figure}

Regarding the LU Boussinesq and Serre-Green-Naghdi models, \modif{the same remarks on the peak spreading and the symmetry of the means and standard deviations apply} -- see Fig~\ref{fig:plots-stats-p2}. Again, the spreading is concentrated in the ``upstream" of the wave, specifically at its peaks and troughs. \modif{In particular, the maximum height value of the LU models varies depending on the stochastic realisation. Such observation gives lines of approach for building stochastic models by selecting the noise $\sigma$ to data, using for instance calibration or data assimilation techniques.}

\modif{Although the Serre-Green-Naghdi and the Boussinesq models give similar results for the parameters $(\mathcal{P}_2)$, we expect that much more differences would be observed between the two models with $\beta = 1$ and $\epsilon = 0.1$. As mentioned before, wave models written in conservative form should be more stable and allow to perform such tests in this configuration.

\section{Conclusion and discussion}

In this work, we investigated the stochastic representation of several shallow water coastal wave models within the LU framework. These stochastic models maintain the same conservation properties formally and thus exhibit physical consistency with their deterministic counterparts. We demonstrated numerically that they induce a pathwise symmetry breaking, accompanied by a restoration of this symmetry in law, which is a property that should be expected in representing turbulent effects. For low noise amplitudes, the studied models have  shown to converge toward the deterministic solutions. For large amplitudes, we need to transition to a numerical scheme in conservation form. This will be the subject of a future study.}



\appendix
\section{Appendix 1: Leibniz formula} \label{sec:Leibniz}
We give in the following the expression of the Leibniz integral rule for a stochastic process.  We want to evaluate the derivative of an integral of the form 
\[
    \int^{b(x,t)}_{a(x,t)} f(x,z,t) \dif z,
\]
with respect to the space,  $x$, and time, $t$, and where $a(x,z,t)$, $b(x,z,t)$ and $f(x,z,t)$ are continuous  $C^2$-semimartingale. For the space variable, no time derivative is at play and the usual formula stands:
\[
    \frac{\partial}{\partial x}\int^{b(x,t)}_{a(x,t)} f(x,z,t) \dif z = \int^{b(x,t)}_{a(x,t)} \!\!\!\!\!\!\!\partial_x f(x,z,t) \dif z +\partial_x b \,f(x,b(x,t),t) - \partial_x a \,f(x,a(x,t),t).
\]
Differentiation with respect to time, involves now a stochastic integration. Understanding the stochastic integral in the Stratonovich setting, we have,
\[
    \dif_t   \int^{b(x,t)}_{a(x,t)} f(x,z,t) \dif z = \int^{b(x,t)}_{a(x,t)} \!\!\!\!\!\!\! \dif_t   f(x,z,t) \dif z +\dif_t  b \,f(x,b(x,t),t) -  \dif_t   a \,f(x,b(x,t),t).
\]
The key argument of proof consists in defining the functions $F(x,y,t)=\int_{z_1}^{y} f(x,z,t)\dif z$ and $G(x,a,b,t) = \int_{a(x,t)}^{b(x,t)} f(x,z,t)\dif z$, then link them with a functional relation, that is
\begin{align*}
    G(x,a,b,t) &=  \int_{0}^{b(x,t)} f(x,z,t)\dif z - \int_{0}^{a(x,t)} f(x,z,t)\dif z\\
    & = F\bigl(x,b(x,t),t\bigr) - F\bigl(x,a(x,t),t\bigr).
\end{align*}
Now we introduce the following chain rule for Stratonovich calculus,   
\begin{theorem}[Generalized It\={o}'s formula -- Stratonovich form]\label{th:stratwentzell}
    Let $\theta(\xx,t), \xx\in\Omega$ be a continuous $C^3$-process and a continuous $C^2$-semimartingale, let $\XX_t$ be a continuous semimartingale with values in $\Omega$. Then, the following formula is satisfied :
    \begin{equation}\label{eq:dStheta}
    \dif \theta(\XX_t,t) = \dif_t \theta(\XX_t,t) + \frac{\partial\theta}{\partial x_i}(\XX_t,t) \circ \dif X_t^i.
    \end{equation}
\end{theorem} 

Upon applying this on $G(x,a,b,t)$, we obtain
\begin{align*}
    \dif_t  G &= \dif_t   F(x,b(x,t),t) + \partial_b F (x,b(x,t),t)\dif_t   b -  \dif_t F(x,a(x,t),t) - \partial_b F (x,b(x,t),t)\dif_t  a \\
    &=  \int^{b(x,t)}_{a(x,t)} \!\!\!\!\!\!\! \dif_t f(x,z,t) \dif z +\dif_t b \,f(x,b(x,t),t) -  \dif_t  a \,f(x,b(x,t),t),
\end{align*}
To interpret the stochastic integral in the It\={o} setting, we would need to adapt the previous chain rule as follows,
\begin{theorem}[Generalized It\={o}'s formula --     It\={o} form]\label{th:itowentzell}
    Let $\theta(\xx,t), \xx\in\Omega$ be a continuous $C^2$-process and a continuous $C^1$-semimartingale, let $\XX_t$ be a continuous semimartingale with values in $\Omega$. Then, $\theta(\XX_t,t)$ is a continuous semimartingale satisfying
    \begin{equation}\label{eq:dtheta}
        \dif \theta(\XX_t,t) = \dif_t \theta(\XX_t,t)+\frac{\partial\theta}{\partial x_i}(\XX_t,t)\dif X_t^i+\alf\frac{\partial^2\theta}{\partial x_i\partial x_j}(\XX_t,t)\dif\Big\langle X^i,X^j\Big\rangle_t+\dif\Big\langle\frac{\partial\theta}{\partial x_i}(\XX,\bcdot),X^i\Big\rangle_t.
    \end{equation}
\end{theorem}

Applying it to the process $G(x,a,b,t)$, we obtain
\begin{multline}
    \dif_t  G =  \int^{b(x,t)}_{a(x,t)} \!\!\!\!\!\!\! \dif_t f(x,z,t) \dif z +\dif_t b \,f(x,b(x,t),t) -  \dif_t  a \,f(x,b(x,t),t) + \\
    \dif\langle  f(x,b(x,t),t) , b\rangle_t \dif t - \dif\langle f(x,a(x,t),t) , a\rangle_t \dif t +\\
    \alf \dif\langle b,b\rangle_t\partial^2_{b^2} F (x,b(x,t),t) - \alf \dif\langle a,a\rangle_t \partial^2_{a^2} F (x,a(x,t),t)\dif t .
\end{multline}
This latter includes several  additional quadratic variation terms, and is more cumbersome to use in formal developments. However, this expression is necessary to access to the mathematical expectation.

\section{Appendix 2: Quadratic covariation}\label{sec:bracket}

In stochastic calculus, the quadratic covariation (or cross-variance) of two processes $X$ and $Y$ is defined as
\begin{equation}\label{eq:defbracket}
    \langle X, Y \rangle_t = \lim_{n\rightarrow 0} \sum_{i=1}^{p_n} (X_{i}^{n} - X_{i-1}^{n})(Y_{i}^{n} - Y_{i-1}^{n}),
\end{equation}
where $0=t_0^n<t_1^n<\cdots<t_{p_n}^n=t$ is a partition of the interval $[0,t]$. The previous limit, if it exists, is defined in the sense of convergence in probability.

Let $X$ and $Y$ two continuous semimartingales, defined as $X_t = X_0 + A_t + M_t, Y_t = Y_0 + B_t + N_t$, $M, N$ being martingales and $A, B$ finite variation processes. In this context, their quadratic covariation \eqref{eq:defbracket} exists, and is given by 
\begin{equation}\label{eq:semibracket}
    \langle X, Y \rangle_t = \langle M, N \rangle_t.
\end{equation}
In particular, the quadratic variation of a standard Brownian motion $B$ (as a martingale) is $\langle B \rangle_t := \langle B , B\rangle_t = t$, by definition.

Quadratic covariations play an important role in stochastic calculus, as they arise in It\={o}'s lemma, which can be interpreted as a stochastic chain rule. In particular, these terms are involved in the It\={o} integration by parts formula,
\begin{equation}\label{eq:IPP}
    \dif_t (XY) = X \dif_t Y + Y \dif_t X + \dif \langle X, Y \rangle_t,
\end{equation}
and in It\={o}'s isometry, expressing the covariance of two It\={o} integrals: let $f$ and $g$ two predictable processes such that $\int_0^t f^2 \dif \langle M,M \rangle_s$ and $\int_0^t g^2 \dif \langle N,N \rangle_s$ are finite, then
\begin{equation}\label{eq:isometry}
    \Exp \Big[ \big(\int_0^t f \dif M_s\big) \big(\int_0^t g \dif N_s\big) \Big] = \Exp \Big[ \int_0^t fg \dif \langle M, N \rangle_s \Big].
\end{equation}

\section{Appendix 3: LU Kordeveg - De Vries equation}\label{sec:KdV}
To derive the LU interpretation of the Kordeveg - De Vries equation (KdV), we will adapt a standard procedure to the stochastic case \cite{Johnson-97}. For simplification purpose, we will assume homogeneity in the transverse direction and consider in the following a 1D version of the Boussinesq model. Assuming that the bottom is flat, the 1D Boussinesq model reads
\begin{subequations}\label{S-B-1D}
\begin{align}
    &D^\hor_t \eta = - h \bigl(\partial_x (\overline{u} -\alf \Upsilon\epsilon \overline{u}_s)\dif t + \Upsilon^{\salf}  \partial_x\overline{\sodbt^x}\bigr),\\
    &\dif_t \overline{u}+ \epsilon \bigl( \overline{u}- \alf \Upsilon\epsilon \overline{u}_s\bigr) \partial_x \overline{u}\dif t + \Upsilon^{\salf} \epsilon \overline{\sodbt^x }\partial_x \overline{u}+ \partial_x \eta \dif t 
    -  \epsilon\beta^2 \bigl(\frac{h^2}{3}  \partial^2_{xx} \dif_t  \overline{u}  \bigr) = {\cal O}(\beta^4,\epsilon \beta^4)
\end{align}
\end{subequations}
The solutions we consider for the above Boussinesq system are assumed to be waves, and we apply the following change of variable: $\xi = x -t + \varphi_t$, where $\varphi_t = f(B_t)$ is a random phase that does not depends on $x$. The wave shape is also assumed to change on a large temporal scale $\tau = \epsilon t$. The depth averaged horizontal velocity $u(\xi, \tau, \varphi )$ and the surface elevation $\eta (\xi, \tau, \varphi )$ are assumed to be smooth functions. The noise $\sodbt^x (\xi, \tau)$ is assumed to be homogeneous, and thus is associated to a constant variance tensor and a zero ISD. In this new formalism the surface elevation equation reads
\begin{equation}
- \partial_\xi \eta  \, \dif t + \epsilon \partial_\tau \eta  \, \dif t  + \partial_\varphi \eta f'(B_t) \dif B_t  + \epsilon u \partial_\xi \eta\, \dif t + \Upsilon^{\salf}\epsilon\sodbt^x \partial_\xi \eta  +h \bigl(\partial_\xi u  \, \dif t +  \Upsilon^{\salf} \partial_\xi \sodbt^x\bigr) =0.
\end{equation}
After converting this equation in It\={o} form, terms of finite variations (ie ``$\dif t$" terms)  and martingale terms (i.e. ``$\dif B_t$" terms)   can be rigorously separated by the Biechteller-Delacherie theorem. For the martingale terms we have
\begin{equation}\label{M1}
    \partial_\varphi \eta \, f'(B_t) \dif B_t =  - \Upsilon^{\salf} \partial_\xi ( h \sdbt^x),
\end{equation}
and for the finite variation terms
\begin{equation}
    -\partial_\xi \eta    + \epsilon \partial_\tau \eta + \epsilon u \partial_\xi \eta    - \alf\Upsilon\epsilon^2 a^{\hor}\partial_{\xi\xi} \eta   + h  \partial_\xi u  =0.
\end{equation}
 
Moreover,for the velocity equation we have 
\begin{equation}
\begin{split}
    &- \partial_\xi u   \, \dif t + \epsilon  \partial_\tau u   \, \dif t  + \partial_\varphi u f'(B_t) \dif B_t  + \epsilon u  \partial_\xi u  \, \dif t + \Upsilon^{\salf} \epsilon \sodbt^x \partial_\xi  u \\ &\hspace*{2cm} -\alf \Upsilon \epsilon^2 \bar a ^{hh} \partial_{\xi \xi}^2 u \, \dif t + \partial_\xi \eta  \, \dif t
    - \frac{1}{h} \epsilon\beta^2 \bigl(\frac{h^3}{3}  \partial^2_{\xi\xi}(- \partial_\xi u+ \epsilon  \partial_\tau u)  \bigr) \, \dif t= {\cal O}(\beta^4,\epsilon \beta^4).
\end{split}
\end{equation}
Then, the martingale terms yield
\begin{equation}\label{M2}
    \partial_\varphi u f'(B_t) \dif B_t = - \Upsilon^{\salf} \epsilon \sodbt^x \partial_\xi  u 
\end{equation}
and for the finite variation terms 
\begin{equation}
    - \partial_\xi u   + \epsilon \partial_\tau u + \epsilon u\partial_\xi u   - \alf\Upsilon\epsilon^2 a^{\hor}\partial_{\xi\xi} u   + \partial_\xi \eta  + \frac{1}{h} \epsilon\beta^2 \bigl(\frac{h^3}{3}(  \partial^3_{\xi\xi\xi} u - \epsilon \partial^3_{\xi\xi\tau} u  \bigr)   =0.
\end{equation}
 
Expanding $u$ and $\eta$ in terms of the small parameter $\epsilon$ as $q= q_0 + \epsilon q_1 + \cdots$, and identifying the equations term by term of corresponding order, we get for the zero order terms,
\[
    \partial_\xi \eta_0 = \partial_\xi u_0  \text{ and hence } \eta_0= u_0.
\]
At order $\epsilon$, for a  noise of magnitude up to $\Upsilon \sim {\cal O}(1)$, we obtain the system 
\begin{align*}
    & - \partial_\xi u_1 +\partial_\tau \eta_0 +\eta_0 \partial_\xi \eta_0 + \partial_\xi \eta_1 +\frac{1}{3} \partial^3_{\xi\xi\xi} \eta_0 =0,\\
    & - \partial_\xi \eta_1 +\partial_\tau \eta_0 + 2 \eta_0 \partial_{\xi} \eta_0 + \partial_\xi u_1 =0,
\end{align*}
from which one gets immediately the classical KdV equation with random phase
\begin{equation}
    \partial_\tau \eta_0  +\frac{3}{2} \eta_0 \partial_\xi \eta_0 + \frac{1}{6} \partial^3_{\xi\xi\xi} \eta_0=0.
\end{equation}
This equation has the structure of a Burger equation with an additional dispersive term. A modified KdV equation is obtained by considering a stronger noise of amplitude $\Upsilon\sim 1/\epsilon$. In that case the second order terms must be kept in the Taylor development of terms at $\epsilon$ order. Eventually, one obtains the dissipative KdV equation (with random phase)
\begin{equation}
    \partial_\tau \eta_0  +\frac{3}{2} \eta_0 \partial_\xi \eta_0 -\alf a^{\hor} \partial_{\xi\xi}\eta_0+ \frac{1}{6} \partial^3_{\xi\xi\xi} \eta_0 =0.
\end{equation}
  
The random phase of both KdV equation  is determined by \eqref{M1} and \eqref{M2}. From the latter, we notice immediately that $f'(B)$ and $\epsilon$ must share the same order for the solution not to be trivial. At order $\epsilon $ from \eqref{M2}, one has 
\[
    \partial_\varphi \eta_0 \, f'(B_t) \dif B_t =  - \Upsilon^{\salf} \partial_\xi \eta_0 \sdbt^{0,x}.
\]
 A simple way to choose $f$ is to impose $f(B_t)= - k \Upsilon^{\salf}\epsilon k^\sigma B_t$ and $\sigma = k^\sigma$, where $k^\sigma$ is associated to the waves form
 \begin{equation}
 \eta_0 = h_0 e^{i (kx - \omega t - \epsilon \Upsilon^{\salf} k k^\sigma B_t)}.
\end{equation}
This solution corresponds to stochastic linear waves solutions as derived in \cite{Dinvay-Memin-PRSA-22}. Note that the derivations above have been done using the strong hypothesis of $u$ and $\eta$ being smooth functions of time -- more precisely of finite variation. Looking now for more general stochastic solutions, one considers the following variables
\begin{equation}
u(\xi, \tau), \quad \sodbt^\hor (\xi, \tau), \quad \eta (\xi, \tau),
\end{equation}
which are not differentiable with respect to $\tau$ (i.e. they are semi-martingale stochastic processes). Making the same assumptions and deriving the equations in Stratonovich form in the same way as previously, we obtain two coupled SPDE's,
\begin{equation}
    - \partial_\xi \eta \dif t  + \epsilon \dif _\tau \eta  + \epsilon u \partial_\xi \eta   \, \dif t + \Upsilon^{\salf}\epsilon \sodbt \partial_\xi \eta  + h \bigl(\partial_\xi u  \, \dif t +  \Upsilon^{\salf} \partial_\xi \sodbt^x\bigr) =0,
\end{equation}
and
\begin{equation}
\begin{split}
    &- \partial_\xi u\dif t   + \epsilon  \dif_\tau u  + \epsilon u  \partial_\xi u  \, \dif t + \Upsilon^{\salf} \epsilon \sodbt^x \partial_\xi  u   \\ &  \hspace*{5cm} + \partial_\xi \eta  \, \dif t
    - \frac{1}{h} \epsilon\beta^2 \bigl(\frac{h^3}{3}  \partial^2_{\xi\xi}(- \partial_\xi u+ \epsilon  \partial_\tau u)  \bigr) \, \dif t= {\cal O}(\beta^4,\epsilon \beta^4).
\end{split}
\end{equation}
At zeroth order the system reads 
\[
    \partial_\xi \eta_0 = \partial_\xi u_0  \text{, hence } \eta_0= u_0 \text{ together with } \partial_\xi \sodbt^{0,x}=0.
\]
At order $\epsilon$, one has
\begin{align*}
    & - \partial_\xi u_1 \dif t  +\dif_\tau \eta_0  +\eta_0 \partial_\xi \eta_0\dif t + \Upsilon^{\salf} \sodbt^{0,x}\partial_\xi \eta_0  + \partial_\xi \eta_1\dif t  +\frac{1}{3} \partial^3_{\xi\xi\xi} \eta_0 \dif t =0,\\
    & - \partial_\xi \eta_1\dif t  +\dif_\tau \eta_0  + \Upsilon^{\salf} \sodbt^{0,x}\partial_\xi \eta_0 + 2 \eta_0 \partial_{\xi} \eta_0 \dif t  + \partial_\xi u_1\dif_t +  \Upsilon^{\salf} \partial_\xi \sodbt^{1,x}  =0.
\end{align*}
Assuming that the noise terms $\sodbt^{1,x}$ does not depend on space, yet making no assumption on $\sodbt^{0,x}$, one obtains a stochastic KdV equation  with transport noise,
\begin{equation}
    \dif_\tau \eta_0  +\frac{3}{2} \eta_0 \partial_\xi \eta_0  \dif t  + \Upsilon^{\salf} \sodbt^{0,x} \partial_\xi \eta_0   + \frac{1}{6} \partial^3_{\xi\xi\xi} \eta_0 \dif t =0.
\end{equation}
Relaxing the spatially constant noise assumption on $\sodbt^{1,x}$ would add an additive stochastic forcing $\alf\Upsilon^{\salf} \partial_\xi \sodbt^{1,x}$. In such a case, we would get a stochastic KdV equation forced by an additive white noise of the form studied in \cite{deBouard-Debussche-98}, and for which existence and unicity of solution have been shown in the Sobolev space $H^1(\reel)$.

Without this additive forcing term, we face a simpler system that boils down  to the deterministic one. As a matter of fact, proceeding to the change of variable suggested by Wadati \cite{Wadati-83}, one has
\begin{equation}
    X = \xi -\frac{2}{3}\Upsilon^{\salf}\int_0^t  {\sodbs}^{\!\!\!0,x},
\end{equation}
with $\bsigma$ being constant over time and space. Now reparametrising $\eta(\xi,\tau)$ as $\eta'(X,t)$, we obtain the unperturbed KdV equation 
\begin{equation}
    \dif_t \eta'_0  +\frac{3}{2} \eta'_0 \partial_X \eta'_0  \dif t    + \frac{1}{6} \partial^3_{XXX} \eta'_0 \dif t =0,
\end{equation}
using the chain rules $\partial_\xi \eta = \partial_X\eta' \partial_\xi X=  \partial_X\eta'$ and $\df_\tau \eta= \partial_X \eta' \df_t X + \df_t \eta'$.
Such an equation admits a solitary travelling wave solution given by
\beq  
    \eta'_0(X,t) = A\, {\rm sech}^2\bigl(\sqrt{\frac{3}{4}A}(X- (1+ \alf A)t\bigr).
\eeq 
Note that considering a solution with a random phase would lead to the presence of a new additive noise term $ \partial_\varphi \eta_0 \, \epsilon f'(B_t) \dif B_t $. Plus, the noise $\sodbt^{0,x}$ may not be constant in space anymore.

\bibliographystyle{elsarticle-num}
\bibliography{glob.bib}

\begin{thebibliography}{10}
\expandafter\ifx\csname url\endcsname\relax
  \def\url#1{\texttt{#1}}\fi
\expandafter\ifx\csname urlprefix\endcsname\relax\def\urlprefix{URL }\fi
\expandafter\ifx\csname href\endcsname\relax
  \def\href#1#2{#2} \def\path#1{#1}\fi

\bibitem{LDLM_2023}
L.~Li, B.~Deremble, N.~Lahaye, E.~M{\'e}min, Stochastic data-driven parameterization of unresolved eddy effects in a baroclinic quasi-geostrophic model, Journal of Advances in Modeling Earth Systems 15~(2) (2023) e2022MS003297.

\bibitem{TML_2023}
F.~L. Tucciarone, E.~M{\'e}min, L.~Li, Primitive equations under location uncertainty: Analytical description and model development (2023).

\bibitem{MLLTC_2023}
E.~M{\'e}min, L.~Li, N.~Lahaye, G.~Tissot, B.~Chapron, Linear wave solutions of a stochastic shallow water model check for updates, in: Stochastic Transport in Upper Ocean Dynamics II: STUOD 2022 Workshop, London, UK, September 26--29, Vol.~11, Springer Nature, 2023, p. 223.

\bibitem{Dinvay-Memin-PRSA-22}
E.~Dinvay, E.~M{\'e}min, \href{https://royalsocietypublishing.org/doi/abs/10.1098/rspa.2022.0050}{Hamiltonian formulation of the stochastic surface wave problem}, Proceedings of the Royal Society A: Mathematical, Physical and Engineering Sciences 478~(2265) (2022) 20220050.
\newblock \href {https://doi.org/10.1098/rspa.2022.0050} {\path{doi:10.1098/rspa.2022.0050}}.
\newline\urlprefix\url{https://royalsocietypublishing.org/doi/abs/10.1098/rspa.2022.0050}

\bibitem{Memin14}
E.~M\'emin, Fluid flow dynamics under location uncertainty, Geophys. \& Astro. Fluid Dyn. 108~(2) (2014) 119--146.
\newblock \href {https://doi.org/10.1080/03091929.2013.836190} {\path{doi:10.1080/03091929.2013.836190}}.

\bibitem{Restrepo_2007}
J.~M. Restrepo, \href{https://journals.ametsoc.org/view/journals/phoc/37/7/jpo3099.1.xml}{Wave breaking dissipation in the wave-driven ocean circulation}, Journal of Physical Oceanography 37~(7) (2007) 1749 -- 1763.
\newblock \href {https://doi.org/10.1175/JPO3099.1} {\path{doi:10.1175/JPO3099.1}}.
\newline\urlprefix\url{https://journals.ametsoc.org/view/journals/phoc/37/7/jpo3099.1.xml}

\bibitem{RRmWB_2011}
J.~M. Restrepo, J.~M. Ramírez, J.~C. McWilliams, M.~Banner, \href{https://journals.ametsoc.org/view/journals/phoc/41/5/2010jpo4298.1.xml}{Multiscale momentum flux and diffusion due to whitecapping in wave–current interactions}, Journal of Physical Oceanography 41~(5) (2011) 837 -- 856.
\newblock \href {https://doi.org/10.1175/2010JPO4298.1} {\path{doi:10.1175/2010JPO4298.1}}.
\newline\urlprefix\url{https://journals.ametsoc.org/view/journals/phoc/41/5/2010jpo4298.1.xml}

\bibitem{BmW_2002}
P.~S. Berloff, J.~C. McWilliams, \href{https://journals.ametsoc.org/view/journals/phoc/32/3/1520-0485_2002_032_0797_mtiogp_2.0.co_2.xml}{Material transport in oceanic gyres. part ii: Hierarchy of stochastic models}, Journal of Physical Oceanography 32~(3) (2002) 797 -- 830.
\newblock \href {https://doi.org/10.1175/1520-0485(2002)032<0797:MTIOGP>2.0.CO;2} {\path{doi:10.1175/1520-0485(2002)032<0797:MTIOGP>2.0.CO;2}}.
\newline\urlprefix\url{https://journals.ametsoc.org/view/journals/phoc/32/3/1520-0485_2002_032_0797_mtiogp_2.0.co_2.xml}

\bibitem{Kraichnan59}
R.~Kraichnan, The structure of isotropic turbulence at very high reynolds numbers, J. of Fluids Mech. 5 (1959) 477--543.

\bibitem{Serre_1953}
F.~Serre, Contribution {\`a} l'{\'e}tude des {\'e}coulements permanents et variables dans les canaux, La Houille Blanche~(6) (1953) 830--872.

\bibitem{GN_1976}
A.~E. Green, P.~M. Naghdi, A derivation of equations for wave propagation in water of variable depth, Journal of Fluid Mechanics 78~(2) (1976) 237--246.

\bibitem{Berthelemy-04}
E.~Barth{\'e}lemy, \href{https://doi.org/10.1007/s10712-003-1281-7}{Nonlinear shallow water theories for coastal waves}, Surveys in Geophysics 25~(3) (2004) 315--337.
\newblock \href {https://doi.org/10.1007/s10712-003-1281-7} {\path{doi:10.1007/s10712-003-1281-7}}.
\newline\urlprefix\url{https://doi.org/10.1007/s10712-003-1281-7}

\bibitem{DaPrato}
G.~D. Prato, J.~Zabczyk, Stochastic equations in infinite dimensions, Cambridge University Press, 1992.

\bibitem{Li-et-al-2022}
L.~Li, E.~M{\'e}min, G.~Tissot, Stochastic parameterization with dynamic mode decomposition, in: STUOD Proceedings, Springer Verlag, 2022.

\bibitem{Schmid10}
P.~Schmid, Dynamic mode decomposition of numerical and experimental data, J. Fluid Mech. 656 (2010) 5--28.

\bibitem{Dufee-QJRMS-22}
B.~Duf{\'e}e, E.~M{\'e}min, D.~Crisan, \href{https://doi.org/10.1002/qj.4247}{Stochastic parametrization: An alternative to inflation in ensemble kalman filters}, Quarterly Journal of the Royal Meteorological Society 148~(744) (2022) 1075--1091.
\newblock \href {https://doi.org/https://doi.org/10.1002/qj.4247} {\path{doi:https://doi.org/10.1002/qj.4247}}.
\newline\urlprefix\url{https://doi.org/10.1002/qj.4247}

\bibitem{Bauer-et-al-JPO-20}
W.~Bauer, P.~Chandramouli, B.~Chapron, L.~Li, E.~M{\'e}min, \href{https://journals.ametsoc.org/view/journals/phoc/50/4/jpo-d-19-0164.1.xml}{Deciphering the role of small-scale inhomogeneity on geophysical flow structuration: A stochastic approach}, Journal of Physical Oceanography 50~(4) (01 Apr. 2020) 983 -- 1003.
\newblock \href {https://doi.org/10.1175/JPO-D-19-0164.1} {\path{doi:10.1175/JPO-D-19-0164.1}}.
\newline\urlprefix\url{https://journals.ametsoc.org/view/journals/phoc/50/4/jpo-d-19-0164.1.xml}

\bibitem{Kunita}
H.~Kunita, Stochastic flows and stochastic differential equations, Cambridge University Press, 1990.

\bibitem{Resseguier-GAFD-I-17}
V.~Resseguier, E.~M{\'e}min, B.~Chapron, {Geophysical flows under location uncertainty, Part I Random transport and general models}, Geophys. \& Astro. Fluid Dyn. 111~(3) (2017) 149--176.

\bibitem{Kadri-CF-17}
S.~Kadri~Harouna, E.~M{\'e}min, \href{https://hal.inria.fr/hal-01394780}{{Stochastic representation of the Reynolds transport theorem: revisiting large-scale modeling}}, {Computers \& Fluids} 156 (2017) 456--469.
\newblock \href {https://doi.org/10.1016/j.compfluid.2017.08.017} {\path{doi:10.1016/j.compfluid.2017.08.017}}.
\newline\urlprefix\url{https://hal.inria.fr/hal-01394780}

\bibitem{Smagorinsky63}
J.~Smagorinsky, General circulation experiments with the primitive equation: I. the basic experiment, Monthly Weather Review 91 (1963) 99--165.

\bibitem{Redi-82}
M.~H. Redi, Oceanic isopycnal mixing by coordinate rotation, Journal of Physical Oceanography 12~(10) (1982) 1154--1158.

\bibitem{Craik-Leibovich-76}
A.~Craik, S.~Leibovich, Rational model for langmuir circulations, J. Fluid Mech. 73 (1976) 401--426.

\bibitem{Leibovich80}
S.~Leibovich, On wave-current interaction theories of langmuir circulations, J. Fluid Mech. 99~(4) (1980) 715--724.

\bibitem{Resseguier2020arcme}
V.~Resseguier, L.~Li, G.~Jouan, P.~Derian, E.~M\'{e}min, B.~Chapron, New trends in ensemble forecast strategy: uncertainty quantification for coarse-grid computational fluid dynamics, Archives of Computational Methods in Engineering (2020) 1886--1784.

\bibitem{Barthelemy-04}
E.~Barth{\'e}lemy, \href{https://doi.org/10.1007/s10712-003-1281-7}{Nonlinear shallow water theories for coastal waves}, Surveys in Geophysics 25~(3) (2004) 315--337.
\newblock \href {https://doi.org/10.1007/s10712-003-1281-7} {\path{doi:10.1007/s10712-003-1281-7}}.
\newline\urlprefix\url{https://doi.org/10.1007/s10712-003-1281-7}

\bibitem{KG_2016}
R.~Kozlovsky, J.~Grobman, The blue garden: coastal infrastructure as ecologically enhanced wave-scapes, Landscape Research (12 2016).
\newblock \href {https://doi.org/10.1080/01426397.2016.1260702} {\path{doi:10.1080/01426397.2016.1260702}}.

\bibitem{Green-Naghdi-76}
A.~Green, P.~Naghdi, A derivation of equations for wave propagation in water of variable depth, Journal of Fluid Mechanics 78 (1976) 237--246.

\bibitem{Serre-53}
F.~Serre, Contribution \`{a} l'\'{e}tude des \'{e}coulements permanents et variables dans les canaux, Houille Blanche 8 (1953) 374--388.

\bibitem{Brecht-et-al-2021}
R.~Brecht, L.~Li, W.~Bauer, E.~M{\'e}min, Rotating shallow water flow under location uncertainty with a structure-preserving discretization, Journal of advances in modelling earth systems 13~(12) (2021).

\bibitem{Bauer2020jpo}
W.~Bauer, P.~Chandramouli, B.~Chapron, L.~Li, E.~M\'{e}min, Deciphering the role of small-scale inhomogeneity on geophysical flow structuration: a stochastic approach, Journal of Physical Oceanography 50~(4) (2020) 983--1003.

\bibitem{DJ_2019}
D.~G. Dritschel, M.~R. Jalali, On the regularity of the green–naghdi equations for a rotating shallow fluid layer, Journal of Fluid Mechanics 865 (2019) 100–136.
\newblock \href {https://doi.org/10.1017/jfm.2019.47} {\path{doi:10.1017/jfm.2019.47}}.

\bibitem{Vallis-17}
G.~Vallis, Atmospheric and Oceanic Fluid Dynamics, Cambridge University Press, 2017.

\bibitem{MS_1985}
J.~Miles, R.~Salmon, Weakly dispersive nonlinear gravity waves, Journal of Fluid Mechanics 157 (1985) 519--531.

\bibitem{website_DN_2024}
V.~Duchêne, P.~Navarro, \href{https://perso.univ-rennes1.fr/vincent.duchene/post/waterwaves1d/}{Waterwaves1d.jl} (Jun 2022).
\newline\urlprefix\url{https://perso.univ-rennes1.fr/vincent.duchene/post/waterwaves1d/}

\bibitem{XT_2016}
A.~Xiao, X.~Tang, High strong order stochastic {Runge-Kutta} methods for {Stratonovich} stochastic differential equations with scalar noise, Numerical Algorithms 72 (2016) 259--296.

\bibitem{GF_2019}
F.~Gugole, C.~L.~E. Franzke, \href{https://doi.org/10.1515/mcwf-2019-0004}{Numerical development and evaluation of an energy conserving conceptual stochastic climate model}, Mathematics of Climate and Weather Forecasting 5~(1) (2019) 45--64 [cited 2024-03-21].
\newblock \href {https://doi.org/doi:10.1515/mcwf-2019-0004} {\path{doi:doi:10.1515/mcwf-2019-0004}}.
\newline\urlprefix\url{https://doi.org/10.1515/mcwf-2019-0004}

\bibitem{PS_2008}
G.~Pavliotis, A.~Stuart, Multiscale methods: averaging and homogenization, Springer Science \& Business Media, 2008.

\bibitem{FBLM_2023}
C.~Fiorini, P.-M. Boulvard, L.~Li, E.~M{\'e}min, A two-step numerical scheme in time for surface quasi geostrophic equations under location uncertainty (2023).

\bibitem{Johnson-97}
R.~S. Johnson, A Modern Introduction to the Mathematical Theory of Water Waves, Cambridge Univ. Press, 1997.

\bibitem{deBouard-Debussche-98}
A.~de~Bouard, A.~Debussche, {On the Stochastic Korteweg de Vries Equation}, Journal of Functional Analysis 154 (1998) 215--251.

\bibitem{Wadati-83}
M.~Wadati, \href{https://doi.org/10.1143/JPSJ.52.2642}{Stochastic korteweg-de vries equation}, Journal of the Physical Society of Japan 52~(8) (1983) 2642--2648.
\newblock \href {http://arxiv.org/abs/https://doi.org/10.1143/JPSJ.52.2642} {\path{arXiv:https://doi.org/10.1143/JPSJ.52.2642}}, \href {https://doi.org/10.1143/JPSJ.52.2642} {\path{doi:10.1143/JPSJ.52.2642}}.
\newline\urlprefix\url{https://doi.org/10.1143/JPSJ.52.2642}

\end{thebibliography}

\end{document}